\documentclass[12pt]{article}

\usepackage{amsthm,amsmath,stmaryrd,bbm,hyperref,geometry,xcolor}
\usepackage{amssymb}
\usepackage[english]{babel}
 \usepackage[utf8]{inputenc}
\usepackage{graphicx}
\usepackage{amsfonts,amssymb}
\usepackage{verbatim}
\usepackage{enumitem}

\usepackage{caption}
\usepackage{subcaption}

\newcommand{\po}{\left(}
\newcommand{\pf}{\right)}
\newcommand{\co}{\left[}
\newcommand{\cf}{\right]}
\newcommand{\cco}{\llbracket}
\newcommand{\ccf}{\rrbracket}
\newcommand{\R}{\mathbb R}

\newcommand{\N}{\mathbb N} 
\newcommand{\dd}{\text{d}}
\newcommand{\na}{\nabla}

\newcommand{\param}{\mathfrak{p}}

\newcommand{\new}[1]{#1}
\newcommand{\rev}[1]{#1}
\newcommand{\rpi}{\rev{\bar{\pi}}}

\newcommand{\modifplb}[1]{#1}

\newtheorem{theorem}{Theorem}

\newtheorem{proposition}[theorem]{Proposition}
\newtheorem{remark}[theorem]{Remark}
\newtheorem{lemma}[theorem]{Lemma}
\newtheorem{corollary}[theorem]{Corollary}

\newtheorem{assu}{Assumption}
\title{HMC and \new{underdamped} Langevin united in the unadjusted convex \new{smooth} case.}
\author{\new{Nicolaï Gouraud, Pierre Le Bris, Adrien Majka and} Pierre Monmarché}

\begin{document}

\maketitle

\begin{abstract}
We consider a family of unadjusted generalized HMC samplers, which includes standard position HMC samplers and discretizations of the underdamped Langevin process. A detailed analysis and optimization of the parameters is conducted in the Gaussian case, which shows an improvement from $1/\kappa$ to $1/\sqrt{\kappa}$ for the convergence rate in terms of the condition number $\kappa$ by using partial  velocity refreshment, with respect to classical full refreshments. A similar effect is observed empirically for two related algorithms, namely Metropolis-adjusted gHMC and kinetic piecewise-deterministic  Markov processes. Then, a stochastic gradient version of the samplers is considered, for which dimension-free convergence rates are established for log-concave smooth targets \new{over a large range of parameters}, gathering  in a unified framework previous results on position HMC and underdamped Langevin and extending them to HMC with inertia.    
\end{abstract}


\medskip

\noindent\textbf{MSC class:} 65C05

\section{Overview}\label{sec:overview}

Let $\rpi$ be a probability law with density proportional to $\exp(-U)$ with $U\in\mathcal{C}^2(\R^d)$. 
We are interested in the question of sampling $\rpi$ with Markov Chain Monte Carlo (MCMC) samplers in a class of unadjusted \new{generalized} Hamiltonian Monte Carlo (\new{g}HMC) chains. These are kinetic chains (i.e. with a state  described by a position $x$ and a velocity $v$)  which alternates between two steps: first, a position Verlet (or leapfrog) numerical integration of the Hamiltonian dynamics $(\dot x,\dot v)=(v,-\na U(x))$, given by  
\begin{equation}\label{Verlet}
\new{K\text{ times}}\left\{
\begin{array}{rcl}
x  & \leftarrow  & x  + \frac\delta2 v \\
v  & \leftarrow & v  - \delta \na U(x )\\
x  & \leftarrow & x  + \frac\delta2 v 
\end{array}\right.
\end{equation}
 (i.e. a step of acceleration performed between two half-steps of free transport) for some time-step $\delta>0$ \new{and  number of iterations $K\in\N$}; second, an auto-regressive Gaussian \new{randomization} of the velocity
\begin{equation}\label{eq:damping}
 v  \leftarrow  \eta v +  \sqrt{1-\eta^2}  G\,,
\end{equation}
where $G$ is a standard $d$-dimensional Gaussian variable (mean $0$, variance $I_d$) and $\eta\in[0,1)$ is the damping (or friction) parameter. \new{We say that \eqref{eq:damping} is a full (resp. partial) refreshment of the velocity if $\eta=0$ (resp. $\eta>0$).}

 Within this class, a chain is thus characterized by its parameters $\delta,K,\eta$. In the classical HMC scheme, $\eta=0$ and $K$ is of order $1/\delta$, so that the process converges as $\delta$ vanishes to an idealized HMC process, which follows the true Hamiltonian dynamics and sees its velocity refreshed anew according to a Gaussian distribution at a fixed period $T>0$. In this classical case where $\eta=0$, the position $(x_n)_{n\in\N}$ alone is a Markov chain, and thus in the rest of this work we refer to this case as the \emph{position HMC} case. On the other hand, if $K=1$ and $1-\eta = \gamma \delta +o(\delta)$ for some $\gamma>0$, the chain is a splitting scheme of the \new{(underdamped)} Langevin process, which is the diffusion process solving
\begin{equation}\label{eq:EDScontinuousLangevin2}
\left\{\begin{array}{lll}
\dd x_t & = & v_t\dd t \\
\dd v_t & = & - \na U(x_t)\dd t - \gamma v_t \dd t + \sqrt{2\gamma}\dd B_t\,.
\end{array}\right.
\end{equation}
\new{Finally, if $K=1$ and $\eta=0$, then, up to a half step-size shift, this is essentially  the Euler-Maruyama scheme of the overdamped Langevin diffusion
\begin{equation}\label{eq:overdampedLangevin}
\dd x_t = -\na U(x_t) \dd t + \sqrt{2}\dd B_t ,
\end{equation}
except that the step-size is $\delta^2/2$.
 }
These samplers are very popular for applications either in statistics, machine learning or molecular dynamics (where \eqref{eq:EDScontinuousLangevin2} also has a physical motivation). By considering a unified framework, a fair comparison of the (unadjusted) Langevin and  HMC processes is possible. Here, \emph{unadjusted} refers to the fact that there is no Metropolis-Hastings step at the end of the Verlet integrator (which would ensure that the equilibrium is correct, but could slow down the convergence). \new{Note that, in this work, following the nomenclature of molecular dynamics and statistical physics, without any specification, \emph{Langevin process} always refer to the kinetic/underdamped process \eqref{eq:EDScontinuousLangevin2}, while \eqref{eq:overdampedLangevin}, which is not the main topic of this work, is referred to as the overdamped Langevin process.}

\medskip

\new{This work is organized as follows. In the rest of Section~\ref{sec:overview} we give an informal presentation of our results (Gaussian target distributions  in Section~\ref{subsec:overview_gaussien},   general convex smooth case with stochastic gradient in Section~\ref{subsec:overview_general}, numerical experiments in Section~\ref{subsec:overview_experiment}). Our theoretical results are stated in Section~\ref{sec:mainresults}, first in the Gaussian case in Section~\ref{sec:results:Gaussian} and second in the general one in Section~\ref{sec:results:Wasserstein}. Section~\ref{sec:numerique} is devoted to the numerical experiments. The proofs are given in the appendix, respectively  in Section~\ref{sec:Gauss:proof} for the Gaussian case and Section~\ref{sec:proofDimFreeCV} for the general one.}

\subsection{Gaussian target distributions}\label{subsec:overview_gaussien}

Our first main contribution is  a detailed analysis and optimi\new{z}ation in term of the parameters $K,\delta,\eta$, in the Gaussian case, namely when $U$ is a quadratic function. Here,   contrary to \cite{Arnold,MonmarcheGuillin,HHS1,Lelievre2012} in continuous time, we don't consider a fixed Gaussian distribution, but rather we address the question of the optimal sampler to tackle a given class of quadratic potentials. This is detailed in Section~\ref{sec:results:Gaussian} below. These results can be informally summarized as follows.

 Given some fixed $m,L>0$, consider the problem of sampling, with the same Markov chain, all $d$-dimensional Gaussian distribution satisfying $m\rev{I_d}\leqslant \nabla^2 U \leqslant L \rev{I_d}$ \new{(where $A\leqslant B$ for symmetric matrices $A,B$ means that $B-A$ is positive)} , with a \new{tolerance} $\varepsilon>0$ on the invariant measure in the $\mathcal W_2$-Wasserstein sense (see \eqref{eq:deftolerance}). The simulation cost is measured by the number of computations of $\nabla U$ \new{which, up to some logarithmic terms in $m,L,d,\varepsilon$ (see in particular Remark~\ref{rem:rate}), is proportional to $1/\rho$ where $\rho$ (given in \eqref{eq:defrho}) is the long-time Wasserstein convergence rate of the chain per gradient computation}. Define the relative tolerance $\varepsilon'=\varepsilon\sqrt{L/d}$, the rescaled time-step $\delta'=\delta\sqrt L$ and the condition number $\kappa=L/m$.
\begin{enumerate}
\item For small values of $\varepsilon'$, uniformly in $\kappa$, the optimal choice of parameters in term of Wasserstein convergence rate is  given by
\begin{equation}\label{eq:1optimchoix}
\delta' = \sqrt{8\varepsilon'}\,,\qquad K = \left\lfloor \frac{\pi}{\delta' \po 1+ 1/\sqrt\kappa\pf}\right\rfloor\,,\qquad \eta = \frac{1-\sin \po \pi/(1+\sqrt\kappa)\pf }{\cos\po \pi/(1+\sqrt\kappa)\pf }\,,
\end{equation}
which gives a convergence rate
\[\rho \underset{\varepsilon'\rightarrow 0}\simeq  \frac{\delta' \po 1+ 1/\sqrt\kappa\pf}{\pi} \ln \po \frac{\cos\po \pi/(1+\sqrt\kappa)\pf }{1-\sin \po \pi/(1+\sqrt\kappa)\pf }\pf \underset{\kappa\rightarrow +\infty }\simeq \frac{\delta'}{\sqrt\kappa}\,.\]
\item Moreover, the choice $\delta' = \sqrt{8\varepsilon'}$, $K=1$ and $\eta=1-\sqrt m \delta=1-\delta'/\sqrt\kappa$ gives a convergence rate $\delta'/\sqrt{\kappa}$, thus equivalent to the optimal choice when $\kappa\rightarrow +\infty$.
\item  Besides, the optimal choice to tune a position HMC (i.e. with $\eta=0$) is given by $\delta'$ and $K$ as in \eqref{eq:1optimchoix}, yielding a convergence rate equivalent to $\pi\delta'/\kappa$ as $\kappa\rightarrow +\infty$ (in accordance with \cite{ChenVempala}).
\item \new{In the overdamped Langevin case ($K=1$, $\eta=0$), it was known from previous work (e.g. \cite{pmlr-v75-wibisono18a}) that the optimal time-step is $\delta' = \sqrt{8\varepsilon'} $ and yields a convergence rate  $\rho = (\delta')^2/\kappa $. }
\item Finally, if $K\delta' \geqslant \pi$, there is no convergence \new{(i.e. $\rho=0$)}.
\end{enumerate}

In particular, we see that  a suitably tuned Langevin dynamics (with $K=1$) is competitive for badly-conditioned  problems with respect to the optimal \new{g}HMC, while being more robust with respect to the choice of parameters. Indeed, if we tune a \new{g}HMC  with the optimal choices \eqref{eq:1optimchoix} for some $m,L$ and then use it to sample a Gaussian process where $\na^2 U$ has an eigenvalue $L' \geqslant (\sqrt{L}+\sqrt{m})^2$ (which, for ill-conditioned problems, is just slightly above $L$), then $K\delta'\geqslant \pi$ and thus the convergence rate goes from optimal to zero. This does not happen in the Langevin case $K=1$, as can be checked from the results of Section~\ref{sec:results:Gaussian}, since in this case the optimal parameters are far from the condition $K\delta' \geqslant \pi$. The second interesting point enlighten\new{ed} by these results is that, when $\kappa$ is large, we see that the optimal choice of $\eta$ in \eqref{eq:1optimchoix} is close to $1$, while a classical position HMC chain (with $\eta=0$) is never competitive in the regime $\kappa\rightarrow +\infty$, being off the optimal convergence rate by a factor $\pi/\sqrt\kappa$ \new{(this is in particular the case of the unadjusted (overdamped) Langevin algorithm (ULA), for which furthermore the dependency in $\varepsilon'$ in $\rho$ is of order $\varepsilon'$ instead of $\sqrt{\varepsilon'}$); i.e., to put it in different terms, the convergence rate is improved from the order $\kappa^{-1}$ to the order  $\kappa^{-1/2}$ by taking $1-\eta$ of order $\kappa^{-1/2} K\delta' = \sqrt{m} K \delta $, with respect to the classical case $\eta=0$ (whenever $1\leqslant K \leqslant 1/\delta'$). \new{Hence, we recover, with a similar scaling of the damping parameter, restricted to Gaussian target measures, a speed-up with respect to classical samplers similar to the one observed recently in \cite{Cao2019OnE} in a more general convex case but only at the level of the $L^2$ norm for the continuous-time Langevin diffusion.} To our knowledge, our analysis in the Gaussian case gives the first result in that direction, with a clear criterion of efficiency involving the discretization error, for a practical generalist algorithm (in the sense that the algorithm is not designed \rev{or finely-tuned specifically} for Gaussian targets) \rev{involving only the computation of the gradient of $\ln \bar{\pi}$ (without requiring to solve a linear system of some optimization problem, as in implicit or preconditionned schemes as in  \cite{hodgkinson2021implicit,klatzer2023accelerated,Pidstrigach} and references within)}. Besides, after a first preprint version of the present work was released, a similar conclusion (i.e. a convergence rate of order $\kappa^{-1/2}$ for a practical algorithm in the Gaussian case) was obtained in \cite{rHMCGaussien}, but now for position HMC ($\eta=0$)  with randomized integration times $T=K\delta$ (also, in total variation rather than Wasserstein distance, but this distinction is not crucial since, following \cite{MonmarcheSplitting},  \cite{BouRabeeEberle} or \cite[Proposition 3]{MonmarcheHMCentropic}, a total variation convergence can be obtained from a Wasserstein convergence, with the same rate). An interesting point is that the optimal random integration time in \cite{rHMCGaussien} is of order $m^{-1/2}$, in accordance with results in continuous-time which advocates a refreshment rate of order $m^{1/2}$, see \cite{LuWang,Doucet}. This is consistent with our own result, in which we see that $K\delta/(1-\eta)$, which is the typical time needed for the process to forget its initial velocity (and thus is the typical distance that would be covered by each coordinate of the position in a flat potential during this time, since each velocity coordinate is of order $1$), is of order $m^{-1/2}$ for the optimal choices of parameters, which is the standard deviation of the lowest frequency modes of the target. Illustrating our theoretical analysis in the Gaussian case with numerical experiments in Section~\ref{sec:numerique}, we highlight the three following relevant time-scales related to the parameters $K,\delta,\eta$:
\begin{itemize}
\item The step size $\delta$ determines the numerical error (the condition $\delta \sqrt{L} <2$ is necessary for the stability of the scheme).
\item The integration time $K\delta$ determines whether there are periodic resonances in the system, which drastically impairs the convergence, particularly of high frequency modes (the sharp condition $\delta K \sqrt{L} < \pi $ ensures that there are no such resonances).
\item In the absence of periodic resonances, the key parameter for the sampling efficiency is the (coordinate-wise) mean free path $K\delta/(1-\eta)$ (which should be of order $m^{-1/2}$ for optimal sampling).
\end{itemize}
In \cite{rHMCGaussien}, the use of randomized integration times  destroy the deterministic periodicity issues, which allows for larger integration times, hence reaching a mean free path of order $m^{-1/2}$ with $\eta=0$. The two approaches are thus similar and can be easily combined, for instance by replacing the deterministic step size $\delta$ by a randomized step-size $\tilde \delta$ of the same order in order to destroy possible resonances at the frequency $1/\delta$, then taking $K=1$ and $\eta=1-\sqrt{m}\delta $, in order to benefit from the robustness of the Langevin scaling with respect to sub-optimal choice of parameters (see Section~\ref{sec:numerique}). This is also the spirit of the randomized midpoint method of \cite{NEURIPS2019_eb86d510} for undajusted underdamped Langevin. The study of randomized step-size is beyond the scope of the present article.} 

\medskip

\new{As a conclusion, the analysis for  Gaussian target distributions indicates that it should be more efficient in practice to replace the widely used position HMC samplers by kinetic processes where the velocity is only partially refreshed between two steps of Hamiltonian dynamics, i.e. with inertia.  This  motivates the study of gHMC with inertia beyond the Gaussian case, which is the topic of the next section.}


\subsection{The general \new{strongly} convex \new{smooth} case with stochastic gradient}\label{subsec:overview_general}

\new{Motivated by the results in the Gaussian case,} our second main contribution is thus  a dimension-free non-asymptotic Wasserstein convergence result for the chain in the general \new{strongly} convex smooth case, namely when $0<m\rev{I_d}\leqslant \na^2 U\leqslant L\rev{I_d}$ but $\na^2 U$ is not necessarily constant \rev{(to be clear, here we mean the long-time convergence of the chain towards its biased equilibrium; by contrast, the convergence of the bias to zero as the step-size goes to zero does depend on the dimension, as we discuss below)}. In practice, this is a restrictive framework which is not much general than the Gaussian case, but the main point is that it means the results are robust in the sense that they do not rely on explicit formulas as in the specific Gaussian case. It can also be thought for multi-modal targets as a way to quantify the sampling of a local mode \rev{(a positive lower bound on $\na^2 U$ holding on some balls around the local minimizers of $U$, where the chain spends most of its time between transitions from one mode to another)}. For these reasons, and because it is convenient to get explicit bounds (in particular in terms of the dimension) this framework has been used in the literature to compare different samplers \cite{BouRabeeSchuh,Dwivedi2,ChenVempala,Chatterji1,dalalyan1,DoucetHMC,durmus2019,Dwivedi,MangoubiSmith,SSZ}, etc. (of course the non-convex case has also drawn much interest, but this is not the topic of the present work so that we refer the interested reader to e.g. \cite{BouRabeeEberle,BouRabeeSchuh,MonmarcheHMCentropic} and references within). In fact, our result is a way to revisit similar results such as \cite{MonmarcheSplitting} (for $K=1$) or \cite{HMC,Holmes2014,Pidstrigach} (for $K=O(\delta^{-1}),\eta=0$), i.e. we show that the Langevin and HMC samplers can be treated in a unified framework, in a single computation. Besides, the method is the same as in previous works, namely it relies on the synchronous coupling of two chains. A slight difference with all previous works except \cite{MonmarcheSplitting,BouRabeeSchuh} is that, similarly to those two works, we establish the long-time convergence of the discrete-time Markov chain toward its equilibrium (independently from the question of the numerical bias on the equilibrium) instead of combining the long-time convergence of the continuous-time ideal process toward the true target with a numerical error analysis to get non-asymptotic efficiency bounds (but not convergence rates) for the discretised chain. As explained in \cite{MonmarcheSplitting}, this different approach is motivated by the theoretical discussion in \cite{QinHobert} on the difficulties of obtaining convergence rates for discretized Markov chains which scale correctly with the time-step. More importantly, that way, we avoid the question of distinguishing whether the limit continuous-time process as $\delta\rightarrow 0$ is the Langevin diffusion \eqref{eq:EDScontinuousLangevin2} or an ideal HMC process where the Hamiltonian dynamics is solved exactly.

  Moreover, in this part, we consider possibly a stochastic gradient variant of the algorithm, namely the (exact) acceleration step in the Verlet scheme \eqref{Verlet} is replaced by
\[v_1 = v_0 - \delta b(x_{1/2},\theta)\]
where $\theta$ is a random variable  with a law $\omega$ on a set $\Theta$ and   $b:\R^d\times \Theta\rightarrow \R^d$. We have in mind the case where $b(x,\theta)$ with $\theta\sim \omega$ is an unbiased estimator of $\na U(x)$ for all $x\in\R^d$, i.e.  $\mathbb E_\omega(b(x,\theta))=\na U(x)$. Assuming that the Jacobian matrix $\nabla_x b(\cdot,\theta)$, seen as a quadratic form, is bounded by $m$ and $L$ for all $\theta$, the extension of the result  to this stochastic case is straightforward (as e.g. in \cite{chaterrjiSto,Dalalyan_sto}). We emphasize that our main contribution here with respect to previous works is neither the method (which is standard) nor the generalization to the stochastic gradient case (which is straightforward), but the fact that we treat simultaneously the cases of $K$ ranging from $1$ to $O(\delta^{-1})$.

 \begin{remark}
 In the stochastic gradient case, if the position Verlet integrator \eqref{Verlet} is replaced by a velocity Verlet one, namely
 \begin{equation*}
 \new{
\begin{array}{rcl}
v  & \leftarrow & v  - \frac\delta2 \na U(x )\\
x  & \leftarrow  & x  + \delta v \\
v  & \leftarrow & v  - \frac\delta2\na U(x ),
\end{array}}
\end{equation*} 
 as in \cite{MonmarcheSplitting}, then $(z_n)_{n\in\N}=(x_n,v_n)_{n\in\N}$ is not a Markov chain anymore: indeed (take $K=1$ to fix ideas), with the velocity Verlet integrator, the force $\nabla U(x_n)$ is used both in the transition from $z_{n-1}$ to $z_n$ and from $z_n$ to $z_{n+1}$. With the stochastic gradient version, this force is replaced by $b(x_n,\theta_n)$, which means that $\theta_n$ intervenes in both transitions. As a consequence, $(z_n,\theta_n)_{n\in\N}$ is a Markov process, but  $(z_n)_{n\in\N}$ alone is not. The problem is the same if $K>1$, the random $\theta$ used in the last Verlet step of a transition has to be used in the first Verlet step of the next transition. Although we expect similar results, the analysis would thus be slightly messier (in particular, in order to study the equilibrium bias, one should consider the first marginal of the invariant measure of $(z_n,\theta_n)_{n\in\N}$ since the invariant measure of the non-Markov chain $(z_n)_{n\in\N}$ is not well defined). This is why, with respect to \cite{MonmarcheSplitting}, we switched to the position Verlet integrator.
 \end{remark}
 
 The results concerning the Wasserstein convergence of the stochastic gradient \new{g}HMC chains are detailed in Section~\ref{sec:results:Wasserstein}. \new{ To summarize them, applying together Theorem~\ref{thm:crude} (for the long-time convergence rate to the biased equilibrium), Proposition~\ref{prop:numerique} (for the discretization error) and Proposition~\ref{prop:erreur_sto} (for the error due to the stochastic gradient approximation), we get that  there exists $\overline{T},\overline{\gamma}>0$ which depend only on $m,L$ such that for all parameters $\delta,K,\eta$ satisfying
 \begin{equation}\label{eq:widerange}
  \delta K \leqslant \overline{T},\qquad 1-\eta \geqslant \overline{\gamma} K \delta\,,  
 \end{equation}
 there exist $C,a>0$ which depends only on $m,L$, the Lipschitz norm of $\na^2 U$ and $(1-\eta)/(K\delta)$ such that,  assuming that the stochastic gradient is given by a Monte Carlo estimation with $p$ independent realizations (see Assumption~\ref{assu:gradient_sto} below), 
 \begin{equation}\label{eq:summary}
\mathcal W_2\po \mathcal Law(x_n), \rpi\pf \leqslant C \po e^{-a n K\delta } \mathcal W_2\po \mathcal Law(x_0), \rpi\pf + \delta^2 d + \sqrt{\frac{\delta (d+\sigma^2)}{p}}\pf \,,  
 \end{equation}
where $\sigma^2$ is the variance of one stochastic gradient realization (at the minimizer of $U$) and $\mathcal W_2$ stands for the $L^2$ Wasserstein distance (if $\na^2 U$ is not Lipschitz the term $\delta^2 d$ is replaced by $\delta \sqrt{d}$). As already mentioned, the main point is that this estimate is thus uniform over a large range of parameters, as \eqref{eq:widerange} goes from the Langevin to the HMC scaling.  Moreover, in terms of $m,L$, the choice of parameters which maximize the rate $a$ that we get are such that $\eta \rightarrow 1$ as $m/L\rightarrow 0$, which is consistent with our conclusion in the Gaussian case. However, these are just upper bounds for the genuine convergence rate, so we cannot deduce the advantage of inertia as clearly as in the Gaussian case (see the discussion after Corollary~\ref{Cor:contraction} below for a more detailed discussion on the dependency in $m,L$ in our results on the convergence rate). 

The stochastic gradient error in \eqref{eq:summary} is similar to the one obtained in previous works on similar algorithms \cite{chaterrjiSto,Dalalyan_sto}, and  sharp (besides, it is new for underdamped Langevin splitting schemes). In the non-stochastic gradient case (i.e. $p=+\infty$ in \eqref{eq:summary}) and when $\na^2 U$ is Lipschitz,   considering $Kn$ the number of gradient computation necessary to get $\mathcal W_2(\mathcal Law(x_n),\rpi) \leqslant \varepsilon$ for some tolerance $\varepsilon>0$,  we get that this can be obtained with $Kn$ of order $\sqrt{d/\varepsilon}\ln(\sqrt{d}/\varepsilon)$ (omitting the dependency in $m,L$) uniformly over a wide range of parameters (from Langevin splitting to position HMC). In terms of $d$ and $\varepsilon$, we recover the same complexity (under the same conditions) that 
position HMC (in \cite[Theorem 1.5]{ChenVempala}, \cite[Table 1]{HMC} or \cite[Theorem 3.8]{BouRabeeSchuh}) or for second order discretizations of the Langevin process such as  \cite[Table 1]{MonmarcheSplitting} for splitting schemes or \cite[Theorem 2]{dalalyan1} for a stochastic Euler scheme (which require the computation of $\na^2 U$, and thus is numerically more expensive). 
 When we don't assume that $\na^2 U$ is Lipschitz, from Proposition~\ref{prop:numerique} below we can see that we only get $\sqrt{d}/\varepsilon\ln(\sqrt{d}/\varepsilon)$, in which case we recover, under the same conditions, the efficiency previously obtained for the overdamped Langevin process (see  \cite[Table 2]{Dwivedi}) or for first order discretization schemes of the Langevin process (as in \cite[Theorem 1]{dalalyan1} or \cite[Theorem 3]{Chatterji1}).  }
 
 \rev{Finally, let us mention that, after the release of the first preprint version of the present study, the works \cite{leimkuhler1,leimkuhler2} have  conducted an analysis much related to the question discussed in this section, with similar techniques but emphasizing different points. In particular, they concern only splitting schemes of the Langevin diffusion;  \cite{leimkuhler1} considers a wide variety of splitting schemes with a particular focus on the restriction on the step-size and the consistency in the high friction regime;  \cite{leimkuhler2} considers the case of stochastic gradients and works with conditions weaker than ours (they assume that the target (non-stochastic) potential is convex but not a uniform contraction property for the stochastic drift as we do in Assumption~\ref{Assu:mLb}). }

\subsection{Numerical analysis of kinetic samplers}\label{subsec:overview_experiment}

\new{Finally, our third contribution is a numerical analysis of the influence of the velocity \new{randomization} mechanism on the efficiency of the algorithms. This allows to consider two families of samplers which are related to unadjusted gHMC, but not covered by our theoretical results: first, Metropolis-adjusted gHMC (MagHMC) and, second,  kinetic piecewise deterministic Markov processes (PDMP) samplers (see e.g. \cite{MonmarcheRecuitPDMP,DoucetHMC,BierkensRobertsZitt,MMZ} and references within), which have been shown to converge in some regimes to Hamiltonian-based dynamics  \cite{DoucetHMC,MMZ}.  In both cases, the first marginal of the invariant measure of the chain is the correct target measure.

Concerning MagHMC, let us notice the following. In he classical case $\eta=0$, the Metropolis adjusted HMC chain is reversible with respect to $\rpi$. This is not the case when $\eta>0$. Indeed, a transition of MagHMC is given by three successive step: first, a Metropolis-Hastings step using as a proposal a Verlet trajectory followed by a reflection of the velocity. Second, a reflection of the velocity. Third, the Gaussian \new{randomization} of the velocity. Each of these steps is reversible with respect to $\mu = \rpi \otimes \mathrm{N}(0,I_d)$ (with $\mathrm{N}(0,I_d)$ the standard Gaussian measure), but the full transition is not (except when $\eta=0$, where moreover the velocity reflection does not appear since it is erased by the refreshment). Hence,  using partial refreshment for MagHMC instead of $\eta=0$ amounts to use a non-reversible chain instead of a reversible one, which goes in the direction of the numerous works on non-reversible sampling in the last decade \rev{\cite{power2019accelerated,diaconis2013spectral,10.1214/16-ECP25,lelievre2013optimal,bierkens2016non,duncan2017nonreversible,miclo2013etude,doi:10.1137/20M1378752} (preceded by the seminal \cite{Diaconisetal} among others)}.

 Concerning PDMP samplers, by replacing the Verlet integration of the Hamiltonian dynamics in an HMC scheme by a such a kinetic process, one gets a rejection-free non-reversible unbiased sampler, namely the target equilibrium is correct but, by contrast to standard Metropolis-adjusted schemes, rejections are replaced by velocity jumps, which should improve the sampling (since the process never stops moving, and is only reflected in the gradient direction). Up to now, most works on these samplers have considered a (usually randomized) full refreshment of the velocity, and the question of the scaling of the refreshment rate has been abundantly discussed. In view of this, it is natural to ask whether the highlights obtained in our theoretical study of the unadjusted gHMC sampler extrapolates to these processes.
 

The conclusion from this empirical study, presented in Section~\ref{sec:numerique}, is that the same phenomenon is indeed observed for the three samplers, concluding on the benefit of samplers with inertia. In particular, it comforts the role of the mean free path $K\delta/(1-\eta)$ as a key variable in the sampling efficiency for kinetic processes. Besides, we observe that, even with fixed deterministic integration times, the inherent randomness of PDMP samplers, which is not sufficient in general to ensure ergodicity, has however the effect of suppressing  periodicity issues.  }

\section{Results}\label{sec:mainresults}

Let us first introduce some notations and give a clear definition of the Markov chains which are the topic of this work. We use the notations of the introduction. For $\delta>0$, $\theta\in\Theta$ and $(x,v)\in\R^{2d}$, let
\begin{equation}\label{eq:VerletSto}
\Phi_\delta^\theta(x,v) = \po x+\delta v - \frac{\delta^2}{2} b\po x +\frac\delta2 v,\theta \pf\, , \, v - \delta  b\po x +\frac\delta2 v,\theta \pf\pf,
\end{equation}
so that a stochastic Verlet step (i.e. \eqref{Verlet} with $\na U$ replaced by $b(\cdot,\theta)$) reads $z_{k+1} = \Phi_\delta^\theta(z_k)$. When $b(\cdot,\theta)=\na U$ for all $\theta$, we simply write $\Phi_\delta^\theta = \Phi_\delta$. Consider the corresponding Markov transition operators $P_{sV}$ and $P_{V}$ given for bounded measurable functions $f$ by
\[P_{sV} f(x,v) = \mathbb E_{\omega} \po f\po \Phi_\delta^\theta(x,v)\pf\pf \,,\qquad P_V f(x,v) = f\po\Phi_\delta(x,v)\pf\,.\]
The Markov transition operator corresponding to the \new{randomization} part is given by
\[\new{P_R} f(x,v) = \mathbb E\po f\po x,\eta v +\sqrt{1-\eta^2}  G\pf \pf\,,\]
with $\eta\in[0,1)$, where $G$ is a standard $d$-dimensional Gaussian variable. We call \new{g}HMC (resp. stochastic gradient \new{g}HMC -- SG\new{g}HMC) chain with parameters $(\delta,K,\eta)$ the Markov chain on $\R^{2d}$ with transition
\begin{equation}\label{def:PPs}
P=\new{P_R} P_V^K \new{P_R}\qquad \text{(resp. $P_s= \new{P_R}P_{sV}^K\new{P_R}$).}
\end{equation}

For $M$ a $2d\times 2d$ positive definite symmetric matrix, we write $\R^{2d}\ni z \mapsto \|z\|_M := \sqrt{z\cdot M z}$ the associated Euclidean norm,  $\mathrm{dist}_M:\R^{2d}\times\R^{2d}  \ni (z,y) \mapsto \|z-y\|_M$ the corresponding distance and, for $r\in[ 1,+\infty]$ 
 and $\nu,\nu'$ two probability measures on $\R^{2d}$, 
\[\mathcal W_{M,r}(\nu,\nu') \ = \ \inf_{\mu \in \Pi(\nu,\nu')} \|\mathrm{dist}_M\|_{L^r(\mu)}\]
the associated $L^r$-Wasserstein distance,  where $\Pi(\nu,\mu)$ is the set of couplings of $\nu$ and $\nu'$, namely the set of probability laws on $\R^{2d}\times\R^{2d}$ with $2d$-dimensional marginals $\nu$ and $\nu'$. When $M$ is the identity we simply write $\|z\|_M=|z|$ and $\mathcal W_{M,r}=\mathcal W_r$. We denote by $\mathcal P_2$ the set of probability measures on $\R^{2d}$ with a finite second moment.

\subsection{The Gaussian case}\label{sec:results:Gaussian}

In this section, we consider $P$ the Markov transition operator of the \new{g}HMC chain given in \eqref{def:PPs} and address the question of the optimal choice of the parameters $\param=(\delta,K,\eta)$   for target distributions which are in a given family of Gaussian laws. Gaussian distributions are often used as a theoretical benchmark for comparing MCMC samplers and scaling parameters since convergence rates have  explicit expressions. When this is made at the level of continuous time processes (as in \cite{HHS1,Lelievre2012,MonmarcheGuillin,Arnold}), the question of the correct normalisations for a fair comparison between different processes makes the results difficult to interpret. Here we don't have this problem since we work directly with the numerical schemes. A comparison of the convergence rates (with respect to the number of gradient evaluations) of two schemes is fair if they have the same accuracy on the target measure.  

As already mentioned in the introduction, let us emphasize that, contrary to e.g. \cite{HHS1,Lelievre2012,MonmarcheGuillin,Arnold} and references within, our goal is not to tune the most efficient sampler for a single specific Gaussian distribution, but rather we want to find the parameters that give the best convergence rate over a family of Gaussian distributions. This is less precise but more robust, and the results are meant to give some indications on a suitable choice of parameters in practice in the general non-Gaussian case.
  
  More precisely, we consider $L>m>0$ and $d\geqslant1$ (fixed in all the section, unless otherwise specified) and denote by $\mathcal M_s(m,L)$ the set of $d\times d$ symmetric matrices $S$ with $m\rev{I_d}\leqslant S\leqslant L\rev{I_d}$. For $S\in\mathcal M_s(m,L)$  and parameters $\param=(\delta,K,\eta)$  with $\delta>0$, $K\in\N$, $\eta\in[0,1)$,  we write $\pi_S$ the centered Gaussian distribution with variance $S^{-1}$ and $P_{S,\param}$ the Markov transition operator corresponding to $P$ in \eqref{def:PPs} with parameters $\param$ and potential $U(x)=x\cdot S x/2$ for all $x\in\R^d$. 

The proofs of all the results stated in this section are given in Section~\ref{sec:Gauss:proof} in the appendix.  
\begin{proposition}\label{prop:Gauss_error}
Let   $S\in\mathcal M_s(m,L)$  and $\param=(\delta,K,\eta)$  be such that $\delta^2L<4$. Denote by  $S(\param)$ the symmetric matrix that has the same eigenvectors as $S$ and such that an eigenvector associated to an eigenvalue $\lambda$ for $S$ is associated to $\lambda (1-\delta^2\lambda/4)$ for $S(\param)$. Then  $\mu_{\param}:=\pi_{S(\param)}\otimes\pi_{I_d}$ is invariant for $P_{S,\param}$ and
\begin{equation}\label{eq:deftolerance}
\varepsilon(\param) := \sup_{S\in\mathcal M_s(m,L)} \mathcal{W}_2\po \pi_S,\pi_{S(\param)}\pf =\sqrt{d} \frac{1 - \sqrt{1-\delta^2L/4}}{\sqrt{L}} \,.
\end{equation}
\end{proposition} 
\new{This result is consistent with \cite[Example 2]{pmlr-v75-wibisono18a} concerned with the overdamped case ($K=1$, $\eta=0$).}
 
This means that two sets of parameters yield the same accuracy (and thus can be fairly compared) if they have the same   time-step $\delta$. Moreover, for a desired accuracy $\varepsilon>0$ with $\varepsilon <\sqrt{d/L}$, the required time-step is
\begin{equation}\label{eq:deltachoix}
\delta =  2  \sqrt{\frac{1-(1-\varepsilon \sqrt {L/d})^2}{L}} \ \underset{\varepsilon \sqrt {L/d}\rightarrow 0} \simeq \frac{\sqrt{8\varepsilon}}{(dL)^{1/4}} \,.
\end{equation}

\begin{remark}
More generally, if we are interested in a accuracy in term  for instance  of total variation  norm or relative entropy, then, again,
since $\pi_{S(\param)}$ is independent from $K$ and $\eta$,   two sets of parameters have the same accuracy is they have the same time-step $\delta$.
\end{remark}

 We are interested in the worst convergence rate of $P_{S,\param}$ for all $S\in \mathcal M_s(m,L)$, namely
 \begin{equation}
 \label{eq:defrho}
 \rho(\param) =  \inf_{S\in\mathcal M_s(m,L) } \lim_{n\rightarrow +\infty} - \frac{1}{Kn} \ln \po \sup_{\nu\in\mathcal{P}_2(\R^{2d})\setminus\{\mu_{\param}\}} \frac{\mathcal{W}_2(\nu P_{S,\param}^n,\mu_{\param})}{\mathcal{W}_2(\nu ,\mu_{\param})} \pf \,,
 \end{equation}
 since $Kn$ is the number of gradient computations needed to perform $n$ transitions of the Markov process.

\begin{remark}\label{rem:rate}
Here we only compare the asymptotic convergence rate, in particular if for instance
\[\sup_{\nu\in\mathcal{P}_2(\R^{2d})\setminus\{\mu_{\param}\}} \frac{\mathcal{W}_2(\nu P_{S,\param}^n,\mu_{\param})}{\mathcal{W}_2(\nu ,\mu_{\param})}  \underset{n\rightarrow +\infty}\simeq C e^{-nK \rho(\param)}\,,\]
we do not take the prefactor $C$ into account, in order to have a simpler optimization problem. This could be an issue if optimizing $\rho(\param)$ lead to values of $\param$ for which $C$ depends badly on $\delta$ (hence on the accuracy) or $d$. We will see that this is not the case. In fact, the reasonable choices for $\param$ are such that, as $\delta$ vanishes, the chain converges to some continuous-time process, either of Langevin or of idealized HMC type, and in these cases the prefactor $C$ converges to some finite value which is independent from $d$ and at most polynomial in $m$, $L$ and $t=n\delta$, which means the prefactor   leads at most to a logarithmic term in $m,L,\varepsilon,d$ in the number of steps needed to achieve a given accuracy $\varepsilon$. See also Remark~\ref{rem:NormeAn} below. 
\end{remark}

\begin{remark}\label{rem:rescaling}
If $(x_n,v_n)_{n\in\N}$ is a Markov chain associated to $P_{S,\param}$ for some parameters $\param=(\delta,K,\eta)$ and a matrix $S\in \mathcal M_s(m,L)$ then, for all $\alpha>0$, $(\sqrt \alpha x_n,v_n)_{n\in\N}$ is associated to $P_{S/\alpha,\tilde{\param}}$  where  $\tilde{\param}=(\sqrt{\alpha}\delta,K,\eta)$. By homogeneity of the Wasserstein distance, $\varepsilon(\tilde{\param}) = \sqrt{\alpha}\varepsilon(\param)$ and $\rho(\param)=\rho(\tilde{\param})$. Moreover, by decomposing along an eigenbasis of $S$, it is clear that $\rho(\param)$ is independent from $d$ (see the next proposition). This means that, for given $L>m>0$, $\varepsilon>0$ and $d\geqslant 1$, considering the question of finding the parameters $\param$ that maximizes $\rho(\param)$ under the constraint that $\varepsilon(\param)\leqslant \varepsilon$, we can work with the following scaling-invariant parameters: the rescaled time-step $\delta\sqrt{L}$, the relative error $\varepsilon \sqrt {L/d}$, and the condition number $L/m$ (while $K$ and $\eta$ are already invariant by scaling). In other words, the optimal parameters $(\delta,K,\eta)$ for some $m,L,\varepsilon$ are such that, necessarily, $(\delta \sqrt L,K,\eta)$ can be expressed as a function of $\varepsilon\sqrt{L/d}$ and $L/m$ alone. Or, in other words, from Proposition~\ref{prop:Gauss_error}, at a fixed $\delta$, the parameters $(K,\eta)$ which maximizes $\rho(\delta,K,\eta)$ can be expressed as a function of $\delta\sqrt L$ and $L/m$.
\end{remark}

\begin{proposition}\label{prop:GaussCVlongtime}
For  $\param=(\delta,K,\eta)$ with $\delta^2 L<4$,
\[\rho(\param)= \frac{-\ln g\po h(K,\delta),\eta\pf }{K } \]
 where, writing $\varphi_\lambda = \arccos(1-\delta^2\lambda /2)$ for $\lambda \in[m,L]$, 
 \begin{eqnarray*}
 h(K,\delta) &= & \sup\{|\cos(K\varphi_\lambda)|,\ \lambda \in [m,L] \}\\
 g(c,\eta) &=& \eta  \vee \frac{(1+\eta^2)c}2  + \sqrt{\po \po \frac{(1+\eta^2)c}2 \pf^2 - \eta^2\pf_+}\,.
 \end{eqnarray*}
\end{proposition}

\begin{figure}
\begin{center}
\includegraphics[scale=0.4]{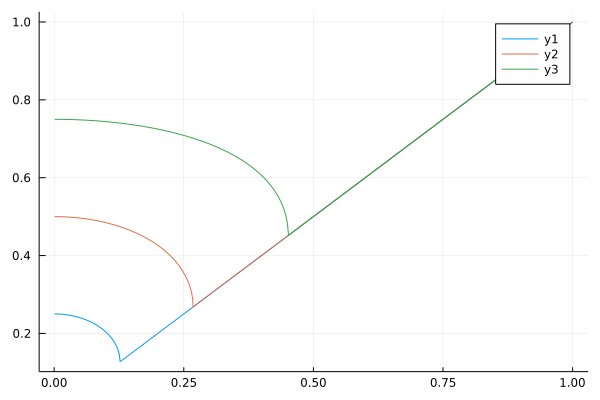}
\caption{$\eta\mapsto g(\eta,c)$ with $c=1/4,1/2,3/4$}\label{fig:g}
\end{center}
\end{figure}

See Figure~\ref{fig:g} for the graph of $\eta\mapsto g(\eta,c)$ for some values of $c$.

In the unadjusted HMC sampler, as $\delta\rightarrow 0$, typically, $K$ is of order $1/\delta$, and $\eta$ is constant ($\eta=0$ in the most classical case). When using a numerical integrator for the Langevin process, typically, $K=1$ and $\eta=e^{-\gamma \delta} = 1-\gamma\delta +o(\delta)$ where $\gamma>0$ is the friction parameter. Considering the regime $\delta\rightarrow 0$ in the next result shows that these are in fact the only two   relevant scalings:

 \begin{proposition}\label{prop:GaussienScaling}
 Let  $(\param_n)_{n\in\N}=(\delta_n,K_n,\eta_n)_{n\in\N}$ be a sequence of parameters with $\delta_n\rightarrow 0$ as $n\rightarrow +\infty$. Then, up to extracting a subsequence (still denoted by $(\param_n)_{n\in\N}$), we are necessarily in one of the three following cases:
 \begin{enumerate}  
 \item There exists $T>0$ and $\eta\in[0,1)$ such that $K_n\delta_n \rightarrow T$  and $\eta_n\rightarrow \eta$. In that case,
   \[\rho(\param_n) \ \underset{n\rightarrow+\infty}\simeq  \  \frac{\delta_n|\ln g(h_*(T),\eta)|}{T}\ :=\ \delta_n \sqrt L \bar  \rho_{HMC}(T,\eta)\]
   with $h_*(T)=\sup\{|\cos(x)|,\ x\in[T\sqrt m,T\sqrt L]\}$. We call this the HMC scaling.
   \item $K_n\delta_n \rightarrow 0$   and there exists $\gamma>0$ such that $\eta_n=1-\gamma K_n\delta_n +o(K_n\delta_n)$. In that case,
       \[\rho(\param_n) \ \underset{n\rightarrow+\infty}\simeq \po\gamma  - \sqrt{(\gamma^2-m)_+} \pf\delta_n \ :=\ \delta_n \sqrt L \bar \rho_{Lang}(\gamma) \,.  \ \]
       We call this the Langevin scaling.
\item $\rho(\param_n) =  o(\delta_n)$ as $n\rightarrow +\infty$.
 \end{enumerate}
 \end{proposition}

\begin{remark}\label{rem:NormeAn}
As mentioned earlier, we see that the relevant scalings for the parameters are such that, as $\delta$ vanishes, the Markov chain converges toward a continuous process, either of HMC or Langevin type. For instance, in the Langevin case, for a fixed $t>0$, we get that
\[\mathcal{W}_2(\nu P_{S,\param}^{\lfloor t/\delta\rfloor},\mu_{\param}) \underset{\delta\rightarrow 0}{\longrightarrow} \mathcal{W}_2(\nu P_t,\rpi)\]
where $(P_t)_{t\geqslant 0}$ is the semi-group associated with the Langevin diffusion \eqref{eq:EDScontinuousLangevin2} with $\na U(x) = Sx$. As can be seen e.g. by applying   Lemma~\ref{prop:Gaussien} from the appendix, we get
\[\mathcal{W}_2(\nu P_t,\rpi) \leqslant|e^{tB}| \mathcal{W}_2(\nu,\rpi)\,,\qquad \text{with}\qquad
B = \begin{pmatrix}
0 & I_d \\ -S & - \gamma I_d
\end{pmatrix}\,. \]
By diagonalizing $S$, it is not hard to see that $|e^{tB}| $ can be upper bounded by an expression of the form $(C_1+C_2 t) e^{-t\sigma(B)}$ where $C_1$ and $C_2$ depend only on $m$ and $L$ and $\sigma(B)$ is the spectral gap of $B$.  As discussed in Remark~\ref{rem:rate}, this means that, in the Langevin scaling, we only lose some logarithmic terms in $m,L$ and $t=n\delta$ by considering only the asymptotic convergence rate and not the contraction constant at a fixed $t$. The case of the HMC scaling is similar.
\end{remark}

 Finally, we turn to the question of optimizing the convergence rate at a fixed accuracy. We consider both the regime $\delta_n \rightarrow 0$ and the case of a fixed $\delta>0$.  For $c\in[0,1]$, write  $\eta_*(c) = (1-\sqrt{1-c^2})/c$, which is where the minimum of $\eta\mapsto g(c,\eta)$ is attained (see  Lemma~\ref{lem:g} in the appendix).
 
 \begin{proposition}\label{prop:Gaussoptim}
  First, consider the notations of Proposition~\ref{prop:GaussienScaling}
 
 \begin{itemize}
 \item \textbf{Langevin scaling:} The maximum of $\bar \rho_{Lang}$ is $\sqrt {m/L}$, attained at $\gamma= \sqrt m$.
 \item \textbf{HMC scaling:} The maximum of $\bar \rho_{HMC}$ is attained at $T = T_*:=\pi/(\sqrt{L}+\sqrt m)$ and $\eta = \eta_*\po \cos (T_*\sqrt m)\pf$. More precisely, for all $T\in[0,T_*]$,  $\eta\mapsto \bar \rho_{HMC}(T,\eta)$ is maximal at $\eta=\eta_*\po \cos (T\sqrt m)\pf$, and then $T\mapsto \bar \rho_{HMC}(T,\eta_*\po \cos (T\sqrt m)\pf)$ is  increasing over $[0,T_*]$ and goes to $\sqrt{m/L}$ as $T\rightarrow 0$.
 
 \item \textbf{position HMC:} The maximum of $T\mapsto  \bar \rho_{HMC}(T,0)$ is also attained at $T=T_*$, with $\bar \rho_{HMC}(T_*,0)=-\ln \po \cos(T_*\sqrt m)\pf /(T_*\sqrt L)$.
 \end{itemize}
 
 Second, we consider a fixed $\delta>0$. Assume that $2\varphi_m < \varphi_L \leqslant \pi/2$ and let $K_*=1+\left \lfloor \pi/(\varphi_m+\varphi_L)\right\rfloor $.
\begin{itemize}
\item The maximum of $(K,\eta) \mapsto \rho(\delta,K,\eta)$ is attained for $\eta = \eta_*(h(K,\delta))$ and one of the values $K\in\{K_*-1,K_*\}$.
\item The maximum of $K \mapsto \rho(\delta,K,0)$ is attained for  one of the values $K\in\{K_*-1,K_*\}$.
\end{itemize} 
 \end{proposition}
 
 \begin{remark}\label{rem:conditionsSupp}
 In the fixed $\delta$ case, the additional assumption is meant for the simplicity of the analysis. We expect the result to hold without it, but then the study of $K\mapsto h(K,\delta)$ is more cumbersome due to some possible pseudo-periodicity issues. Moreover, the additional condition typically holds in practice: for badly conditioned problems, $4m<L$ and thus  $2\varphi_m < \varphi_L$ for $\delta$ small enough, and $\varphi_L \leqslant \pi/2$ is simply $\delta^2 L <2$, which is close to the condition $\delta^2 L<4$ which is necessary for the bare stability of the chain, and from Proposition~\ref{prop:Gauss_error} is equivalent to say that the required relative error $\varepsilon\sqrt{L/d}$ is less than $1/2$.
 \end{remark}
 
 As a conclusion of the analysis in the Gaussian case, we can express these last results in terms of the scaling-invariant parameters of Remark~\ref{rem:rescaling}. For a given condition number $\kappa=L/m$ and a given small relative error $\varepsilon'= \varepsilon \sqrt{L/d}$, the optimal choice of parameters is to take a rescaled time-step $\delta'=\delta \sqrt L = 2\sqrt{1-(1-\varepsilon')^2}$ (this is \eqref{eq:deltachoix}) and, considering the second part of Proposition~\ref{prop:Gaussoptim},
 \begin{equation}\label{eq:bestKeta}
 K \underset{\varepsilon'\rightarrow 0}\simeq \frac{\pi}{\delta' \po 1+ 1/\sqrt\kappa\pf}\,,\qquad \eta \underset{\varepsilon'\rightarrow 0}\simeq \frac{1-\sin \po \pi/(1+\sqrt\kappa)\pf }{\cos\po \pi/(1+\sqrt\kappa)\pf }\,,
 \end{equation}
these equivalents as $\varepsilon'\rightarrow 0$ being uniform in $\kappa>1$. It also corresponds to the best choice in the HMC scaling in the first part of Proposition~\ref{prop:Gaussoptim}. As $\varepsilon'\rightarrow 0$ and $\kappa\rightarrow +\infty$, this yields a convergence rate
\[\rho(\delta,K,\eta) \underset{\varepsilon'\rightarrow 0}\simeq  \frac{\delta' \po 1+ 1/\sqrt\kappa\pf}{\pi} \ln \po \frac{\cos\po \pi/(1+\sqrt\kappa)\pf }{1-\sin \po \pi/(1+\sqrt\kappa)\pf }\pf \underset{\kappa\rightarrow +\infty }\simeq \frac{\delta'}{\sqrt\kappa}\,.\]
 This can be compared, on the one hand, to the best choice in the Langevin scaling, which would be any $K$ such that $K\delta\rightarrow 0$ as $\varepsilon'\rightarrow 0$ (we can simply consider $K=1$ for instance), and then $\eta = 1- \sqrt{m}\delta K = 1 - \delta'K/\sqrt\kappa$ and, on the other hand, to the best position HMC choice, namely $\eta=0$ and $K$ as in \eqref{eq:bestKeta}. In the best Langevin case, from Proposition~\ref{prop:Gaussoptim}, we get a rate $\delta'/\sqrt\kappa$ when $\varepsilon'\rightarrow 0$, which is thus equivalent as $\kappa\rightarrow +\infty$ to the rate obtained with the optimal choice of parameters. In the best position HMC case, we get
\[ \rho(\delta,K,0)  \underset{\varepsilon'\rightarrow 0}\simeq  \frac{\delta' \po 1+ 1/\sqrt\kappa\pf}{\pi} \ln \po \frac1{\cos\po \pi/(1+\sqrt\kappa)\pf } \pf \underset{\kappa\rightarrow +\infty }\simeq \frac{\delta' \pi}{\kappa}\,,\]
 which is off to the optimal rate by a factor $\pi/\sqrt{\kappa}$.
 
 Finally, let us notice that, from Proposition~\ref{prop:GaussCVlongtime}, we see that $\rho(\param)=0$ as soon as $h(K,\delta)=1$, since $g(1,\eta)=1$ for all $\eta\in[0,1)$. Assuming for simplicity (as in Remark~\ref{rem:conditionsSupp})  that $4<\kappa $ and that  $\delta'$ is small enough, it means that  $\rho(\param)=0$ as soon as $\pi \leqslant K\varphi_L\simeq K\delta'$, which is close to the optimal choice given by \eqref{eq:bestKeta}. This explains the sensitivity of the optimal \new{g}HMC parameters when $\kappa$ is large, as discussed in the introduction, by comparison with the optimal parameters in the Langevin case ($K=1$).
 
  \new{
\begin{remark}
The fact $\rho(\mathfrak{p})=0$ when $K\delta' \geqslant \pi$ is due to the existence of $S\in\mathcal M_s(m,L)$ such that $K\delta'$ is a period of $(\dot x,\dot v)=(v,-Sx)$, and similarly for the Verlet scheme. In particular, for the corresponding \new{g}HMC chain, $x_n=x_0$ for all $n\in\mathbb N$ almost surely. However, one may argue that, in practice along the computation of $z_{n+1}$ from $z_n=(x_n,v_n)$, a trajectory $\left( \varphi^{\circ j}(x_n,\eta v_n + \sqrt{1-\eta^2} G_n)\right)_{1\leqslant j \leqslant K}$ is computed so that, denoting by $x_{n+j/K}$ the corresponding position, an estimator for $\int fd \rpi$ can be
\[J_N := \frac1N \sum_{n=1}^N \frac 1K\sum_{j=1}^K f(x_{n+j/K}),\]
for which the position $x_{n+j/K}$ is not constant as soon as $\delta' < \pi$. However, this doesn't solve the problem. Indeed, considering a non-negative $f$, we see that $J_N \geqslant f(x_0)/K$ for all $N \in\mathbb N$ in the case of a periodic orbit. Since $f$ may be a $\mathcal C^\infty$ with $f(x_0)=1$ and $\int fd\rpi$ arbitrarily small, this clearly prevents convergence. In fact, letting $K\rightarrow \infty$ for a fixed $f$ wouldn't work either, since the previous argument is easily adapted to show that the idealized \new{g}HMC where the exact Hamiltonian dynamics is solved also suffers from periodicity. More generally, in practice, even if $K\delta'$ is not exactly a period of the dynamics, we still get that the convergence may be arbitrarily slow as $K\delta'$ get close to a period of the dynamics, since in that case $x_n$ moves slowly and so does the law of the orbits $(x_{n+j/K})_{j=1..K}$. 

On the other hand, non-convergence can only concern the high-frequency modes (i.e. the dimensions with variance less  than $(K\delta/\pi)^2$)  
\end{remark}
}

\subsection{Dimension-free convergence rate}\label{sec:results:Wasserstein}

In this section, we consider the SG\new{g}HMC chain with transition operator $P_{sV}$ given by \eqref{def:PPs}. We work under the following condition on $b$, which is the extension of the usual convex/smooth condition to the stochastic gradient case.

\begin{assu}\label{Assu:mLb}
There exist $m,L>0$ such that for all $\theta\in\Theta$ and $x,u\in\R^d$,
\[  u\cdot \na_x b(x,\theta)u \ \geqslant \ m |u|^2 \,,\qquad  |\na_x b(x,\theta)u| \ \leqslant \ L|u|.\]
\end{assu}

This condition is satisfied for instance if $b(x,\theta) = \sum_{i\in \theta} \na U_i(x)$ where $\theta \in \Theta\subset \cco 1,N\ccf^n$ for some $n\leqslant N$ and $(U_i)_{i\in\cco 1,N\ccf}$ are  strongly convex and gradient-Lipschitz $\mathcal C^2$ functions, for instance $U_i(x)= \| A_i x-y_i\|^2$ for some vectors $y_i$ and injective matrices $A_i$. This corresponds to the method of subsampling with mini-batches.

\bigskip


 Given two initial conditions $z_0,z_0'\in\R^{2d}$, we say that $(z_k,z_k')_{k\in\N}$ is a parallel (or synchronous) coupling of two SG\new{g}HMC chains if both $(z_k)_{k\in\N}$ and $(z_k')_{k\in\N}$ are SG\new{g}HMC chains and they are constructed using at each step the same Gaussian variables $(G_k,G_k')$ (in the two \new{randomization} steps of a transition) and the same auxiliary variables $(\theta_{k,j})_{j\in \cco 1,K\ccf}$ (in the $K$ stochastic gradient Verlet steps).
 
 \new{In the next result, which is the main one of this section and is proven in the appendix in Section~\ref{sec:preuveParallel}, three points are distinguished: the goal of the two first is to give dimension-free convergence rates with the correct dependency in $K\delta$ and $1-\eta$ over a wide range of parameters, from position HMC to Langevin (the first point is just the particular case $\eta=0$, which is simpler). The third point only considers a subset of these parameters  (specifically, it is restricted to the strong inertia case, in the sense that $K\delta$ is sufficiently small and $1-\eta$ is of order $K\delta$), the interest with respect to the second point being that, as discussed below, in this regime we get a better dependency in terms of the condition number.
 }
 
 \begin{theorem}\label{thm:crude}
Under Assumption~\ref{Assu:mLb}, suppose that $2\delta \sqrt{L} \leqslant 1 $. Let $(z_n,z_n')_{n\in\N}$ be a parallel coupling of two SG\new{g}HMC  chains.
\begin{enumerate}
\item If $\eta=0$\rev{, $K\delta\leqslant 1/(2\sqrt{L})$  and $\delta  \leqslant K\delta m/(31L)$} then, denoting by $x_n$ and $x_n'$ the first $d$-dimensional coordinates of $z_n$ and $z_n'$, almost surely, for all $n\in\N$,
\begin{equation}\label{eq:rho_cas0}
|x_n-x_n'| \leqslant \rev{\po 1 - \frac{m}{20}(\delta K)^2 \pf^n} |x_0-x_0'|\,.
\end{equation}
\item For any $\eta\in[0,1)$, assume that
\begin{equation}\label{eq:conditioncrude}
1-\eta^2   \geqslant  \frac{\new{7} K\delta L^{3/2}}{m}\,,\qquad\text{and}\qquad  e^{\new{3}K\delta\sqrt{L}/\new{2}} K\delta   \leqslant  \frac{m^{3/2} }{40L^2(1-\eta^2)}\,.
\end{equation}
Then, there exists a symmetric matrix $M$ with,  for all $z=(x,v)\in\R^{2d}$,
  \begin{equation}\label{eq:norme_equivalence2}
  \frac{2}{3}\po |x|^2 + \frac{2m}{L^2}|v|^2\pf \leqslant  \|z\|_M^2 \leqslant \frac{4}{3}\po |x|^2 + \frac{2m}{L^2}|v|^2\pf 
  \end{equation}
and such that,  for all $n\in\N$,
\begin{equation}\label{eq:contractMz}
\|z_n-z_n'\|_M \leqslant \po 1 -\rho\pf^n \|z_0-z_0'\|_M
\end{equation}
with
\begin{equation}\label{eq:rho_cas1}
 \rho = \frac{ (\delta K)^2m}{40(1-\eta^2)}\,.
\end{equation}
  \item \new{For $\eta>0$, assume that 
\begin{equation}\label{eq:conditioncrude2} 
  \overline{\gamma}:=\frac{1-\eta^2}{K\delta\sqrt{L}\eta}\geqslant 2\qquad\text{and}\qquad K\delta   \leqslant \frac{\eta m}{ L^{3/2} \overline{\gamma}(\modifplb{21}+\modifplb{9}\overline{\gamma}^2)}\,.
  \end{equation}
Then, there exists a symmetric matrix $M$ with,  for all $z=(x,v)\in\R^{2d}$,
  \begin{equation}\label{eq:norme_equivalence2bis}
  \frac{1}{2}\po |x|^2 + \frac{1}{L}|v|^2\pf \leqslant  \|z\|_M^2 \leqslant \frac{3}{2}\po |x|^2 + \frac{1}{L}|v|^2\pf 
  \end{equation}
and such that,  for all $n\in\N$, \eqref{eq:contractMz} holds with 
\begin{equation}\label{eq:rho_cas2}
\rho = \frac{\eta m}{\modifplb{6} \sqrt{L}\overline{\gamma}} K\delta \,.
\end{equation}
}
\end{enumerate}
\end{theorem}

From this result, the proof of the following is standard (see e.g. \cite{MonmarcheSplitting}) and is thus omitted.

\begin{corollary}\label{Cor:contraction}
Under the conditions and with the notations of Theorem~\ref{thm:crude}, for all probability measures $\nu,\mu$ on $\R^{2d}$, $p\geqslant 1$ and $n\in\N$,
\[\mathcal W_{M,p}\po \nu P_{s}^n ,\mu P_s^n\pf \leqslant (1-\rho)^n \mathcal W_{M,p}\po \nu ,\mu\pf\,.\]
Moreover, $P_s$ admits a unique invariant law $\tilde \mu_s$, which has finite moments of all orders.
\end{corollary}

Thanks to the norm equivalence \eqref{eq:norme_equivalence2}, this implies a similar contraction for the standard $\mathcal W_p$ distance. A convergence in total variation can be deduced by combining this Wasserstein convergence with a Wasserstein/total variation regularization result, as \cite[Proposition 3]{MonmarcheSplitting} in the case $K=1$ or \cite[Lemma 16]{BouRabeeEberle}  in the position HMC case.

\medskip

Let us now discuss the explicit conditions and bounds of Theorem~\ref{thm:crude}. First, as in Remark~\ref{rem:rescaling}, due to the scaling properties of the SG\new{g}HMC (as can be see in the proof of Theorem~\ref{thm:crude}), the conditions \eqref{eq:conditioncrude} and the bound on $\rho$ can be expressed in term of the rescaled time-step $\delta\sqrt{L}$ and of the condition number $L/m$. In the rest of this discussion we assume that $L=1$.

The main point of Theorem~\ref{thm:crude}  is the following: the conditions on the parameters and the convergence rate are independent from the dimension and, for a fixed $\delta$, the convergence rate with respect to the number of gradient computations, i.e. $-\ln(1-\rho)/K $, is lower bounded uniformly over a wide range of values of $K$ and $\eta$. Moreover, as $\delta$ vanishes, both in the regime $K\delta \rightarrow T>0$ small enough with $\eta\in[0,1)$ independent from $\delta$ (small enough depending on $T$) and in the regime $K\delta\rightarrow 0$ with $1-\eta^2 = \gamma K\delta + o(K\delta)$ for some $\gamma>0$ large enough,  this convergence rate  is of order $\delta$, which is optimal as seen in Section~\ref{sec:results:Gaussian} in the Gaussian case. 

Notice that the second part of the theorem yields $\rho = (\delta K)^2m/40$ in the case $\eta=0$. The loss of a factor \rev{$1/2$} with respect to the first part is a simple consequence of some rough bounds used in the proof in order to keep  simple the condition \eqref{eq:conditioncrude} and other computations. \rev{More importantly, applying the general result for $\eta=0$ doesn't give sharp rates in terms of the condition number because of the conditions~\eqref{eq:conditioncrude} on the parameters, since} \new{the optimal value of $\rho/K$ based on \eqref{eq:rho_cas0} is obtained by saturating the constraint $5\delta K \leqslant m$, in which case $\rho/K = 2\delta m^2$. \rev{By contrast, the first item of the theorem gives an optimal value of $\rho/K$ equal to $mT\delta /(20L)$ (linear in the inverse of the condition number, which is sharp as seen in the Gaussian case) with an integration time $T=K\delta = 1/(2\sqrt{L})$, recovering the analogue of the result of Chen and Vempala \cite{ChenVempala}  in the continuous-time limit (our additional condition $\delta  \leqslant T m/(31L)$ disappears in this continuous-time limit where $\delta$ goes to zero with $T$ being fixed).} 

In the Langevin case $K=1$, in order to maximize $\rho/K$ when $\rho$ is given by \eqref{eq:rho_cas1}, the choice of $\eta$ should saturate the first part of \eqref{eq:conditioncrude}, and this gives $\rho/K = \delta m^2 / 280$. Again, the dependency is sharp for $\delta$ but not for $m$, as a convergence rate of order $m$ is obtained in \cite{dalalyan1} (for the continuous-time underdamped Langevin diffusion). This is why we also provided the alternative estimate \eqref{eq:rho_cas2}. For $K=1$, taking $\eta$ so that $\overline{\gamma}=2$, we get that $\rho/K$ is of order $m\delta$. In other words, in terms of the condition number, in the general strongly convex case, we have obtained for the Langevin splitting scheme the same dependency as the continuous-time process in \cite{dalalyan1} (in particular we improve upon the rates of \cite{MonmarcheSplitting}, which are of order $m^2$). Besides, we get the same result in the general gHMC scaling as long as $K\delta$ is of order $m$ (this is the second part of \eqref{eq:conditioncrude2}), which is small for badly conditioned problems but independent from $\delta$, namely this goes beyond the underdamped Langevin process and corresponds rather to an HMC scaling with a strong inertia.  } \rev{Besides, concerning the comparison between the two results \eqref{eq:rho_cas1} and \eqref{eq:rho_cas2}, notice that, at first sight, the condition on $1-\eta$ in the first part of  \eqref{eq:conditioncrude} may seem more restrictive than the first part of \eqref{eq:conditioncrude2}. However, the conditions \eqref{eq:conditioncrude2} should be thought as follows: first, choose a value of $\bar{\gamma}\geqslant 2$; then, choose $K\delta$ small enough to satisfy the second part of \eqref{eq:conditioncrude2} with this value of $\bar{\gamma}$; finally, take $\eta = 1 - c K\delta$ with $c$ small enough for the first part to hold. Indeed, if we start by fixing $\eta$ to some value in $(0,1)$, and then look for $K\delta$ to fulfill the conditions of \eqref{eq:conditioncrude2}, the first condition requires $K\delta$ to be small enough but the second condition requires it to be large enough (because $\bar{\gamma}$ is proportional to $1/(K\delta)$ for a fixed $\eta$), and thus there may be no solution.  }


 In all cases, the fact that $K\delta$ has to be small enough is consistent with the analysis of the Gaussian case of Section~\ref{sec:results:Gaussian}, since we should avoid the periodicity issues of the Hamiltonian dynamics. The fact that $\eta$ has to be small enough -- namely that the damping has to be strong enough --  corresponds to the condition that the friction $\gamma$ in the Langevin case has to be large enough, as discussed in \cite{MonmarcheSplitting,MonmarcheContraction}. This condition does not appear in the Gaussian case, where we get a convergence rate for all values of $\eta\in[0,1)$. However, as proven in \cite[Proposition 4]{MonmarcheContraction} for the continuous-time Langevin process \eqref{eq:EDScontinuousLangevin2}, it is in fact necessary in non-Gaussian cases in order to get the existence of a Euclidean norm $\|\cdot\|_M$ such that the associated Wasserstein distance is contracted by the transitions of the chain. More specifically, it is proven in \cite{MonmarcheContraction} that, for the continuous-time process \eqref{eq:EDScontinuousLangevin2}, if $\gamma \leqslant \sqrt{L}-\sqrt{m}$ then there exists $U\in\mathcal C^2(\R^d)$ with $m\rev{I_d}\leqslant \na^2 U \leqslant L\rev{I_d}$ such that there exists no $M,p$ such that $\mathcal W_{M,p}$ is contracted by the associated semi-group.   In other words, although the condition~\eqref{eq:conditioncrude} in the general case are not sharp, conditions of the form $K\delta \leqslant c_1$ and $1-\eta^2 \geqslant c_2$ where $c_1$ depends on $m,L$ and $c_2$ depends on $m,L$ and $K\delta$ are necessary to get a contraction of a norm. \new{In the high inertia case, the bound $\overline{\gamma} \geqslant 2 \sqrt{L}$ in \eqref{eq:conditioncrude2}  is of the same order in terms of $m,L$ as the optimal condition to get a contraction (since $\overline{\gamma}$ goes to the parameter $\gamma$ of the continuous-time process \eqref{eq:EDScontinuousLangevin2} as $\delta\rightarrow 0$  when, for instance, $K=1$ and $\eta = e^{-\delta \gamma/2}$).}

\new{
\begin{remark}
Concerning the dependency of the convergence rate in the condition number $m/L$, as mentioned in Section~\ref{subsec:overview_gaussien}, an interesting result of \cite{Cao2019OnE} is that, when choosing a damping parameter $\gamma$ of order $\sqrt{m}$ (which is  consistent with our analysis in the Gaussian case) the convergence rate in $L^2$ of the continuous-time kinetic Langevin diffusion scales like $\sqrt{m}$ (as we get in the Gaussian case for the Wasserstein distances, but we only get a rate of order $m$ in the general case; our proof would be easily adapted  in the case of a perturbation of the Gaussian case, namely assuming that the Lipschitz constant of $\na^2 U$ is small enough, and it is unclear whether this is a necessary condition to get a rate of order $1/\sqrt{m}$). In fact, as in Remark~\ref{rem:rescaling}, by rescaling, we can assume without loss of generality that $m=1$ in \cite[Theorem 1(i)]{Cao2019OnE}, so that this result reads: \emph{``The convergence rate in $L^2$ of the continuous-time kinetic Langevin diffusion with $\gamma=1$ is bounded below uniformly over all potential $U$ such that $\na^2 U \geqslant \rev{I_d}$" }. As we see, the key point of this result is that it doesn't involve any upper bound of $\na^2 U$, in particular the result applies to convex potentials with unbounded Hessian matrix. It is clear that this cannot hold for the discrete schemes  considered in the present work (neither in $L^2$ nor in Wasserstein distances), because they are not stable when $\na U$ is not Lipschitz (which doesn't mean that the Wasserstein contraction rate cannot be of order $\sqrt{m}$ when $\gamma=\sqrt{m}$ in the general strongly convex smooth case, however this is not known). It is unclear how the results of \cite{Cao2019OnE} could be used for some numerical schemes, and a related question is whether it could be adapted to other distances than the $L^2$ norm (for instance the relative entropy which behaves better with discretization errors \cite{VempalaWibisono,Chatterji}).  Notice that  \cite{Cao2019OnE} relies on the $L^2$ hypocoercivity method of \cite{albritton2021variational}, however with his initial modified norm method of Villani was already able to get the $L^2$ convergence for unbounded Hessian matrices in \cite{Villani2009} (so it may be possible to get the $\sqrt{m}$ scaling of \cite{Cao2019OnE} in the convex case with this approach, although this is not entirely clear), but not for the relative entropy. Both approaches of \cite{albritton2021variational} and \cite{Villani2009} in $L^2$ rely on Hilbert analysis. To our knowledge, the only result for hypocoercivity in relative entropy with unbounded Hessian matrices has been obtained in \cite{Cattiaux} and, as the analysis is more involved, it is unclear whether it would be possible to get sharp rates with this approach in the convex case.
\end{remark}
}

\medskip

From Theorem~\ref{thm:crude}, \cite[Theorem 23]{MonmarcheSplitting} yields Gaussian concentration inequalities for ergodic averages along the trajectory of the chain (see also \cite{Ollivier,Holmes2014} for exponential concentration bounds given a weaker $\mathcal W_1$ contraction). This yields non-asymptotic confidence intervals for the MCMC estimators of averages with respect to its biased equilibrium $\tilde\mu$. It can also be used to bound the quadratic risk of the estimator. However, this result requires the initial condition to satisfies a so-called log-Sobolev inequality. We now state a bound on this quadratic risk, based on a direct argument, which only requires a finite second moment for the initial condition.

\begin{proposition}\label{prop:risque_quadra}
Under the conditions and with the notations of Corollary~\ref{Cor:contraction}, denote by $\tilde \pi_s$ the first $d$-dimensional marginal of $\tilde \mu_s$\new{, and let $M,\rho$ be as in the second part of Theorem~\ref{thm:crude}}. For all $1$-Lipschitz function $f$ on $\R^d$, all probability measure $\nu_0$ on $\R^{2d}$ and all $n,n_0,N\in\N$, considering  $(Z_k^i)_{k\geqslant 0}=(X_k^i,V_k^i)_{k\geqslant 0}$ for $i\in\cco 1,N\ccf$ to be $N$ independent  SG\new{g}HMC chain with $Z_0^i\sim \nu_0$,
\[\mathbb E\po \left|\frac 1{nN}\sum_{i=1}^N\sum_{k=n_0+1}^{n_0+n}f(X_k^i)-\tilde \pi_{s}f\right|^2\pf \ \leqslant \ \po    \frac{3(1-\rho  )^{2 n_0}}{(\rho  n)^2     }   \mathcal W_{M,2}^2(\nu_0,\mu_s)+ \frac{6 }{nN \rho } \mathrm{Var}_M(\mu_s) \pf\]
where
\begin{eqnarray*}
\mathrm{Var}_M(\mu_s)  &:=& \frac12 \int_{\R^{2d}\times\R^{2d}}\|z-z' \|_M^2 \mu_s(\dd z)\mu_s(\dd z') \\
& \leqslant &  
\frac{22}{\rho} e^{25\delta K L /(\rho\sqrt{m})  } \po    (1-\eta^2) d + \frac{6 m}{L^{5/2}}  \delta K \new{\mathrm{V}_b }\pf \,,
\end{eqnarray*}
with  
 $\mathrm{V}_b  = \inf_{x\in\R^d}\mathbb E_{\omega} \po |b(x,\theta)|^2\pf$. 
\end{proposition}
The proof is given in the appendix in Section~\ref{sec:preuveEmpirical}.

In particular, in the case where $1-\eta^2 = \gamma \delta K$ for $\gamma$ large enough, we see that $\mathrm{Var}_M(\mu_s)  \leqslant C(d+\mathrm{V}_b )$ where $C$ only depends on $m,L$ and $\gamma$.

Notice that, under Assumption~\ref{Assu:mLb}, the vector field $\bar b(x) = \mathbb E_{\omega}(b(x,\theta))$ is such that the ODE $\dot x_t = \bar b(x_t)$ satisfies $|x_t| \leqslant e^{-mt}|x_0|$ for all $t\geqslant 0$, and thus the flow admits a unique fixed point, namely there exists a unique $x_*\in\R^d$ such that $\bar b(x_*)=0$. In particular, $\mathrm{V}_b $ is smaller than the variance of $b(x_*,\theta)$. 

\medskip

\new{Finally, Theorem~\ref{thm:crude}  and Proposition~\ref{prop:risque_quadra} only give the convergence of the algorithm to its biased equilibrium $\tilde \mu_s$. Hence, in order to obtain non-asymptotic efficiency bounds, it remains to control the error between  the marginal $\tilde \pi_s$ and the target measure $\rpi$.   However, once a Wasserstein contraction as in Corollary~\ref{Cor:contraction}  is obtained, a bound on the invariant measures easily follows from a finite-time bound (see  Section~\ref{section:erreurNum}). The fact that such bounds can be obtained uniformly over cases ranging from $K=1$ to $K=O(\delta^{-1})$ is clear since it follows from the fact that the \new{randomization} step does not induce any error.  

In the two next results (Propositions~\ref{prop:numerique} and \ref{prop:erreur_sto}), we consider: $\mu$ the probability measure with density proportional to $e^{-H}$ with $H(x,v) = U(x)+|v|^2/2$, $\tilde \mu$ the invariant measure of the (non-stochastic) \new{g}HMC chain $P$ and $\tilde\mu_s$ the invariant measure of the SG\new{g}HMC chain $P_s$.

The first result bounds the numerical error due to the Verlet integration of the Hamiltonian dynamics.

\begin{proposition}\label{prop:numerique}
Under Assumption~\ref{Assu:mLb}, assume furthermore that there exist $U\in\mathcal C^2(\R^d)$,  $\rho>0$ and $M$  a symmetric definite positive $2d\times 2d$ matrix such that   $b(\cdot,\theta) = \na U$ for all $\theta\in\Theta$ and \eqref{eq:contractMz} holds for the parallel coupling of the \new{g}HMC chain with transition matrix $P$. Then
 \begin{equation}\label{eq:erreur_num_prop1}
 \mathcal W_2(\mu,\tilde\mu) \leqslant  \tilde C \delta \sqrt{d}\qquad\text{with}\qquad \tilde C=  \frac{3\max(L,1)}{\sqrt{m}} (6|M||M^{-1}|)^{2 K\delta\sqrt{L}/\rho}    \,.
  \end{equation}
 If additionally there    exists $L_2>0$ such that $|\na^2 U(x)-\na^2 U(x')| \leqslant L_2|x-x'|$ for all $x,x'\in\R^d$ then
 \begin{equation}\label{eq:erreur_num_prop2}
 \mathcal W_2(\mu,\tilde\mu) \leqslant     \tilde C \delta^2  d  \qquad\text{with}\qquad \tilde C=  \frac{11\max(L,1)(L+L_2/\sqrt{L})}{m}(6|M||M^{-1}|)^{2 K\delta\sqrt{L}/\rho}  \,.
 \end{equation}
 \end{proposition}

 The proof is given in the appendix in Section~\ref{section:erreurNum}.

Again, the point is that, as $\delta$ vanishes, both in the regimes   $K\delta\rightarrow T>0$ with a fixed $\eta$ or $K\delta\rightarrow 0$ with $\eta=1-\gamma K\delta +o(K\delta)$, the constant $\tilde C$ is bounded uniformly in $\delta$. Notice that, in \eqref{eq:erreur_num_prop1} and \eqref{eq:erreur_num_prop2}, in the context of Theorem~\ref{thm:crude}, the operator norms $|M|$ and $|M^{-1}|$ are bounded in terms of $m$ and $L$ thanks to \eqref{eq:norme_equivalence2} or \eqref{eq:norme_equivalence2bis}.


\begin{remark}
As in \cite[Lemma 32]{MonmarcheSplitting}, due to the scaling properties of the $\mathcal W_2$ distance for separable laws, $d$ can be improved to $\sqrt{d}$ in \eqref{eq:erreur_num_prop2} in the separable case $U(x) = \sum_{i=1}^d U_i(x_i)$. In particular we recover the sharp scaling $\delta^2 \sqrt{d}$ of the Gaussian case \rev{(notice that the image of a gHMC chain by an orthonormal transformation is still a gHMC chain -- in other words the canonical basis plays no role in the algorithm -- and in particular in the Gaussian case the condition $U(x) = \sum_{i=1}^d U_i(x_i)$ always holds in an orthonormal eigenbasis of the covariance matrix).}
\end{remark}

The second result bounds the error due to the stochastic gradient method. For clarity, we assume that the stochastic gradient is given by an unbiased Monte Carlo method:

\begin{assu}\label{assu:gradient_sto}
The stochastic gradient is of the form
\begin{equation}\label{eq:forme_b}
b(x,\theta) = \frac{1}{p}\sum_{j=1}^p \na_x U_{\theta_i}(x)
\end{equation}
for some $p\geqslant 1$ where  $\theta= (\theta_1,\dots,\theta_p)$ are i.i.d. random variables on some set $\Theta_1$, and $U(x) = \mathbb E\po  U_{\theta_1}(x)\pf$. There exist $L\geqslant m>0$ such that, for all $\theta_1\in\Theta_1$ and $x\in\R^{2d}$, $m\rev{I_d}\leqslant \na_x^2 U_{\theta_1}(x) \leqslant L\rev{I_d}$, and $\na U(0)=0$. Finally, there exists $x\in\R^{2d}$ with $\mathbb E |\na_x U_{\theta_1}(x) - \na_x U(x)|^2 < \infty$.
\end{assu}

\begin{proposition}\label{prop:erreur_sto}
Under Assumption~\ref{assu:gradient_sto}, 
\[\mathcal W_2(\tilde\mu_s,\tilde\mu)   
  \leqslant \tilde C_s \frac{\sqrt{ \delta}}{\sqrt{p}} \sqrt{ 
  d L   +  \mathbb E\po  |\na U_{\theta_1}(0)|^2\pf } \]
  with
  \begin{equation}\label{eq:papoue}
  \tilde  C_s =  (6|M||M^{-1}|)^{2 K\delta\sqrt{L} /\rho}\frac{\sqrt{ \modifplb{30}   e^{4K\delta \sqrt{L}} \po C_2  +  2  \pf }}{\min(L,1)}\,.
  \end{equation}
  \end{proposition}

 The proof is given in the appendix in Section~\ref{section:erreurNum}. We recover the sharp dependency in $\sqrt{\delta/p}$, see \cite{chaterrjiSto,Dalalyan_sto,VZT}. As previously, in the context of Theorem~\ref{thm:crude},  in \eqref{eq:papoue},   $|M||M^{-1}|$ is bounded thanks to  \eqref{eq:norme_equivalence2} or \eqref{eq:norme_equivalence2bis}.

\medskip

In the non-stochastic case, combining Propositions~\ref{prop:risque_quadra} and \ref{prop:numerique}, we see that, in order to get a quadratic risk smaller than some $\varepsilon>0$ for a given $1$-Lipschitz observable $f$, in the non-stochastic gradient case, disregarding the dependency in $m,L$ and $L_2$, assuming that the initial condition is such that $\mathcal W_2(\nu_0,\rpi)$ is of order $\sqrt{d}$ (which is feasible, see   e.g.  \cite[Section 3.4]{MonmarcheSplitting}) we can take a time-step $\delta$ of order $\sqrt{\varepsilon/d}$, a burning time $n_0=0$ and a number of iterations $n$ of order $d/(\varepsilon\rho)$. In other words, we can obtain 
\[\mathbb E\po \left|\frac 1n\sum_{k=1}^{n}f(X_k)-\rpi f\right|^2\pf \ \leqslant\varepsilon \]
with a number of gradient computations $Kn = \mathcal O((d/\varepsilon)^{3/2})$  for all the choices of parameters such that $\rho/K$ is of order $\delta$ as $\delta$ vanishes, ranging from position HMC to Langevin splitting.

If the objective is to have $\mathcal W_2(\mathcal Law(x_n),\rpi)\leqslant \varepsilon$ (which is the most standard criterion in previous related works, although it is rightly argued in \cite{dalalyan2} that,  due to the scaling properties of the Wasserstein distance, it is more natural to ask for $\mathcal W_2(\mathcal Law(X_n),\rpi)\leqslant \varepsilon \mathcal W_2(\delta_0,\rpi)$)  then, from Corollary~\ref{Cor:contraction} and Proposition~\ref{prop:numerique}, in the non-stochastic case, taking $\delta $ of order $\sqrt{\varepsilon/d}$, it is sufficient to take $n$ of order $\ln(\sqrt{d}/\varepsilon)/\rho$, i.e. $Kn=\mathcal O(\sqrt{d/\varepsilon} \ln(\sqrt{d}/\varepsilon))$. In the separable case, we get a complexity $Kn=\mathcal O(d^{1/4}/\sqrt{\varepsilon} \ln(\sqrt{d}/\varepsilon))$. Hence, we recover in this unified framework a complexity similar to  some previous works in particular cases, such as \cite{MonmarcheSplitting,BouRabeeSchuh,ChenVempala,HMC}, and better than those based on the overdamped Langevin process or first order discretization schemes of the underdamped Langevin process \cite{Dwivedi,dalalyan1,Chatterji1}. 

In the stochastic case, combining Corollary~\ref{Cor:contraction} and Propositions~\ref{prop:numerique} and \ref{prop:erreur_sto}, again,  $n = \ln(\sqrt{d}/\varepsilon)/\rho$ computations of $b$ are sufficient  to achieve $\mathcal W_2(\mathcal Law(x_n),\rpi)\leqslant \varepsilon$ but, now, with the restriction
\begin{equation}\label{eq:contrainte}
\delta^2 d + \sqrt{\delta d/p} = \mathcal O( \varepsilon)\,.
\end{equation}
Besides, here, the cost of one computation of $b(x,\theta)$ is of order $p$. Since the convergence rate $\rho$ given by Theorem~\ref{thm:crude} is of order $1/\delta$, the total complexity is  $\ln(\sqrt{d}/\varepsilon) p/\delta$, where $p,\delta$ satisfy the constraint \eqref{eq:contrainte}. We see that the best complexity we can get is $\ln(\sqrt{d}/\varepsilon)d/\varepsilon^2$, and it is achieved by taking $p $ of order $d\delta/\varepsilon^2$ and any $\delta$ with $\delta = \mathcal O(\sqrt{\varepsilon/d})$ and $  \delta^{-1} = \mathcal O(d/\varepsilon^2)$. We see that, in the stochastic case, without further refinement, there is no particular gain in using second-order numerical schemes, since the same analysis  can be done  for first-order schemes, with \eqref{eq:contrainte} replaced by $
\delta \sqrt{d} + \sqrt{\delta d/p} = \mathcal O( \varepsilon)$, so that taking $p$ of order $d\delta/\varepsilon^2$ and $\delta$ of order $\varepsilon/\sqrt{d}$, we see that $p/\delta$ is again of order $d/\varepsilon^2$.
}

\section{\new{Numerical experiments} }\label{sec:numerique}

\begin{figure}
\begin{subfigure}{0.49\textwidth}
    \centering
    \includegraphics[scale=0.38]{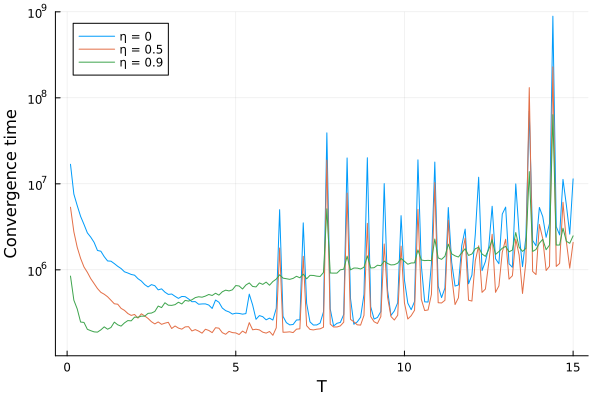}
    \caption{ugHMC ($T=K\delta$, $K$ varies)}
    \label{fig:unadjusted_T}
    \end{subfigure}
\begin{subfigure}{0.49\textwidth}
    \centering
    \includegraphics[scale=0.38]{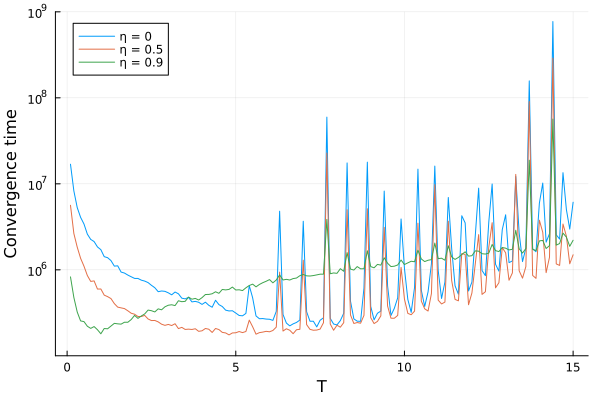}
    \caption{MagHMC ($T=K\delta$, $K$ varies)}
    \label{fig:adjusted_K}
\end{subfigure}
\begin{subfigure}{0.49\textwidth}
    \centering
    \includegraphics[scale=0.38]{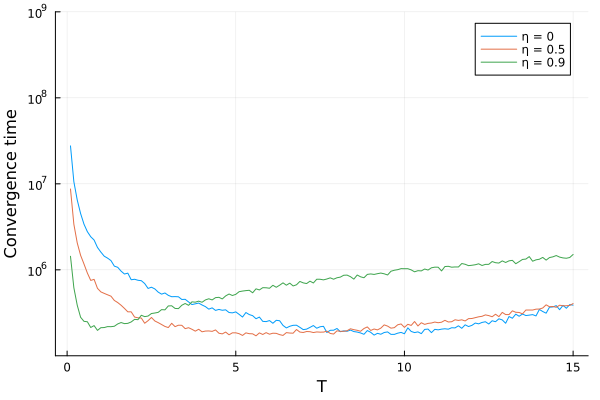}
    \caption{PDMP}
    \label{fig:PDMP-T}
  \end{subfigure} 
  \begin{subfigure}{0.49\textwidth}
    \centering
    \includegraphics[scale=0.38]{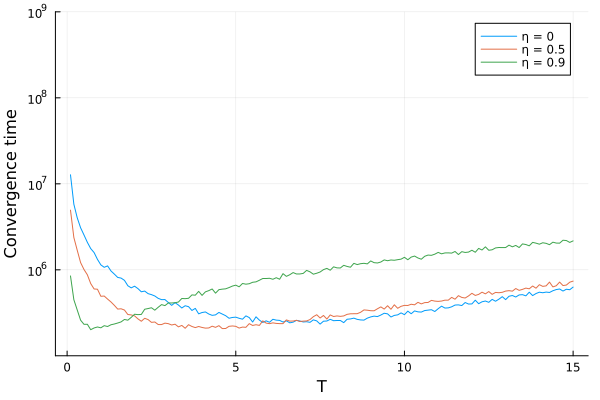}
    \caption{randomized ugHMC ($T=K\mathbb E(\tilde \delta)$, $K$ varies)}
    \label{fig:unadjusted-randomized.}
    \end{subfigure}
    \caption{Influence of the integration time $T$}\label{figure:T}
\end{figure}

\new{In this section, we study empirically the role of the velocity  \new{randomization} mechanism in the sampling efficiency. By comparison with the theoretical results,  it enables to study other samplers, namely Metropolis-adjusted gHMC (MAgHMC) chains and kinetic piecewise deterministic Markov processes  (PDMP). In both cases, the \new{randomization} step \eqref{eq:damping} with damping parameter $\eta$ is unchanged, but the unadjusted Verlet scheme \eqref{Verlet} is replaced by another transition. In the MAgHMC, it is replaced by, consecutively, first, a Metropolis-Hastings step with transition proposal given by $K$ Verlet steps with step-size $\delta$ followed by a velocity reflection and target measure $\rpi\otimes\mathrm{N}(0,I_d)$ and, second, a velocity reflection. At the end, it gives the following transition: 
\[(x_{n+1},v_{n+1}) = \left\{\begin{array}{ll}
\Phi_\delta^K (x_n,v_n) & \text{with probability } \exp \po - \po H \po \Phi_\delta^K (x_n,v_n) \pf - H(x_n,v_n)  \pf_+\pf\\
(x_n,-v_n) & \text{otherwise,}
\end{array}\right.\]
where $H(x,v)=U(x)+|v|^2/2$ stands for the total energy.

In the PDMP case, given an integration time $T>0$, the unadjusted Verlet transition is replaced by 
\[(x_{n+1},v_{n+1}) = (X_T,V_T)\]
where $(X_t,V_t)_{t\in[0,T]}$ is a trajectory of the Markov process introduced in \cite{MMZ}, which depends on a parameter $a>0$, with generator $\mathcal L$ given by
\[\mathcal L f(x,v) = v\cdot \na_x f(x,v) + \lambda_a(x,v) \po \int_{\R^d} f(x,w) q_a(x,v,\dd w) - f(x,v)\pf \]
where, writing $\xi(x) = \na U(x)/|\na U(x)|$ if $\na U(x)\neq 0$ and $\xi(x)=0$ otherwise and considering a one-dimensional standard Gaussian variable $G$,
\begin{eqnarray*}
\lambda_a(x,v) & = &  \mathbb{E}\po \po v\cdot \na U(x) + a |\na U(x)| G\pf_+\pf\\
  \int_{\R^d} f(x,w) q_a(x,v,\dd w) & = & \mathbb{E}\po f\po x, v - 2\frac{v\cdot \xi(x) + aG}{1+a^2}\xi(x)\pf  \frac{\po v\cdot \na U
  (x) + a |\na U(x)| G\pf_+}{\lambda_a(x,v)}\pf\,.
\end{eqnarray*}
This process can simulated exactly, as detailed in \cite[Section 5]{MMZ}, and it admits $\mu = \rpi\otimes\mathrm{N}(0,I_d)$ as an invariant measure. When $a=0$, we recover the so-called bouncy particle sampler \cite{PetersdeWith,Doucet2018,MonmarcheRecuitPDMP}. As $a\rightarrow \infty$, the process converges toward the Hamiltonian dynamics, see  \cite{MMZ}. In our experiments below we take $a=1$, which corresponds to the process introduced in \cite{MichelDurmusSenecal}.

The settings of the numerical experiments are the following. The dimension is $d=100$ \rev{(varying the dimension has little effect, see appendix for the results of the same experiments in dimension $d=50$)}, the target distribution is an ill-conditioned Gaussian measure corresponding to $U(x) = \sum_{k=1}^{100}\frac{|x_k|^2}{2k}$ (in particular, using the notations of Section~\ref{sec:results:Gaussian}, $L=1$ and $\kappa=1/m=100$)\rev{, the initial distribution is a centered isotropic Gaussian distribution with variance $50$ for the positions and, independently, a standard Gaussian distribution for the velocities}. To measure the efficiency, \rev{as a proxy for the Wasserstein $2$ distance, after $k$ iterations of a sampler  we use as an error function $\varphi(k)= \mathcal W_2( \mathrm{N}(0,\tilde \Sigma_k^2),\bar{\pi})$, using the formula for the Wasserstein distance between Gaussian distributions (see \cite[Proposition 7]{W2Gauss}) where $\tilde \Sigma_k^2$ is the empirical covariance matrix of $n=50$ independent realizations of the chain. Notice that the initial distribution and the target measures being centered, by the symmetries of the problem the distribution of the samplers remain centered for all times, so that we only track the convergence of the covariance matrix. Moreover, in the unadjusted case,  the sampler is a Gaussian process, so that $\varphi(k)$ is a consistent estimator (as $n\rightarrow \infty$) of the true Wasserstein distance between the distribution of the sampler at time $k$ and $\bar{\pi}$, but this is not the case for the two other samplers.}
We fix a tolerance at $\varepsilon = \rev{5.0}$ (notice that the variance of target distribution $\rpi$ is  5050, which gives a relative tolerance $\varepsilon/\mathcal W_2(\rpi,\delta_0)\simeq \rev{0.07}$),  count the number of calls to the gradient routine before the function $\varphi$ goes under $\varepsilon$ and outputs the mean cost over the \rev{$n=50$} runs. In all graphs below, this convergence time is plotted at a logarithmic scale. We study:
\begin{itemize}
    \item The influence of the integration time $T$ on the convergence time for ugHMC, MagHMC and PDMP algorithms, each for three damping parameters: $\eta \in\{ 0, 0.5,0.9\}$. In the case of gHMC algorithms, we fix the step size $\delta = 0.1$ and make $K$ vary in $\cco 1,150\ccf$, with $T=K\delta $. We take the same values $T \in [0.1,15]$ in the PDMP case. Moreover, we also consider a randomized version of ugHMC, which is exactly the same as ugHMC except that, after each velocity randomization, a random step size $\tilde \delta $ is generated, uniformly over $[0,2\delta]$, and this step size is used for the  $K$ next Verlet steps. In this last case we plot the convergence time as a function of $T=K\delta = K\mathbb E(\tilde \delta)$, where $\delta=0.1$ is fixed and $K$ varies from $1$ to $150$.    The results are  provided in Figure \ref{figure:T} \rev{(the scaling of the axes is the same in the four cases)}.
    \item The influence of $\eta$, the damping parameter, on the convergence time for ugHMC, MagHMC and PDMP algorithms, for different integration times. For ugHMC and MagHMC, we consider  $T=K\delta \in \{0.1,7.7,10\}$  (where $T=7.7$ corresponds to a periodic resonance, as can be seen in Figures~\ref{fig:unadjusted_T} and \ref{fig:adjusted_K}), where $\delta=0.1$ is fixed and $K$ varies. In the PDMP case, since there is no resonance, we take $T\in\{0.1,7,10,15\}$. In the three cases, $\eta$ goes from $0$ to $0.998$, with more values between $0.95$ and $0.998$, particularly in the case $T=0.1$, to better capture the variations in this interval.  The results are provided in Figure~\ref{figure:eta} \rev{(the scaling of the axes is the same in the three cases).}
\end{itemize}
}

\begin{figure}
\begin{subfigure}{0.49\textwidth}
    \centering
    \includegraphics[scale=0.38]{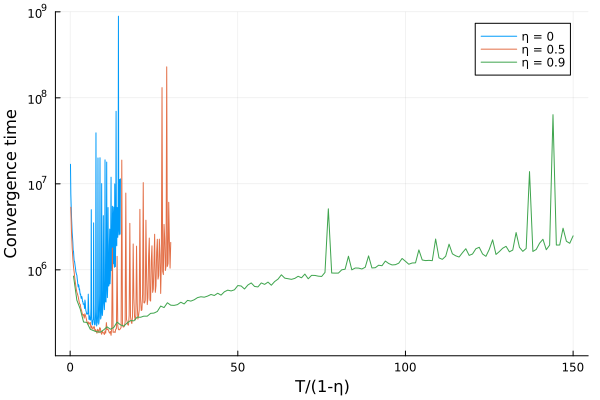}
    \caption{ugHMC}
    \label{fig:uT1-eta}
     \end{subfigure}
     \begin{subfigure}{0.49\textwidth}
    \centering
    \includegraphics[scale=0.38]{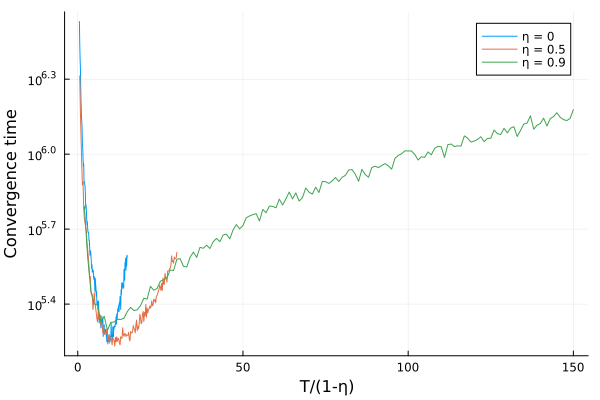}
    \caption{PDMP}
    \label{fig:PDMPT1-eta}
  \end{subfigure}
  \caption{Influence of $T/(1-\eta)$}\label{figure:Tsur1-eta}
\end{figure}

\begin{figure}
\begin{subfigure}{0.49\textwidth}
    \centering
    \includegraphics[scale=0.38]{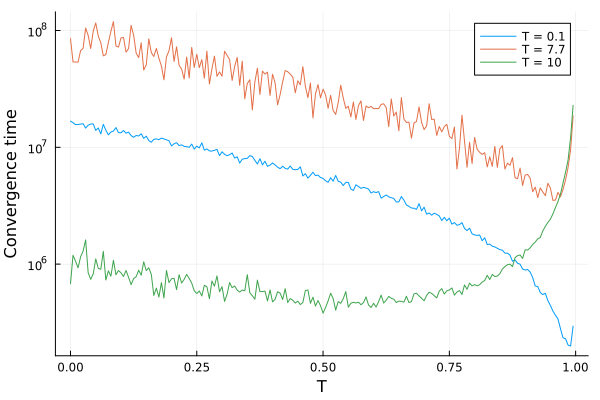}
    \caption{ugHMC}
    \label{fig:eta-unadjusted}
    \end{subfigure}
\begin{subfigure}{0.49\textwidth}
    \centering
    \includegraphics[scale=0.38]{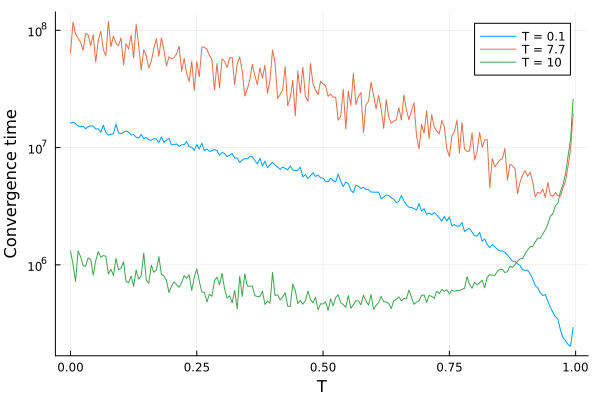}
    \caption{MagHMC}
    \label{fig:eta-adjusted}
\end{subfigure}
\centering
\begin{subfigure}{0.49\textwidth}
    \centering
    \includegraphics[scale=0.38]{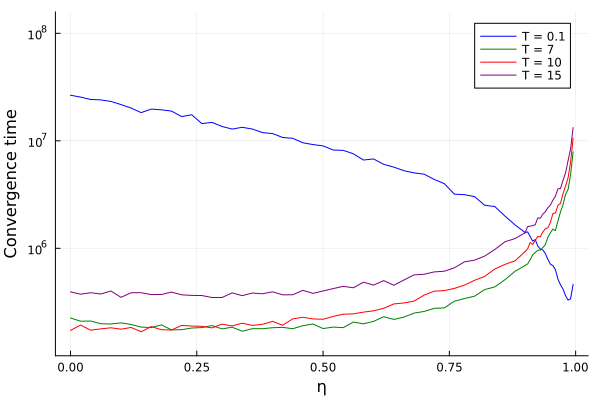}
    \caption{PDMP}
    \label{fig:eta-pdmp}
  \end{subfigure}  
  \caption{Influence of $\eta$}\label{figure:eta}
\end{figure}

\new{ 
Let us now comment these results. The cases  ugHMC and MagHMC are similar, in particular, the efficiency drops dramatically as $T$ approaches values corresponding to the period of the Hamiltonian dynamics for some of the coordinates. These values are the same for ugHMC and MagHMC and they do not depend on $\eta$. However, we see that as $\eta$ increases, the peaks are smaller, i.e. increasing $\eta$ reduces the sensibility to periodicity issues. On the contrary, for PDMP and randomized ugHMC, there is no such resonances. Putting aside these peaks, the picture is similar for the four samplers: for each value of $\eta$, for small values of $T$, the convergence time decreases until it reaches an optimal value of $T$, after which it increases. Moreover, as $\eta$ increases, the optimal value of $T$ is reached sooner (which also has an impact on the sensitivity of the optimal parameters: in the Langevin case, the optimal integration time is small and thus far from the first resonance peaks, while for $\eta=0$ the optimal time is reached between two peaks, so that a slight change of potential may completely destroy the efficiency). In fact, we can rescale the $x$-axis to better understand the phenomenon here: in Figures~\ref{fig:uT1-eta} and \ref{fig:PDMPT1-eta}, we reproduce respectively the results of  Figures~\ref{fig:unadjusted_T} and \ref{fig:PDMP-T}, but now in terms of $K\delta/(1-\eta)$ \rev{(the $y$ axes are not rescaled here for clarity, since comparing these two graphs is not relevant)}. We see that  the initial part of the curves (before reaching the optimal value of $T/(1-\eta)$) now coincides for the three values of $\eta$ (for both ugHMC and PDMP; we do not plot MagHMC and randomized ugHMC, which give similar results). The interpretation of this is that, for small values of $T$, in the absence of periodicity issues, the efficiency is mostly explained by the ability of the chain to perform long ballistic trajectories to cover, in each direction, distances which are of the order of the corresponding standard deviation. 

Concerning the dependency on $\eta$ at fixed $T$ (in Figure~\ref{figure:eta}), we observe the following.  In all cases, the convergence time decreases until an optimal $\eta$ is reached and then increases, in accordance with the theoretical analysis (e.g. Figure~\ref{fig:g}). For ugHMC and MagHMC at $T=7.7$, namely close to a resonance, the behaviour is more erratic, but this general trend still holds. In accordance to the previous discussion on the dependency in $T$, we see that the optimal $\eta$ is larger when $T$ is smaller, and goes to $1$ as $T$ goes to $0$.  }

\section*{Acknowledgments}
 P. Monmarch{\'e} would like to thank Nawaf Bou-Rabee, Andreas Eberle and Katharina Schuh for the organization of a workshop in Bonn on the long-time convergence of Markov processes and related topics which has been the occasion of many fruitful discussions.
 
 This work has received funding from the European Research Council (ERC) under the European Union’s Horizon 2020 research and innovation program (grant agreement No 810367), project EMC2, and from  the  French ANR grant SWIDIMS (ANR-20-CE40-0022).

\bibliographystyle{plain}
\bibliography{../biblio} 

\appendix

\section{Proofs in the Gaussian case}\label{sec:Gauss:proof}

This section is devoted to the proofs of the results stated in Section~\ref{sec:results:Gaussian}.

\begin{proof}[Proof of Proposition~\ref{prop:Gauss_error}]
The coordinates corresponding to eigenvectors of $S$ are independent and thus we focus first on the case $d=1$ with an eigenvalue $\lambda>0$. Assuming $\delta^2 \lambda  < 4$, a Verlet  iteration is then
\begin{equation}\label{eq:verlet_Gauss}
\Phi_V(z) = \begin{pmatrix}
1-\delta^2 \lambda/2 & \delta - \delta^3   \lambda /4 \\ -\delta\lambda   & 1- \delta^2   \lambda/2 
\end{pmatrix} z  = 
\begin{pmatrix}
\cos(\theta) & \frac1\nu \sin(\theta) \\ -\nu\sin(\theta)& \cos(\theta)
\end{pmatrix} z  
\end{equation}
with
\[\cos(\theta) = 1 - \delta^2  \lambda/2 \qquad \sin(\theta) = \delta \sqrt{\lambda   \po 1 - \delta^2  \lambda/4\pf}\qquad \nu = \sqrt{\frac{\lambda}{1-\delta^2 \lambda/4}}\,.\]
The Verlet step \eqref{eq:verlet_Gauss} is an approximation of the Hamiltonian flow with Hamiltonian $H(x,v)=(\lambda|x|^2+|v|^2)/2$, but it is the exact solution at time $\theta$ of the flow associated to the Hamiltonian  $\tilde H(x,v)=(\nu|x|^2+|v|^2/\nu)/2$. As a consequence, any probability density that is a function of $\tilde H$ is invariant by $\Phi_V^{\circ K}$. In particular, the Gaussian density proportional to $\exp\po-(\nu^2 |x|^2+|v|^2)/2\pf$ is invariant by both the Verlet and \new{randomization} steps. Using the formula for the $\mathcal W_2$ distance between Gaussian distributions (see e.g. \cite[Proposition 7]{W2Gauss}), the Wasserstein distance between the theoretical target measure $\pi_\lambda$  and the position equilibrium $\pi_{ \nu^2}$ is thus 
\[\mathcal{W}_2\po \pi_\lambda,\pi_{ \nu^2}\pf = \left|\frac{1}{\sqrt{\lambda}} - \frac{1}{\sqrt{ \nu^2}} \right|= \frac{1 - \sqrt{1-\delta^2\lambda/4}}{\sqrt{\lambda}} \,.\]
 Then, using that $x\mapsto (1-\sqrt{1-x^2})/x$ is increasing on $[0,1]$ and decomposing along the eigenspaces of $S$,
\[\varepsilon(\param) = \sup_{\lambda\in[m,L]} \frac{1 - \sqrt{1-\delta^2\lambda/4}}{\sqrt{\lambda}}  = \frac{1 - \sqrt{1-\delta^2L/4}}{\sqrt{L}}  \,. \]
\end{proof}
\begin{remark}
If the position Verlet integrator is replaced by the velocity one as in \cite{MonmarcheSplitting}, \eqref{eq:verlet_Gauss} still holds, with the same $\theta$ but with $\nu= \sqrt{\lambda(1-\delta^2 \lambda/2)}$ and thus the $\mathcal W_2$ error are equivalent for the two integrators as $\delta \rightarrow 0$. The value of $\nu$ does not intervene in  the asymptotic convergence rate of the process (see below), which means all the results of the rest of the section are true with the velocity Verlet integrator.
\end{remark}

 \begin{proposition}\label{prop:Gaussien}
Consider $Q$ the transition operator of an autoregressive Gaussian chain on $\R^{2d}$ with transitions $z'=Az+BG$, where $G$ is a standard Gaussian variable independent from $z$ and $A,B$ are matrices. Suppose that the spectrum of $A$ is included in $\{\lambda\in \mathbb{C}, |\lambda|<1\}$. Then the process has a unique invariant probability measure $\mu$, which is a  Gaussian distribution, and for all $n\in \mathbb{N}$,
\[\sup_{\nu\in\mathcal{P}_2(\R^{2d})\setminus\{\mu\}} \frac{\mathcal{W}_2(\nu Q^n,\mu)}{\mathcal{W}_2(\nu ,\mu)}=|A^n|\,.\]
Moreover, for all $p\geqslant1,n\in\mathbb N$ and probability laws $\nu,\nu'$ on $\R^{2d}$,
\begin{equation}\label{eq:WGauss}
\mathcal{W}_p(\nu Q^n,\nu' Q^n) \leqslant |A^n|\mathcal{W}_p(\nu,\nu')\,.
\end{equation}
\end{proposition}

\begin{proof}
For two probability measures $\nu,\nu'$ on $\R^{2d}$, let $\nu_c \in \Pi(\nu,\nu')$, $(Z,Z')\sim\nu_c$ and $(G_k)_{k\in\N}$ be an i.i.d. sequence of standard Gaussian variables independent from $(Z,Z')$. Set $Z_n = AZ_N+B G_n$ and $Z_n'=AZ_n'+BG_n$ for all $n\in\N$, so that $(Z_n,Z_n')$ is a coupling of $\nu Q^n$ and $\nu' Q^n$ for all $n\in\N$. Then, for all $p\geqslant 1$,
\[\mathcal{W}_p^p(\nu Q^n,\nu' Q^n)  \leqslant \mathbb E \po |Z_n -Z_n'|^p\pf = \mathbb E \po |A^n(Z_0 -Z_0')|^p\pf \leqslant |A^n|^p  \mathbb E \po |Z_0 -Z_0'|^p\pf\,.\]
Taking the infimum over $\nu_c \in\Pi(\nu,\nu')$ gives \eqref{eq:WGauss}. 

The spectral condition on $A$ implies that $|A^n|<1$ for $n$ large enough, which means that, for all $p\geqslant 1$, $Q^n$ admits a unique invariant measure within the set of probability measures with a finite $p^{th}$ moment (which is complete for the $\mathcal W_p$-distance), and thus a unique invariant measure $\mu$ in the set of probability measures with all finite moments. Since $(\mu Q)Q^n = (\mu Q^n)Q = \mu Q$ and since $\mu Q$ obviously has all finite moments if $\mu$ has, by uniqueness of the invariant measure of $Q^n$ we get that $\mu$ is invariant for $Q$.  If $\nu$ is a Gaussian distribution, then so is $\nu Q^n$ for all $n\in\N$ and thus so is $\mu$ by taking $n\rightarrow \infty$. 

 If $\mu'$ is an invariant measure for $Q$ (without assuming it has moments) considering  as before $(Z_n,Z_n')_{n\in\N}$ in the case $\nu=\mu,\nu'=\mu'$, we still get that
\[|Z_n-Z_n'| \leqslant |A^n||Z_0-Z_0'| \underset{n\rightarrow \infty}\longrightarrow 0 \qquad a.s.\]
which proves that $\mu'=\mu$.

Finally, let $\nu$ be a Gaussian distribution with the same covariance matrix as $\mu$ and a mean $z\in\R^{2d}$. Then $\nu Q^n$ has the same covariance matrix as $\mu$ for all $n$ (since it is a fixed
point of the recurrence equation satisfied by the covariance matrix of $\nu Q^n$) and a mean $A^n z$, while $\mu$ is centered. Using again the formula for the $\mathcal W_2$-distance of Gaussian distributions of \cite[Proposition 7]{W2Gauss}, we get
\[\mathcal W_2(\nu Q^n,\mu )= |A^n z| \]
for all $n\in \N$. Conclusion follows by taking the supremum of $\mathcal W_2(\nu Q^n,\mu )/\mathcal W_2(\nu ,\mu )$ over $z\in\R^{2d}$.
\end{proof}

\begin{proof}[Proof of Proposition~\ref{prop:GaussCVlongtime}]
A transition of the chain is
\begin{equation}\label{eq:transit-gaussian}
z' = \begin{pmatrix}
\cos(K\theta) & \frac1{\nu}\sin(K\theta)\eta \\ -\nu\sin(K\theta)\eta &\cos(K\theta)\eta^2
\end{pmatrix} z + \sqrt{ 1-\eta^2 }  \begin{pmatrix}
\frac1{\nu}\sin(K\theta)\eta & 0\\
\cos(K\theta)\eta & 1
\end{pmatrix}
\begin{pmatrix}  G\\ G'
\end{pmatrix} := A z+B \begin{pmatrix}  G\\ G'
\end{pmatrix}
\end{equation}
As a consequence, from Proposition~\ref{prop:Gaussien}, for any parameter $\param=(\delta,K,\eta)$,
\[\rho(\param)= \inf_{\lambda\in[m,L]} \lim_{n\rightarrow +\infty}\frac{-\ln|A^n|}{Kn} = \inf_{\lambda\in[m,L]}  \frac{-\ln \max\{|a|,a\in\sigma(A) \}}{K } \]
where $\sigma(A)$ stands for the spectrum of $A$ given by \eqref{eq:transit-gaussian}. Its eigenvalues are
\[\cos(K\theta)(1+\eta^2)/2 \pm \sqrt{\cos^2(K\theta)(1+\eta^2)^2/4-\eta^2}\,.\]
They are real iff $\cos^2(K\theta) \geqslant 4\eta^2/(1+\eta^2)^2$, and in that case the spectral radius of $A$ is
\[ |\cos(K\theta)|(1+\eta^2)/2 + \sqrt{\cos^2(K\theta)(1+\eta^2)^2/4-\eta^2}.\]
If $\cos^2(K\theta) \leqslant 4\eta^2/(1+\eta^2)^2$, the eigenvalues are complex conjugate with modulus $\eta$. In other words, in both cases, the spectral radius of $A$ is $g(|\cos(K\theta)|,\eta)$.  This concludes since $g$ is a non-decreasing function of its first variable.
\end{proof}
 
\begin{lemma}\label{lem:g}
The function 
  \[g(c,\eta) = \eta  \vee \frac{(1+\eta^2)c}2  + \sqrt{\po \po \frac{(1+\eta^2)c}2 \pf^2 - \eta^2\pf_+}\]
  is  increasing in its first variable. It is such that, for a fixed $c\in[0,1]$, $\eta\rightarrow  g(c,\eta)$ is decreasing on $[0,\eta_*(c)]$ and increasing on $[\eta_*(c),1]$ with $\eta_*(c) = (1-\sqrt{1-c^2})/c$. 
\end{lemma}
\begin{proof}
 It is clear that $g$ is an increasing function of its first variable, with $g(1,\eta)=1$ for all $\eta\in[0,1]$. Moreover, for a fixed $c\in[0,1]$, denoting by $\eta_*(c)=(1-\sqrt{1-c^2})/c$ the solution in $[0,1]$ of $2\eta=(1+\eta^2)c$, we see that $g(c,\eta)=\eta$ for $\eta \geqslant \eta_*(c)$ and, for $\eta\leqslant \eta_*(c)$, writing $y=c(1+\eta^2)/2-1/c$ (so that $y\leqslant -\sqrt{1/c^2-1} $) we get
\[g(c,\eta) = y +1/c + \sqrt{y^2 +1-1/c^2}\,, \]
which is decreasing in $y \in (-\infty,-\sqrt{1/c^2-1})$, so that $\eta\mapsto g(c,\eta)$ is decreasing on $[0,\eta_*(c)]$.
\end{proof}

 \begin{proof}[Proof of Proposition~\ref{prop:GaussienScaling}]
Notice that, as $n\rightarrow +\infty$, for any $\lambda>0$, $\arccos(1-\delta_n^2 \lambda/2) \simeq \delta_n \sqrt{\lambda}$.

 Up to extracting a subsequence, we can assume that $K_n\delta_n$ has a limit in $[0,+\infty]$ and that $\eta_n$ has a limit in $[0,1]$. Moreover, if $K_n\delta_n \rightarrow 0$ and $\eta_n\rightarrow 1$, up to extracting a subsequence we can assume that $(1-\eta_n)/K_n\delta_n$ has a limit in $[0,+\infty]$. We now distinguish all the possible cases for the values of these limits.

  \begin{itemize}
 \item if $K_n\delta_n \rightarrow +\infty$, since $m<L$, then $h(K_n,\delta_n)=1$ for $n$ large enough, so that $\rho(\param_n)= 0$.
 \item  if  $K_n\delta_n \rightarrow T$ for some $T>0$, then $h(K_n,\delta_n)\rightarrow h_*(T)$.  If, moreover, $\eta_n\rightarrow \eta<1$, then
 \[\rho(\param_n) \simeq \frac{-\ln g(h_*(T),\eta)}{K_n} \simeq  \frac{\delta_n|\ln g(h_*(T),\eta)|}{T}\]
 while, if  $\eta_n\rightarrow 1$,   $g(h(K_n,\delta_n),\eta_n)\rightarrow 1$ and thus 
 \[\rho(\param_n) \simeq \frac{-\delta_n \ln g(h(K_n,\delta_n),\eta_n)}{T} = o(\delta_n)\,. \]
 \item if $K_n\delta_n \rightarrow 0$, then $h(K_n,\delta_n) = 1 - K_n^2\delta_n^2m/2+o(K_n^2\delta_n^2)$. First, if $\eta_n \rightarrow \eta<1$, a Taylor expansion  yields 
  \[\rho(\param_n) \simeq  \frac{K_n\delta_n^2m(1+\eta^2)}{2(1-\eta^2)} =o(\delta_n)\,. \]
 Second, if $\eta_n\rightarrow 1 $, we have to distinguish  whether $(1-\eta_n)/(K_n\delta_n) $ converges to $0$, $+\infty$ or some $\gamma>0$.  Moreover, up to an additional extraction, we assume either that $\eta_n \geqslant \eta_*(h(K_n,\delta_n))$ for all $n\in\N$ or that $\eta_n \leqslant \eta_*(h(K_n,\delta_n))$ for all $n\in\N$. First, writing $x_n=K_n\delta_n$,
 \[\eta_*(h(K_n,\delta_n)) = \frac{1-\sqrt{1-\po 1 - x_n^2 m/2+o(x_n^2)\pf^2}}{1 - x_n^2 m/2+o(x_n^2)} =1-x_n \sqrt{m}+o(x_n)\,.\]
 Now, we consider four cases.
 \begin{itemize}
 \item If $(1-\eta_n)/(K_n\delta_n)  \rightarrow 0$ then $g(h(K_n,\delta_n),\eta_n) = \eta_n$ for $n$ large enough, and thus
 \[\rho(\param_n) = \frac{-\ln(\eta_n)}{K_n} \simeq  \frac{1-\eta_n}{K_n} = o(\delta_n) \]
   \item If $\eta_n \geqslant \eta_*(h(K_n,\delta_n))$ for all $n\in\N$ and $(1-\eta_n)/(K_n\delta_n) \rightarrow \gamma$ with $\gamma>0$,   then, necessarily, $\gamma\leqslant\sqrt{m}$ and
  \[\rho(\param_n) = \frac{-\ln(\eta_n)}{K_n} \simeq  \frac{1-\eta_n}{K_n}  \simeq  \gamma \delta_n\,.\]
   \item If $\eta_n \leqslant \eta_*(h(K_n,\delta_n))$ for all $n\in\N$ and $(1-\eta_n)/(K_n\delta_n) \rightarrow \gamma$ with $\gamma>0$, then, necessarily, $\gamma\geqslant\sqrt{m}$ and, writing $y_n=(1-\eta_n)/\gamma$ (so that $y_n\simeq x_n$),
   \begin{eqnarray*}
   \frac{1+\eta_n^2}{2}h(K_n,\delta_n)&=& \frac{2-2\gamma y_n+\gamma^2y_n^2}2\po 1-y_n^2m/2+o(y_n^2)\pf\\
   &=& 1 -\gamma y_n + (\gamma^2-m)y_n^2/2+o(y_n^2)
   \end{eqnarray*}
   and thus 
    \begin{eqnarray*}
   g(h(K_n,\delta_n),\eta_n) &=& 1- \gamma y_n+o(y_n) \\
   & &  + \ \sqrt{\po 1 -\gamma y_n + (\gamma^2-m)y_n^2/2 +o(y_n^2)\pf^2 - (1-\gamma y_n)^2}\\
   &=& 1-  \po\gamma  - \sqrt{\gamma^2-m} \pf y_n+o(y_n) 
   \end{eqnarray*}
so that, finally,
\[\rho(\param_n) \simeq \frac{  \po\gamma  - \sqrt{\gamma^2-m} \pf y_n}{K_n} \simeq  \po\gamma  - \sqrt{\gamma^2-m} \pf\delta_n \,. \]
\item If $(1-\eta_n)/(K_n\delta_n) \rightarrow +\infty$ then $\eta_n \leqslant \eta_*(h(K_n,\delta_n))$ for all $n$ large enough. Consider an arbitrarily large $\gamma>\sqrt m$ and, for $n\in\N$ large enough, $\tilde\eta_n = 1 - \gamma \delta_n K_n$. Using that $\eta\mapsto \rho(\delta,K,\eta)$ is decreasing on $[0,\eta_*(h(K,\delta))$, we get that $\rho(\param_n) \leqslant \rho(\delta_n,K_n,\tilde\eta_n)$ for $n$ large enough. Using the result in the previous case, $\limsup \rho(\param_n)/\delta_n \leqslant \inf\{ \gamma  - \sqrt{\gamma^2-m},\ \gamma>\sqrt m\}=0$,   in other words $\rho(\param_n)=o(\delta_n)$.
 \end{itemize}
 \end{itemize}

 \end{proof}
 
 \begin{proof}[Proof of Proposition~\ref{prop:Gaussoptim}]
 The case of the Langevin scaling is elementary. Considering HMC cases, notice that the infimum of the function $h_*$ is clearly attained in $[0,T_*]$ where $T_*=\inf\{t>0, \cos(t\sqrt{m})=-\cos(t\sqrt{L})\} = \pi/(\sqrt L + \sqrt m)$. Moreover, on this interval $h_*(T) = \cos(T\sqrt m)$, with $T\sqrt m \leqslant \pi/2$.
 
 First, consider the cases where both $\eta$ and $T$ are optimized. The minimal value of $\eta\mapsto g(h_*(T),\eta)$ is $\eta_*(h_*(T))$, obtained at $\eta=\eta_*(h_*(T))$. The question is thus to chose $T$ to maximise
\[ \frac{|\ln \eta(h_*(T)|}{T} := \frac{u(T\sqrt m)}T\,,\qquad u(t) =  \ln\cos(t) - \ln(1-\sin(t))\,. \]
On $[0,\pi/2[$, $u'(t) = 1/\cos(t)$  is increasing and thus $u$ is strictly convex. Since $u(0)=0$, we get that $t\mapsto u(t)/t$ is increasing on $[0,T_*\sqrt{m}]\subset [0,\pi/2]$. In particular, the optimal choice in the HMC scaling is $T=T_*$.  The comparison with the Langevin scaling comes from
\[\lim_{T\rightarrow 0} \frac{u(T\sqrt{m})}{T} = \sqrt{m}u'(0) = \sqrt m.\]

In the position HMC case, we need to find the maximum over $T\in[0,T_*]$ of
\[\frac{-\ln g(h_*(T),0)}T = \frac{-\ln h_*(T)}T =\frac{-\ln \cos(T\sqrt{m})}T\]
Using that again $t\mapsto -\ln \cos(t)$ is convex over $[0,\pi/2]$, we see that again the optimal choice is $T=T_*$. 

The case of a fixed $\delta>0$ is similar, in particular we know that the maximum of $\eta \mapsto \rho (\delta,K,\eta)$ is attained at $\eta=\eta_*(h(\delta,K))$, which means we are lead to show that $K\mapsto -\ln( \eta_*(h(\delta,K)))/K$ and $K\mapsto -\ln( h(\delta,K)) /K$ (the case $\eta=0)$ both reach their maximum for $K=K_*$ or $K=K_*-1$.  From the previous study of the HMC scaling, we already know that these functions are increasing over $K\in \cco 1,K_*-1\ccf$.  Indeed, notice that, for $K <\pi/\varphi_L$,  $h(K,\delta)\neq \cos(K\varphi_m)$ if and only if $-\cos(K\varphi_L) > \cos(K\varphi_m)$, i.e. $K > \pi/(\varphi_m+\varphi_L)$, which means
\[ \inf\{k\in \N, h(\delta,k)\neq \cos(k\varphi_m)\} = \min \po 1+\left \lfloor \frac{\pi}{\varphi_m+\varphi_L}\right\rfloor , \left\lceil \frac{\pi}{\varphi_L}\right\rceil \pf = 1+\left \lfloor \frac{\pi}{\varphi_m+\varphi_L}\right\rfloor = K_*\]
(where we used that $\varphi_L \leqslant \pi/2$). As a consequence, for $k<K_*$, by definition of $K_*$, $h(\delta,K)=\cos(K\varphi_m)=h_*(K\varphi_m/\sqrt m)$. Hence, it only remains to study the case $K>K_*$.

Due to the condition $2\varphi_m < \varphi_L$,
 three cases are possible:
\begin{itemize}
\item if $K\varphi_L \leqslant \pi$ then $h(\delta,K)=\max(|\cos(K\varphi_L)|,\cos(K\varphi_m))$.
\item if $ \pi \in [K\varphi_m ,K\varphi_L ]$ then $h(\delta,K)=1$.
\item if $K\varphi_m \geqslant \pi$ then the length of $[K\varphi_m,K\varphi_L]$ is larger than $\pi$ so, again, $h(\delta,K)=1$.
\end{itemize}
As a consequence, under the additional assumptions, the behavior of $K\mapsto h(\delta,K)$ is the following: for $K< K_* $, $h(\delta,K)=\cos(K\varphi_m)$ is decreasing in $K$; from $K_*$ to $\lceil \pi/\varphi_L\rceil-1$ (this set being potentially empty) $h(\delta,K)=-\cos(K\varphi_L)$ is increasing in $K$; for $K\geqslant  \lceil \pi/\varphi_L\rceil$, $h(\delta,K)=1$. In particular, the minimum of $K\mapsto h(\delta,K)$ is attained on $\cco 1,K_*\ccf$ and, a fortiori, so are the maxima of $K\mapsto \rho(\delta,K,\eta_*(h(\delta,K)))$ and of $K\mapsto \rho(\delta,K,0)$, which concludes.

 \end{proof}

 \section{Proofs in the general case}\label{sec:proofDimFreeCV}

This section is devoted to the proofs of the results stated in Section~\ref{sec:results:Wasserstein}.

 \new{
\begin{remark}\label{rem:rescaling_supp}
 For simplicity, whenever Assumption~\ref{Assu:mLb} is assumed, by rescaling,  we will suppose that $L=1$. Indeed, as in Remark~\ref{rem:rescaling}, if $Z=(X,V)$ is an SG\new{g}HMC with parameters $\delta,K,\eta$ associated with $b$ satisfying Assumption~\ref{Assu:mLb} for some $m,L$, then  $Z'=(\sqrt{L} X,V)$ is an SG\new{g}HMC with parameters $\delta'=\delta\sqrt{L},K,\eta$ and associated with $\hat b(x,\theta) =  b(  x/\sqrt{L},\theta)/\sqrt{L}$, which satisfies Assumption~\ref{Assu:mLb} with $\hat m=m/L$ and $\hat L=1$. Moreover, if $x\mapsto \na_x b(x,\theta)$ is $L_2$-Lipschitz (as in Proposition~\ref{prop:numerique}) then $x\mapsto \na_x \hat b(x,\theta)$ is $L_2/L^{3/2}$-Lipschitz. Hence, once a result is proven with $L=1$, the same result in the general case is obtained by replacing $\delta$ by $\delta \sqrt{L}$, $m$ by $m/L$, $L_2$ by $L_2/L^{3/2}$ and the standard $\mathcal W_2$ distance by the $\mathcal W_2$ distance associated to the distance $\sqrt{L|x-x'|^2 +|v-v'|^2}$, which is larger than $\min(L,1)$ times the former.
\end{remark} 
}

 \subsection{Parallel coupling}\label{sec:preuveParallel}

Let $\mathcal M(m,L)$ be the set of $d\times d$ matrices $A$ with   $m|u|^2 \leqslant u\cdot A u$  and $|Au|\leqslant L|u|$ for all $u\in\R^d$.

\begin{lemma}\label{lem:VerletMatrice}
Under Assumption~\ref{Assu:mLb}, suppose that $2 \delta \leqslant 1 = L $.  Let $z_0,z_0'\in\R^{2d}$ and $\theta_1,\dots,\theta_K\in\Theta$. For $k\in\cco1,K\ccf$, set
\[z_{k}=\Phi_\delta^{\theta_k}(z_{k-1})\,,\qquad z_{k}'=\Phi_\delta^{\theta_k}(z_{k-1}')\,.\]
Then, there exist   $\alpha,\beta,\gamma \in \mathcal M(m,1)$ such that 
\[\left | z_{K}-z_{K}' - \begin{pmatrix}
1- (\delta K)^2  \alpha/2   & K \delta  \\ -\delta K  \beta &   1- (\delta K)^2   \gamma/2
\end{pmatrix} (z_0-z_0')\right|\leqslant C_{K,\delta}    |z_0-z_0'|\]
where 
\begin{equation}\label{eq:CK}
C_{K,\delta} = \frac14 e^{\new{3(K-1)\delta /2}}\delta^3 \po K^3 + K-1\pf
\end{equation}
\end{lemma}
\begin{proof}
Denote $\Delta z_n=(\Delta x_n,\Delta v_n) = z_n-z_n'$ for $n\in\cco 0,K\ccf$. For all $n\in\cco 0,K-1\ccf$, there exists $y_n\in\R^d$ such that $b(x_n+\delta/2 v_n,\theta_{n+1})-b(x_n'+\delta/2 v_n',\theta_{n+1})=\na_x b(y_n,\theta_{n+1})(\Delta x_n+\delta/2\Delta v_n)$. Denote $B_n = \na b(y_n,\theta_{n+1})$. Fix $k\in\cco 0,K-1\ccf $. Then 
\begin{equation}
\label{eq:rev}
\Delta z_{k+1} = \begin{pmatrix}
1 -\delta^2/2 B_k & \delta - \delta^3/4B_k\\
-\delta B_k  & 1 - \delta^2/2 B_k
\end{pmatrix}
\Delta z_k \,, 
\end{equation}
and in particular
\[\left| \Delta z_{k+1} - \begin{pmatrix}
1 -\delta^2/2 B_k & \delta \\
-\delta B_k  & 1 - \delta^2/2 B_k
\end{pmatrix}
\Delta z_k \right| \leqslant \frac{\delta^3}4 |\Delta z_k| \]
Suppose by induction that there exist $C_k\geqslant 0$ and matrices $\alpha_k,\beta_k,\gamma_k \in \mathcal M(m,1)$ such that
\begin{align*}
\left|\Delta z_{k} - 
\begin{pmatrix}
1- (\delta k)^2 \alpha_k/2 & \delta k\\ -\delta k \beta_k & 1 -(\delta k)^2 \gamma_k/2
\end{pmatrix}   \Delta z_0\right| \leqslant C_k|\Delta z_0|\,,
\end{align*}
which is clearly true for $k=0$ with $C_0=0$ and, according to the previous computation, is also true for $k=1$ with $C_1= \delta^3 /4$. 
In particular, 
\[|\Delta z_k| \leqslant \po 1 +   \delta k +(\delta k)^2/\new{2} +C_k\pf |\Delta z_0|\,, \]
where we used that  for any $\alpha\in\mathcal M(m,1)$ and $w\in\R^d$,
\begin{align*}
|(1- (\delta k)^2 \alpha/2)w|^2 &\leqslant  \po 1 - (\delta k)^2 m+  (\delta k)^4/4\pf |w|^2 \leqslant \po 1  +  (\delta k)^4/4\pf|w|^2\,,
\end{align*}
\new{and then that $\sqrt{1+x} \leqslant 1+\sqrt{x}$ for $x>0$.}  We decompose
\begin{multline*}
\begin{pmatrix}
1 -\delta^2 B_k/2 & \delta  \\
-\delta B_k  & 1 - \delta^2 B_k/2
\end{pmatrix}\begin{pmatrix}
1- (\delta k)^2 \alpha_k/2 & \delta k\\ -\delta k \beta_k & 1 -(\delta k)^2 \gamma_k/2
\end{pmatrix}\\
 = \begin{pmatrix}
1- \delta^2 (k+1)^2 \alpha_{k+1}/2 & \delta (k+1)\\ -\delta (k+1) \beta_{k+1} & 1 -\delta^2 (k+1)^2 \gamma_{k+1}/2
\end{pmatrix} + R
\end{multline*}
with
\[\alpha_{k+1}= \frac1{(k+1)^2}\po B_k + k^2 \alpha_k + 2k\beta_k\pf \qquad \beta_{k+1}= \frac1{k+1}\po B_k+ k\beta_k \pf \]
\[\gamma_{k+1}= \frac1{(k+1)^2}\po (2k+1) B_{k} + k^2 \gamma_k \pf\]
which are all three in $\mathcal M(m,1)$ and
\[R=\begin{pmatrix}
\delta^4 k^2 B_k  \alpha_{k}/4 & -\delta^3 k\po B_k+k \gamma_k\pf/2\\ \delta^3k^2 B_{k} \alpha_k/2     +\delta^3k B_{k}\beta_k /2 & \delta^4 k^2 B_{k} \gamma_k/4
\end{pmatrix}\,, \]
which can be bounded as
\[|R| \leqslant  \delta^4 k^2    /4 + \delta^3 k(k+1)  /2\,.\]
We get that
\begin{align*}
&\left|\Delta z_{k+1} - \begin{pmatrix}
1- \delta^2 (k+1)^2 \alpha_{k+1}/2 & \delta (k+1)\\ -\delta (k+1) \beta_{k+1} & 1 -\delta^2 (k+1)^2 \gamma_{k+1}/2
\end{pmatrix}\Delta z_0\right|\\
 \leqslant& |R\Delta z_0|+ \left|\Delta z_{k+1} - \begin{pmatrix}
1 -\delta^2/2 B_k & \delta \\
-\delta B_k  & 1 - \delta^2/2 B_k
\end{pmatrix}  \Delta z_k\right|\\
&+ \left|  \begin{pmatrix}
1 -\delta^2/2 B_k & \delta \\
-\delta B_k  & 1 - \delta^2/2 B_k
\end{pmatrix} \po\Delta z_k- \begin{pmatrix}
1- (\delta k)^2 \alpha_k/2 & \delta k\\ -\delta k \beta_k & 1 -(\delta k)^2 \gamma_k/2
\end{pmatrix}   \Delta z_0\pf\right| \\
\leqslant & \po\delta^4 k^2    /4 + \delta^3 k(k+1)  /2\pf|\Delta z_0|+\delta^3 |\Delta z_k|/4 + (1+\delta +\delta^2/\new{2})C_k|\Delta z_0|\\
\leqslant &  C_{k+1} |\Delta z_0|
\end{align*}
with, using that $\delta  \leqslant 1/2$,
\begin{align*}
C_{k+1} &= \delta^4 k^2    /4 + \delta^3 k(k+1)  /2+\delta^3 \po 1 +   \delta k +(\delta k)^2/\new{2} +C_k\pf/4 + (1+\delta +\delta^2/\new{2})C_k\\
&\leqslant \delta^3    (3k^2+3k+1)/4 + (1+\new{3\delta /2})C_k
\end{align*}  
From this,
\[C_K \leqslant \ \frac{\delta^3   }4  \sum_{k=0}^{K-1} (3k^2+3k+1)(1+\new{3\delta /2})^{K-1-k} \ \leqslant \ \frac{\delta^3   }4  e^{\new{3\delta (K-1)/2}} \po K^3 + K -1\pf\,.\]
\end{proof}

 \begin{lemma}\label{lem:VerletMatrice+damping}
 Let $(z_k,z_k')_{k\in\N}$ be a parallel coupling of two SG\new{g}HMC chains. Under Assumption~\ref{Assu:mLb}, suppose that $2\delta  \leqslant 1 = L$.  Then, there exist   $\alpha,\beta,\gamma \in \mathcal M(m,1)$ such that
 \[\left | z_{1}-z_{1}' - \begin{pmatrix}
1- (\delta K)^2  \alpha/2   & K \delta \eta \\ -\delta K\eta \beta &  \po 1- (\delta K)^2   \gamma/2\pf \eta^2
\end{pmatrix} (z_0-z_0')\right|\leqslant C_{K,\delta}    |z_0-z_0'|\,,\]
where $C_{K,\delta}$ is given by \eqref{eq:CK}.
 \end{lemma}
 \begin{proof}
 Let $G$ be the Gaussian variable used in the first \new{randomization} step, and let $z_{1/3}=(x_{1/3},v_{1/3})$ with
 \[x_{1/3}=x_0\,,\qquad v_{1/3} = \eta v_0 + (1-\eta^2) G \]
 and similarly for $z'$. Let $z_{2/3}$ and $z_{2/3}'$ be the states of the chains after the $K$ stochastic Verlet steps. The effect of the synchronous \new{randomization} is thus that
 \begin{equation*}
z_{1/3}-z_{1/3}' = \begin{pmatrix}
 1 & 0 \\ 0 & \eta
 \end{pmatrix}(z_0-z_0')\,,\qquad z_{1}-z_{1}' = \begin{pmatrix}
 1 & 0 \\ 0 & \eta
 \end{pmatrix}(z_{2/3}-z_{2/3}')\,. 
 \end{equation*}
 In particular, $|z_{1/3}-z_{1/3}'|\leqslant |z_0-z_0'|$ and, for any matrix $B$,
  \[\left | z_{1}-z_{1}' - \begin{pmatrix}
1   & 0 \\ 0 & \eta 
\end{pmatrix} B (z_{1/3}-z_{1/3}')\right|\leqslant  | z_{2/3}-z_{2/3}' - B(z_{1/3}-z_{1/3}')|\,. \]
Applying Lemma~\ref{lem:VerletMatrice} with initial conditions $(z_{1/3},z_{1/3}')$ concludes.
 \end{proof}

The main step of the proof of Theorem~\ref{thm:crude} is the following.

\begin{proposition}\label{prop:contraction}
Under Assumption~\ref{Assu:mLb}, suppose that $2\delta   \leqslant 1 = L$. Let $(z_n,z_n')_{n\in\N}$ be a parallel coupling of two SG\new{g}HMC  chains.
\begin{enumerate}
\item If $\eta=0$\rev{, $K\delta\leqslant 1/2$  and $\delta  \leqslant mK\delta /31$} then, denoting by $x_n$ and $x_n'$ the first $d$-dimensional coordinates of $z_n$ and $z_n'$, almost surely, for all $n\in\N$,
\begin{equation}
\label{eq:rev19}
|x_n-x_n'| \leqslant \rev{\po 1 - \frac{m}{20}(\delta K)^2 \pf^n} |x_0-x_0'|\,.
\end{equation}
\item For any $\eta\in[0,1)$,  set $a=m/2$, $c=K\delta \eta/(1-\eta^2)$, 
\begin{equation}\label{eq:MrhoProp}
    M=\begin{pmatrix}
1  & c \\ c & a
\end{pmatrix} \,\quad
\rho  =   \min\po \frac{1-\eta^2}{\new{5}}, \frac{ m}{40}(K\delta)^2 \po \frac{2\new{\eta^2}}{1-\eta^2} +  1-\eta^2 \pf \pf  \,.
\end{equation}
and let $C_{K,\delta}$ be  given by \eqref{eq:CK}. Assume that
\begin{equation}\label{eq:condition}
(K\delta)^2 \co \frac{\new{6}    \eta^2}{m(1-\eta^2)}    + \frac{\new{2}\eta }{\new{5}} +   \frac{1 }{m^2}  \cf  \leqslant \frac{1-\eta^2}{\new{5}}\qquad\text{and}\qquad \frac{2}{\sqrt m} C_{K,\delta}\leqslant \rho\,.
\end{equation}
Then, for all $n\in\N$,
\[\|z_n-z_n'\|_M \leqslant \po 1 -\rho\pf^n \|z_0-z_0'\|_M\,.\]
Moreover, for all $z=(x,v)\in\R^{2d}$,
  \begin{equation}\label{eq:norme_equivalence}
  \frac{2}{3}\po |x|^2 +  2m |v|^2\pf \leqslant  \|z\|_M^2 \leqslant \frac{4}{3}\po |x|^2 +  2m |v|^2\pf 
  \end{equation}
  \item \new{For $\eta>0$, assuming that $\overline{\gamma}:=(1-\eta^2)/(K\delta\eta)\geqslant 2$ and that $K\delta \leqslant \eta m /[ \overline{\gamma}(\modifplb{21}+\modifplb{9}\overline{\gamma}^2)]$, then, for all $n\in\N$,
\[\|z_n-z_n'\|_M \leqslant \po  1 - \frac{\eta m}{\modifplb{6} \overline{\gamma}} K\delta  \pf^n \|z_0-z_0'\|_M\,,\]
where
\[M = \begin{pmatrix}
1 & \overline{\gamma}^{-1} \\ \overline{\gamma}^{-1} & 1 
\end{pmatrix}\in\mathcal M(1/2,3/2)\,.\]}
\end{enumerate}

\end{proposition}

\begin{proof}
\noindent \textbf{Full resampling.} Assume $\eta=0$.  \rev{In this case, using only the bound from Lemma~\ref{lem:VerletMatrice+damping} is not enough to recover the sharp condition on the integration time of \cite{ChenVempala}, and thus we adapt the arguments of Chen and Vempala to the discrete time setting. From Lemma~\ref{lem:VerletMatrice+damping} (with $\Delta v_0=0$), using that $(k^3+k-1)\leqslant 5k^3/4$ for all $k\geqslant 1$, 
\[
|\Delta x_1 - \Delta x_0 |   \leqslant   \frac12 |  (K\delta)^2 \alpha \Delta x_0|  + \frac{5}{16} e^{\new{3K\delta /2}}  (K\delta)^3     |\Delta x_0 | \leqslant \frac14 |\Delta x_0 | 
\]
using that $K\delta\leqslant 1/2$. Denoting $(x_{(k)},v_{(k)})$ and $(x_{(k)}',v_{(k)}')$ the intermediary Verlet steps of the two chains for $k\in\cco 0,K\ccf$ (so that $x_1=x_{(K)}$ and $x_1'=x_{(K)}'$) and $\Delta x_{(k)} = x_{(k)}-x_{(k)}$ (and similarly for the velocities), the previous bound holds when $K$ is replaced by $k\leqslant K$, from which
\begin{equation}
\label{eq:rough}
\frac34 |\Delta x_{0}| \leqslant |\Delta x_{(k)}| \leqslant \frac{5}{4} |\Delta x_{0}|
\end{equation}
for all $k\leqslant K$.

Denote 
\[\rho(k) = \frac{\Delta x_{(k)}\cdot B_k \Delta x_{(k)} }{|\Delta x_{(k)}|^2} \in [m,1]\,.\]
By diagonalizing $B_k\in \mathcal M(m,1)$ and then using \eqref{eq:rough}, we see that 
\begin{equation}
\label{eq:coconv}
|B_k \Delta x_{(k)}|^2 \leqslant \Delta x_{(k)} \cdot B_k \Delta x_{(k)}  = \rho(k) |\Delta x_{(k)}|^2 \leqslant \frac{25}{16} \rho(k) |\Delta x_0|^2\,. 
\end{equation}
Using this and \eqref{eq:rev},
\[|\Delta v_{(k+1)}| = |\po 1-\delta^2 B_k/2\pf  \Delta v_{(k)} - \delta B_k \Delta  x_{(k)}| \leqslant (1+\delta^2/2) |\Delta v_{k}| + \delta \frac54  \sqrt{ \rho(k)} |\Delta x_0|  \,,  \]
and by induction and then Cauchy-Schwarz, 
\begin{equation}
\label{eq:rev3}
|\Delta v_{(k)}|^2  \leqslant \co \delta \frac54  e^{k\delta^2/2} |\Delta x_0| \sum_{j=0}^{k-1} \sqrt{\rho(j)} \cf^2 \leqslant \frac{25}{16} k\delta^2 e^{k\delta^2} |\Delta x_0|^2 \sum_{j=0}^{k-1} \rho(j)  \leqslant 2 k\delta^2  |\Delta x_0|^2 Q(k)\,,  
\end{equation}
where, in the last inequality, we introduced $Q(k)= \sum_{j=0}^{k-1} \rho(j)$ and used that $\delta\leqslant 1/62$ and  $k\delta \leqslant 1/2$ for $k\leqslant K$ so that  $25 e^{k\delta^2}/16\leqslant 2 $. In particular, using that $k\delta^2 Q(k) \leqslant k^2\delta^2 \leqslant 1/4$, we get 
\begin{equation}
\label{eq:rough_v}
|\Delta v_{(k)}|^2   \leqslant \frac12   |\Delta x_0|^2 \,.     
\end{equation} 

Now, from \eqref{eq:rev} again, sorting by powers of $\delta$ and bounding rougly terms of order larger than $2$ with $|B_k|\leqslant 1$ and either \eqref{eq:rough} or \eqref{eq:rough_v} (or ignoring them when they are negative),
\begin{eqnarray}
\lefteqn{\Delta x_{(k+1)} \cdot \Delta v_{(k+1)}}\nonumber \\
  &=& \po (1-\delta^2 B_k/2) \Delta x_{(k)} + (\delta - \delta^3 B_k/4) \Delta v_{(k)} \pf \cdot \po -\delta B_k \Delta x_{(k)} + (1-\delta^2 B_k/2) \Delta v_{(k)}\pf \nonumber \\
  & \leqslant & \Delta x_{(k)} \cdot \Delta v_{(k)} - \delta \rho(k) |\Delta x_{(k)}|^2 + \delta |\Delta v_{(k)}|^2 \nonumber \\
   & & + \frac{\delta^3}2 | \Delta x_{(k)}|^2 +  \po \frac{3\delta^2}2  + \frac{3\delta^4}{4} \pf  | \Delta x_{(k)}| | \Delta v_{(k)}|
    +   \frac{\delta^3}4(1+\delta^2/2)  | \Delta v_{(k)}|^2 \nonumber \\
      &    \leqslant &  \Delta x_{(k)} \cdot \Delta v_{(k)} + \co - \frac{9}{16}\delta \rho(k)  +   \frac53 k\delta^3    Q(k)  + \frac32 \delta^2 \cf  |\Delta x_0|^2\,\nonumber
\end{eqnarray}
where in the last line we used the lower-bound of \eqref{eq:rough} (using that $\rho\geqslant 0$) and \eqref{eq:rev3}. By induction,
\begin{eqnarray}
\Delta x_{(k)} \cdot \Delta v_{(k)} & \leqslant &  |\Delta x_0|^2 \sum_{j=0}^{k-1} \co - \frac{9}{16}\delta \rho(j)  +   \frac53 j\delta^3    Q(j)  + \frac32 \delta^2 \cf \nonumber \\
 & = & |\Delta x_0|^2  \po - \frac{9}{16}\delta Q(k) +  \frac32k  \delta^2+   \frac53  \delta^3  \sum_{j=0}^{k-1} j   Q(j) \pf \label{eq:rev5}
\end{eqnarray}

Finally, from \eqref{eq:rev} again, treating higher order terms as before,
\begin{eqnarray*}
|\Delta x_{(k+1)}|^2 & =& |(1-\delta^2 B_k/2) \Delta x_{(k)} + (\delta - \delta^3 B_k/4) \Delta v_{(k)}|^2 \\
& \leqslant  & |\Delta x_{(k)}|^2 + \delta \Delta x_{(k)} \cdot \Delta v_{(k)} + \delta^2 |\Delta v_{(k)}|^2   \\
& &  + 2\po |\Delta x_{(k)} | + \delta |\Delta v_{(k)} |\pf \po \delta^2/2 |\Delta x_{(k)} | + \delta^3/4 |\Delta v_{(k)} |\pf   + \po \delta^2/2 |\Delta x_{(k)} | + \delta^3/4 |\Delta v_{(k)} |\pf^2  \\
& \leqslant & |\Delta x_{(k)}|^2 + \delta \Delta x_{(k)} \cdot \Delta v_{(k)} + \frac52 \delta^2 |\Delta x_0|^2 \,. 
\end{eqnarray*}
By induction and then using \eqref{eq:rev5},  
\begin{eqnarray*}
|\Delta x_{(K)}|^2
& \leqslant & |\Delta x_{0}|^2+ \delta \sum_{j=0}^{K-1} \Delta x_{(k)} \cdot \Delta v_{(k)} + \frac52 \delta^2 K |\Delta x_0|^2   \\
& \leqslant & \co 1 + \frac{5}{2} K\delta^2 + \delta  \sum_{k=0}^{K-1} \po - \frac{9}{16}\delta Q(k) +  \frac32k  \delta^2+   \frac53  \delta^3  \sum_{j=0}^{k-1} j   Q(j) \pf \cf |\Delta x_{0}|^2  \\
& \leqslant & \co 1 + \frac{5}{2} K\delta^2 + \frac34 K^2   \delta^3 +  \delta  \sum_{k=0}^{K-1} \po - \frac{9}{16}\delta Q(k) +   \frac53  \delta^3 Q(k) \sum_{j=0}^{k-1} j    \pf \cf |\Delta x_{0}|^2  \,,
\end{eqnarray*}
where we used that $Q$ is non-decreasing in the last line. Now, for $k\leqslant K$, using that $Q(k) \geqslant km\geqslant 0$, 
\[ Q(k) \po - \frac{9}{16}  +   \frac53  \delta^2   \sum_{j=0}^{k-1} j \pf  \leqslant   Q(k) \po - \frac{9}{16}  +   \frac56  (K\delta)^2    \pf   \leqslant - \frac{2km}{5}\,,   \] 
from which finally 
\begin{eqnarray*}
|\Delta x_{(K)}|^2 &\leqslant& \co 1 + \frac{5}{2} K\delta^2 + \frac34 K^2   \delta^3 -  \delta^2  \sum_{k=0}^{K-1}  \frac{2km}{5}  \cf |\Delta x_{0}|^2  \\
&  = & \co 1 + \po  \frac{5}{2} + \frac{m}{5}\pf  K\delta^2 + \frac34 K^2   \delta^3  -  (K\delta)^2   \frac{m}{5}   \cf |\Delta x_{0}|^2 \\
& \leqslant & \co 1 + \po  \frac{5}{2} + \frac{1}{5} + \frac{3}{40} \pf  K\delta^2   -  (K\delta)^2   \frac{m}{5}   \cf |\Delta x_{0}|^2 \,,
\end{eqnarray*}
where we used that $K\delta \leqslant 1/2$ and $m\leqslant 1$.
Using that $\delta \leqslant m K\delta /30$ concludes the proof of~\eqref{eq:rev19}, since it gives 
\[|\Delta x_{1}| = |\Delta x_{(K)}| \leqslant \sqrt{  1     -  (K\delta)^2   \frac{m}{10}}  |\Delta x_{0}|\leqslant  \po 1     -  (K\delta)^2   \frac{m}{20} \pf |\Delta x_{0}|\] 
}

\noindent \textbf{General case.} Let $\eta\in[0,1)$. 
First, notice that the condition \eqref{eq:condition} and $m\leqslant 1$ imply that
\begin{equation}\label{eq:bornescontraction}
5K\delta \leqslant m\,,\quad (5c)^2 \leqslant 2a\,,\quad c\leqslant \frac15\quad \text{and}\quad \frac23a \leqslant \frac23\begin{pmatrix}
 1 & 0 \\ 0 & a
\end{pmatrix}  \leqslant M \leqslant \frac43 \begin{pmatrix}
 1 & 0 \\ 0 & a
\end{pmatrix} \leqslant \frac43\,.
\end{equation}
In particular, the last part is exactly \eqref{eq:norme_equivalence}.

From Lemma~\ref{lem:VerletMatrice+damping}, there exist $\alpha,\beta,\gamma\in\mathcal M(m,1)$ such that, denoting
\[ A = \begin{pmatrix}
1- (\delta K)^2  \alpha/2   & K \delta \eta \\ -\delta K\eta \beta &  \po 1- (\delta K)^2   \gamma/2\pf \eta^2
\end{pmatrix} \qquad \tilde z = \Delta z_1 - A\Delta z_0\,,    \]
we have $|\tilde z|\leqslant C_{K,\delta}|\Delta z_0|$.  Remark that, for small values of $K\delta$, $C_{K,\delta}$ is of order $(K\delta)^3$. Hence, in the computation of $\|A\Delta z_0\|_M$, it is not useful to keep track of the terms of order larger than $(K\delta)^2$ and we bound them roughly. \new{In the following we write $s=K\delta$.} We end up with
\[A^T M A = M + \begin{pmatrix}
B_{11} & B_{12} \\ B_{12}^T & B_{22}
\end{pmatrix} + R\]
with, organising by powers of $K\delta$,
\begin{eqnarray*}
B_{11} &=& - \new{s \eta } c(\beta+\beta^T) -\new{s}^2\co  \frac{\alpha+\alpha^T}{2} - a\eta^2\beta^T \beta  \cf\\
B_{12} & = & c(\eta^2-1) + \new{s} (\eta-a\eta^3 \beta^T) - \new{s^2 c\eta^2 \co \beta^T +  \frac{\gamma+\alpha^T}{2}\cf } \\
B_{22} &= & a (\eta^4-1) + 2  \new{s} c\eta^3 + \new{s}^2 \co \eta^2 - a\eta^4 \frac{\gamma+\gamma^T}{2}\cf
\end{eqnarray*}
\[R = \begin{pmatrix}
\new{s}^3 c\eta \frac{\alpha^T \beta + \beta^T \alpha}{2} + \frac14\new{s}^4 \alpha^T \alpha  
& \quad & \frac12\new{s}^3 \po   a \eta^3 \beta^T\gamma - \eta\alpha^T\pf + \frac14 \new{s}^4c\eta^2 \alpha^T\gamma   \\
\frac12\new{s}^3 \po   a \eta^3  \gamma^T \beta - \eta\alpha\pf + \frac14 \new{s}^4c\eta^2 \gamma^T  \alpha
& \quad  &
-\new{s}^3 c\eta^3 \frac{\gamma+\gamma^T}{2} + \frac14 \new{s}^4 a\eta^4 \gamma^T \gamma 
\end{pmatrix}.\]
The choice $c=\new{s} \eta/(1-\eta^2)$ cancels out the leading term of $B_{12}$. The remaining terms of $B_{12}$ are bound as follows:
\[\begin{pmatrix}
0 & -\new{s} a\eta^3 \beta^T  \\ -\new{s}  a\eta^3 \beta  & 0 
\end{pmatrix} \leqslant  \new{s} a\eta^3  \begin{pmatrix}
\theta & 0 \\ 0 & \frac{1}{\theta}
\end{pmatrix}
\]
for any $\theta>0$, and similarly for the term in $\new{s}^2$ with a parameter $\theta'>0$. We get
\[\begin{pmatrix}
B_{11} & B_{12} \\ B_{12}^T & B_{22}
\end{pmatrix} \leqslant \begin{pmatrix}
C_1 & 0 \\ 0 & a C_2
\end{pmatrix}\]
with
\begin{eqnarray*}
C_1 &=& -  \new{s} \co 2 \new{\eta} cm -  a\eta^3  \theta\cf -\new{s}^2\co  m - a\eta^2  - \new{2} c\eta^2 \theta'  \cf\\
aC_2 &=& a (\eta^4-1) +   \new{s} \po 2c\eta^3 +  a\eta^3  /\theta\pf + \new{s}^2 \co \eta^2 - a\eta^4m  + \new{2 c\eta^2}/\theta'\cf\,.
\end{eqnarray*}
In order to have a negative leading term in $C_1$ we take $\theta = cm/(a\eta^\new{2}  )$, and for simplicity in the term in $\new{s}^2$ of $C_1$ we take $\new{4 c\eta^2} \theta' = m$, which yields
\[C_1 = -  \new{s\eta } cm  -\new{s}^2\co \frac12  m - a\eta^2   \cf = -  \new{s\eta }cm -\frac12 \new{s}^2 m(1-\eta^2) = -\new{s}^2 m \po \frac{\new{\eta^2} }{1-\eta^2} + \frac{1-\eta^2}{2}\pf \]
with the choice $a=m/2$. These choices yield
\begin{eqnarray*}
C_2 &=&    \eta^4-1 +   \new{s}\po \frac{4  \new{s}  \eta^4}{m(1-\eta^2)} +    \frac{\new{\eta^4}(1-\eta^2)}{2\new{s}}\pf + \new{s}^2 \co \frac{2 \eta^2}m - \eta^4m  + \po \frac{\new{4  c\eta^2}}m\pf^2 \cf \\
 &=&    -(1-\eta^2) \po 1 + \eta^2 -  \frac12 \new{\eta^4} \pf  +     \new{s}^2 \co \frac{4     \eta^4}{m(1-\eta^2)} + \frac{2 \eta^2}m - \eta^4m  + \po \frac{\new{4 c\eta^2}}m\pf^2 \cf \\
 &\leqslant &  -   \po 1 - \eta^2\pf   +     \new{s}^2 \co \frac{6    \eta^2}{m(1-\eta^2)}     +  \frac{ 1  }{\new{m^2}} \cf \,,
\end{eqnarray*}
where in the last term we used that $c \leqslant 1/5$.  Now we bound the higher order terms in $\new{s}$ as
\[R \leqslant \begin{pmatrix}
D_1  
& 0   \\
0 & aD_2
\end{pmatrix}\]
with, using  that $c\leqslant 1/5$, $a\leqslant 1/2$ and $5\new{s} \leqslant 1$, 
\begin{eqnarray*}
D_1 &=& \new{s}^3 c\eta + \frac14\new{s}^4   + \frac12\new{s}^3 \po   a \eta^3  + \eta \pf + \frac14 \new{s}^4c\eta^2   \\
 & \leqslant & 2\new{s}^3    \\ 
aD_2&=& -\new{s}^3 c\eta^3 m + \frac14 \new{s}^4 a\eta^4  + \frac12\new{s}^3 \po   a \eta^3  + \eta \pf + \frac14 \new{s}^4c\eta^2  \\
&\leqslant  & \frac{2}m a\eta \new{s}^3    \,.
\end{eqnarray*}
We have arrived to 
\[A^T M A \leqslant M + \begin{pmatrix}
 E_1 & 0 \\
0 &  aE_2
\end{pmatrix}
 \]
 with, using that $\new{\eta^2}/(1-\eta^2) +(1-\eta^2)/2 \geqslant 1/2 $ for all $\eta\in[0,1)$ and $5\new{s} \leqslant m$,
 \[E_1 = -\new{s}^2 m \po \frac{\new{\eta^2}}{1-\eta^2} + \frac{1-\eta^2}{2} - \frac{2}m \new{s}    \pf   \ \leqslant \ -\new{s}^2 \frac{ m}{10} \po \frac{2\new{\eta^2}}{1-\eta^2} +  1-\eta^2 \pf   \]
 and, using again that $5\new{s} \leqslant m$ and then \eqref{eq:condition},
 \begin{eqnarray*}
 E_2 &=&
  -  \po 1 - \eta^2 \pf  +     \new{s}^2 \co \frac{  \new{6} \eta^2}{m(1-\eta^2)}    + \new{\frac{2\eta }{5}} +   \frac{ 1 }{\new{m^2}} \cf 
\ \leqslant \ - \new{\frac45} (1-\eta^2)\,.
 \end{eqnarray*}
 We have obtained 
\[A^TM A \leqslant M +\max(E_1,E_2) \begin{pmatrix}
 1 & 0 \\ 0 & a
\end{pmatrix} \leqslant (1-4\rho) M\,,
\] 
 and then
  \[\|\Delta z_1\|_M  \leqslant \ \|A\Delta z_0\|_M  + \|\tilde z\|_M \leqslant \po 1 - 2\rho + \sqrt{2/a}C_{K,\delta}\pf \|\Delta z_0\|_M\,.\]
  The conclusion follows from the second part of \eqref{eq:condition}.   
  
  \medskip

\noindent \textbf{\new{Strong inertia case}.} \new{In this last part,  we follow computations similar to the previous case but focusing only on the terms of order $0$ and $1$ in $s$, namely we start by bounding, thanks to Lemma~\ref{lem:VerletMatrice+damping},
\begin{equation}\label{eq:ajout}
\left|\Delta z_1 - \begin{pmatrix}
1 & s \eta \\ -s \eta  \beta & \eta^2
\end{pmatrix} \Delta z_0\right| \leqslant \po \frac{s^2}{2} + \frac{5 s^3 e^{\modifplb{3s/2}}}{16}\pf |\Delta z_0|\,, 
\end{equation}
where $\beta\in\mathcal M(m,1)$. (Notice that, at this stage, it is clear that this bound is already too rough in the case $\eta=0$ to get a contraction afterwards; hence, from now on, assume that $\eta>0$). Since
\[ \begin{pmatrix}
1 & s \eta \\ -s \eta  \beta & \eta^2
\end{pmatrix} = I_{2d} + s\eta  \begin{pmatrix}
0 & 1 \\ - \beta & -\overline{\gamma}
\end{pmatrix} := I_{2d} + s\eta J_\beta\,,\]
 we recover the settings of \cite[Proposition 4]{MonmarcheContraction}. Indeed, for a positive symmetric matrix $M$ and $z\in\R^{2d}$, 
\[\left\| \po I_{2d} + s\eta J_\beta\pf   z\right\|_M^2 = \|z \|_M^2  + \modifplb{2} s\eta z \cdot M J_\beta z + (s\eta)^2 \| J_\beta z\|_M^2\,.\]
Focusing on the leading term in $s$, the goal is thus to find $M$ such that $z\cdot M J_\beta z \leqslant -\rho \|z\|_M^2$ for some $\rho>0$ for all $z\in\R^{2d}$ and all $\beta\in\mathcal M(m,1)$. It is proven in \cite[Proposition 4]{MonmarcheContraction} that such an $M$ exists only if $\overline{\gamma} \geqslant 1 - \sqrt{m} $,  in particular it is not possible to take $\overline{\gamma}$ of order $\sqrt{m}$ when $m$ is small, although this is the optimal choice in the Gaussian case. Hence, we recover in the present discrete-time settings the same limitation as observed in \cite{MonmarcheContraction} in the continuous-time case on the approach based on modified Euclidean distances contractions through parallel couplings. Besides, it is proven in \cite[Proposition 4]{MonmarcheContraction}  that, assuming that $\overline{\gamma}\geqslant 2$ and setting
\[M = \begin{pmatrix}
1 & \overline{\gamma}^{-1} \\ \overline{\gamma}^{-1} & 1 
\end{pmatrix} \]
then $z\cdot M J_\beta z \leqslant -m/(3\overline{\gamma}) \|z\|_M^2$ for all $z\in\R^{2d}$ and $\beta \in \mathcal M(m,1)$, and moreover $M \in \mathcal M(1/2,3/2)$, from which
\[\| J_\beta z\|_M^2 \leqslant |J_\beta|^2 |M| |M^{-1}| \|z\|_M^2 \leqslant 3(1+ \overline{\gamma}^2)\|z\|_M^2\,.\]  
As a conclusion, assuming that $m/(3\overline{\gamma}) \geqslant 3 s (1+\overline{\gamma}^2)$ (which implies that $s   \leqslant 1/90$, hence $15 s e^{\modifplb{3s/2}}/16 \leqslant 1/2$)
\begin{eqnarray*}
\|\Delta z_1 \|_M & \leqslant & \left\| \po I_{2d} + s\eta J_\beta\pf   \modifplb{\Delta z_0}\right\|_M + |M||M^{-1}| \po \frac{s^2}{2} + \frac{5 s^3 e^{\modifplb{3s/2}}}{16}\pf \|\Delta z_0\|_M \\
& \leqslant & \po \sqrt{1 - \frac{\modifplb{2}s\eta m}{3\overline{\gamma}} + 3s^2\eta^2(1+\overline{\gamma}^2)} + \frac{3 s^2}{2} + \frac{15 s^3 e^{\modifplb{3s/2}}}{16}\pf\|\Delta z_0\|_M  \\
& \leqslant & \po  1 - \frac{s\eta m}{\modifplb{3}\overline{\gamma}} + \frac{3}2 s^2\eta^2(1+\overline{\gamma}^2) + 2 s^2\pf\|\Delta z_0\|_M \\
& \leqslant & \po  1 - \frac{s\eta m}{\modifplb{6}\overline{\gamma}} \pf\|\Delta z_0\|_M \,,
\end{eqnarray*}
if $s \leqslant \eta m /[ \overline{\gamma}\modifplb{(21+9\overline{\gamma}^2)}]$, which is stronger than the previous condition.
}
  \end{proof}
  
\begin{proof}[Proof of Theorem~\ref{thm:crude}]

Thanks to the rescaling discussed in Remark~\ref{rem:rescaling},  
 Theorem~\ref{thm:crude} is proven by applying Proposition~\ref{prop:contraction}.  Indeed, let $\delta'=\sqrt{L} \delta$, $L'=1$ and $m'=m/L$. The conditions of Theorem~\ref{thm:crude} implies that $2\delta'   \leqslant 1 =  L'$.

 In the case $\eta=0$, the conditions \rev{ $K\delta'\leqslant 1/2$  and $\delta'  \leqslant m'K\delta' /31$}   ensures that Proposition~\ref{prop:contraction} applies to $Z'$, with a contraction rate \rev{$m'(\delta' K)^2/20 = m(\delta K)^2/20   $}, which concludes the proof in this case. 
 
 In the general case $\eta\in[0,1)$, \new{using the expression \eqref{eq:MrhoProp},} we bound
\begin{eqnarray*}
\rho'  & := &   \frac1{\new{5}(1-\eta^2)}\min\po  (1-\eta^2)^2,(K\delta')^2 \frac{ m'}{\new{8}} \po 2\new{\eta^2} +  (1-\eta^2)^2 \pf\pf\  \geqslant \  \frac{m' (K\delta')^2}{40(1-\eta^2)}\,,
\end{eqnarray*}
where we used that $ 2x +  (1-x^2)^2 \geqslant 1 $ for all $x\in[0,1]$ and that, from \eqref{eq:conditioncrude},
\[1-\eta^2   \geqslant  \frac{8K\delta'}{m'} \geqslant K\delta' \sqrt{\frac{m'}{10}} \,.\]
Second, using that $m'\leqslant 1$, we bound
\[
(K\delta')^2 \co \frac{   \new{6}\eta^2}{m'(1-\eta^2)}    + \new{\frac{2\eta }{5}} +  \frac{1  }{(m')^2} \cf   \ \leqslant \ \frac{\new{37} (K\delta')^2}{\new{5}(m')^2(1-\eta^2)}
\]
and thus the first part of \eqref{eq:condition} (applied with the parameters of $Z'$) is implied by
the first part of \eqref{eq:conditioncrude}. Third, using that $\delta' K \leqslant 2\sqrt{m'}/L'$, we bound
\begin{eqnarray*}
C_{K,\delta}' & := & \frac14 e^{\new{3(K-1)\delta' /2}}(\delta')^3  \po K^3 + K-1\pf\ \leqslant \  \frac12 e^{\new{3K\delta'/2}}( K\delta')^3  
\end{eqnarray*}
and thus (using the lower bound on $\rho'$) the second part of \eqref{eq:condition} (applied for $Z'$) is implied by
\[  e^{\new{3K\delta'/2}} K\delta'  \leqslant  \frac{(m')^{3/2} }{40(1-\eta^2)}\,.\]
As a conclusion, the conditions \eqref{eq:conditioncrude} satisfied by the parameters of $Z$ implies the conditions \eqref{eq:condition} for the parameters of $Z'$. Applying Proposition~\ref{prop:contraction} to $Z'$ concludes the proof of Theorem~\ref{thm:crude}. Notice that this gives a matrix $M$ which satisfies 
\[  \frac{2}{3}\po L |x|^2 + \frac{2m}{L}|v|^2\pf \leqslant  \|z\|_M^2 \leqslant \frac{4}{3}\po  L |x|^2 + \frac{2m}{L}|v|^2\pf \]
but then we can replace this $M$ by $M/L$ to get \eqref{eq:norme_equivalence2} in the conclusion since this doesn't impact the contraction property. \new{We also do this in the third case of Theorem~\ref{thm:crude}.}
\end{proof}

\subsection{Empirical averages}\label{sec:preuveEmpirical}

\begin{proof}[Proof of Proposition~\ref{prop:risque_quadra}]
If $f:\R^d \rightarrow \R$ is $1$-Lipschitz then, seen as a function on $\R^{2d}$ with $f(x,v)=f(x)$, using \eqref{eq:norme_equivalence}, we get $|f(x,v)-f(x',v')| \leqslant |x-x'| \leqslant \sqrt{3/2}\|(x,v)-(x',v')\|_M$, i.e. $f$ is $\sqrt{3/2}$ Lipschitz on $\R^{2d}$ with respect to $\|\cdot\|_M$. The characterization of the $\mathcal W_{M,1}$ distance in term of Lipschitz functions and Corollary~\ref{Cor:contraction} then imply that, for all $k\in\N$, $P_s^kf$ is  $\sqrt{3/2}(1-\rho)^k$-Lipschitz with respect to $\|\cdot\|_M$.

Let $\pi_0$ be a coupling of $\nu_0$ and $\mu_s$ and let $(Z_k^i,\tilde Z_k^i)_{k\in\N}$ for $i\in\cco 1,N\ccf$ be $N$ independent parallel couplings of two SG\new{g}HMC with an initial condition $(Z_0^i,\tilde Z_0^i)$ distributed according to $\pi_0$. Write $g = f -\pi_s f$. We bound
\begin{multline*}
\mathbb E\po \left|\frac 1{nN}\sum_{i=1}^N\sum_{k=n_0}^{n_0+n-1} g(Z_k^i)  \right|^2\pf\\    
\leqslant 2 \mathbb E\po \left|\frac 1{nN}\sum_{i=1}^N\sum_{k=n_0}^{n_0+n-1}\po g(Z_k^i) - g(\tilde Z_k^i)\pf\right|^2\pf   + \frac{2}N \mathbb E\po \left|\frac 1{n}\sum_{k=n_0}^{n_0+n-1} g(\tilde Z_k^1) \right|^2\pf \,,
\end{multline*}
where we used that the variables $1/n\sum_{k=n_0}^{n+n_0-1} g(\tilde Z_k^i)$ are centered and i.i.d. for $i\in\cco 1,N\ccf$.

Using the Jensen inequality, the deterministic contraction of Theorem~\ref{thm:crude} and the fact $g$ is $\sqrt{3/2}$-Lipschitz for $\|\cdot\|_M$,
\begin{align*}
\mathbb E\po \left|\frac 1{nN}\sum_{i=1}^N\sum_{k=n_0}^{n_0+n-1}\po g(Z_k^i) - g(\tilde Z_k^i)\pf\right|^2\pf  
& \leqslant \mathbb E\po \left|\frac 1{n}\sum_{k=n_0}^{n_0+n-1}\po g(Z_k^1) - g(\tilde Z_k^1)\pf\right|^2\pf\\   
& \leqslant \frac32    \po \frac 1{n} \sum_{k=0}^{n-1} (1-\rho  )^k \pf^2  (1-\rho  )^{2 n_0} \mathbb E \po \|Z_{0}^1-\tilde Z_{0}^1\|_M^2\pf \,.
\end{align*}
We can then take the infimum over $\pi_0$ to recover $\mathcal W_{M,2}^2(\nu_0,\mu_s)$.

Second, using that $\tilde Z_k^1 \sim \mu_s$ for all $k$ and the Markov property,
\begin{align*}
\mathbb E\po \left|\frac 1{n}\sum_{k=n_0}^{n_0+n-1} g(\tilde Z_k^1)  \right|^2\pf  
& \leqslant \frac2{n^2} \sum_{k=n_0}^{n_0+n-1} \sum_{j\leqslant k} \mathbb E \po g(\tilde Z_k^1)g(\tilde Z_j^1)\pf \\
& = \frac2{n^2} \sum_{k=n_0}^{n_0+n-1} \sum_{j\leqslant k}  \mu_s \po g P_s^{k-j} g\pf \\
& \leqslant  \frac2{n^2} \sum_{k=n_0}^{n_0+n-1} \sum_{j\leqslant k}  \sqrt{\mu_s\po g^2 \pf \mu_s\po (P_s^{k-j}g)^2\pf }\,. 
\end{align*}
Now, for $h$ a $r$-Lipschitz function for $\|\cdot\|_M$ with $\mu_s h =0$,
 \begin{eqnarray*}
\mu_s(h^2) &  = &  \frac12 \int_{\R^{2d}\times\R^{2d}}|h(z)-h(z')|^2 \mu_s(\dd z)\mu_s(\dd z') \ \leqslant \   r^2 \mathrm{Var}_M(\mu_s) \,,
 \end{eqnarray*}
 from which
\[\mathbb E\po \left|\frac 1n\sum_{k=n_0}^{n_0+n-1} g(\tilde Z_k)  \right|^2\pf  \leqslant \frac{3  \mathrm{Var}_M(\mu_s)}{n^2} \sum_{k=n_0}^{n_0+n-1} \sum_{j\leqslant k} (1-  \rho)^{k-j} \leqslant \frac{3  \mathrm{Var}_M(\mu_s)}{n \rho }\,. \]

It remains to bound the variance. Using the scaling properties of the SG\new{g}HMC as in the proof of Theorem~\ref{thm:crude}, for the remaining of the proof we assume without loss of generality that $L=1$ (and, at the end of the computation, we will have to replace $m$ and $\delta$ by $m/L$ and $\delta\sqrt{L}$). For  $n\in\N$, a fixed $x_*\in\R^{d}$ and $z_*=(x_*,0)$, using the variational formula of the variance, we bound
\begin{eqnarray*}
\sqrt{\mathrm{Var}_M(\mu_s)}  & \leqslant & \sqrt{\int_{\R^{2d}}\left\|z-z_*\right\|_M^2 \mu_s(\dd z)}\\
& =& \mathcal W_{M,2}(\mu_s,\delta_{z_*})\\
& \leqslant &   \mathcal W_{M,2}(\mu_s,\delta_{z_*}P_s^n) +  \mathcal W_{M,2}(\delta_{z_*}P_s^n,\delta_{z_*})\\
& \leqslant &   (1-\rho)^n\mathcal W_{M,2}(\mu_s,\delta_{z_*}) +  \mathcal W_{M,2}(\delta_{z_*}P_s^n,\delta_{z_*})
\end{eqnarray*}
where we used Corollary~\ref{Cor:contraction} and the fact $\mu_s = \mu_sP_s^n$. Chosing $n=n_*=\lceil -\ln 2/\ln(1-\rho)\rceil$,
\begin{equation}\label{eq:Varmus}
\mathrm{Var}_M(\mu_s) \leqslant \frac{\mathcal W_{M,2}^2(\delta_{z_*}P_s^{n_*},\delta_{z_*})}{\po 1 - (1-\rho)^{n_*} \pf^2 } \leqslant    \frac{16}{3} P_s^{n_*} \varphi(z_*)\,,
\end{equation}
with $\varphi(z)=|x-x_*|^2 + 2m|v|^2$, where we used \eqref{eq:norme_equivalence}. 

 For any $z=(x,v)\in\R^{2d}$,
\[\new{P_R}\varphi(z) = |x-x_*|^2 + 2m \co \new{\eta^2} |v|^2 + (1-\eta^2) \mathbb E (|G|^2)\cf \leqslant \varphi(z)+ 2m(1-\eta^2) d\,,\]
and
\[P_{sV} \varphi(z)= \mathbb E_{\omega} \po \left|x-x_* + \delta v - \frac{\delta^2}{2} b\po x +\frac\delta2 v,\theta \pf\right|^2 + 2m \left|v-\delta  b\po x +\frac\delta2 v,\theta \pf\right|^2 \pf .\]
Writing $v' = v-\delta b(x+\delta v/2,\theta)$, using Assumption~\ref{Assu:mLb}, we bound
\begin{eqnarray*}
|v'|^2 &\leqslant & \po |v|+ \delta | b(x_*,\theta)| + \delta  \po  |x-x_*| + \delta/2 |v|\pf\pf ^2 \\
& \leqslant & (1+5\delta) |v|^2 + 4\delta | b(x_*,\theta)|^2 + 4\delta  |x-x_*|^2
\end{eqnarray*}
where we used that   $\delta\leqslant 1$ to simplify the higher order terms in $\delta$, and thus (using also that, from \eqref{eq:conditioncrude}, $\delta\leqslant 2\sqrt{m}$  and $m\leqslant 1$ to simplify),
\begin{eqnarray*}
P_{sV} \varphi(z)  & \leqslant & \mathbb E_{\omega} \po |x-x_*|^2  +  \delta |x-x_*||v+v'| + \frac{\delta^2}4 |v+v'|^2 + 2m |v'|^2\pf \\
&\leqslant &  \mathbb E_\omega \po \po 1+ \frac{\delta}{2\sqrt{m}} \pf |x-x_*|^2 + \po \frac{\delta}{\sqrt{m}} + \frac{\delta^2}{2m}\pf m|v|^2 + \po 2 + \frac{\delta}{\sqrt{m}}+ \frac{\delta^2}{2m} \pf m |v'|^2 \pf\\
&\leqslant & \po 1+  \delta   \frac{25}{\sqrt{m}}  \pf |x-x_*|^2 + \po 2+  \delta   \frac{24}{\sqrt{m}}  \pf m  |v|^2 + 24m \delta \mathbb E_{\omega} \po |b(x_*,\theta)|^2\pf \\
& \leqslant & \po 1+  \delta   \frac{25}{\sqrt{m}}  \pf \varphi(z) + 24 m \delta \mathbb E_{\omega} \po |b(x_*,\theta)|^2\pf \,.
\end{eqnarray*}
By induction,
\[P_{sV}^K \varphi(z) \leqslant  e^{25 \delta K /\sqrt{m} }\po \varphi(z) + 24 m \delta K  \mathbb E_{\omega} \po |b(x_*,\theta)|^2\pf \pf\,. \]
Gathering the previous bounds,
\[P_s\varphi(z) = \new{P_R} P_{sV}^{\new{K}} \new{P_R} \varphi(z) \leqslant  e^{25 \delta K /\sqrt{m}  } \po  \varphi(z) +  4 (1-\eta^2) d + 24 m  \delta K  \mathbb E_{\omega} \po |b(x_*,\theta)|^2\pf\pf  \,, \] 
and thus
\[P_s^{n_*}\varphi(z_*) \leqslant  4 n_* e^{25 n_*\delta K /\sqrt{m}  } \po    (1-\eta^2) d + 6 m  \delta K  \mathbb E_{\omega} \po |b(x_*,\theta)|^2\pf\pf \,, \] 
where we used that $\varphi(z_*)=0$. Plugging this estimate in \eqref{eq:Varmus}, taking the infimum over $x_*\in\R^d$ and using that $n_* \leqslant 1+ \ln(2)/\rho \leqslant 1/\rho$ (indeed, the first condition  in \eqref{eq:conditioncrude} implies that $\rho \leqslant 1/(7^2\times 40)$) 
 concludes the proof.
 
\end{proof}

\subsection{\new{Equilibrium bias}}\label{section:erreurNum}

\new{
This section is devoted to the proofs of Propositions~\ref{prop:numerique} and \ref{prop:erreur_sto}. In each case, it is based on a finite-time error bound which is then combined with the long-time contraction established in Theorem~\ref{thm:crude} to conclude.

\begin{proof}[Proof of Proposition~\ref{prop:numerique}]
Again, following the rescaling of Remark~\ref{rem:rescaling}, we only consider the case $L=1$. Recall that $\Phi_\delta$ denotes the the Verlet step \eqref{eq:VerletSto} in the non-stochastic case $b(x,\theta)=\na U(x)$. Let $(\varphi_t)_{t\geqslant 0}$ be the Hamiltonian flow associated to $U$, namely $(x_t,v_t):=\varphi_t(x,v)$ solves 
\[\dot x_t = v_t\,,\quad \dot v_t = -\na U(x_t),\quad (x_0,v_0)=(x,v)\,.\]
Using that $|\na U(x)|\leqslant |x|$ since $\na U(0)=0$, the standard numerical analysis of the Verlet integrator (see e.g. \cite[Lemma 29]{MonmarcheSplitting})  yields, for all $z=(x,v)\in\R^{2d}$,
\begin{equation}\label{eq:Verlet_num1}
|\Phi_\delta(z) - \varphi_\delta(z)|  \leqslant  2 \delta^2 |z|   
\end{equation}
and, if we assume additionnally that $\na^2 U$ is Lipschitz,
\begin{equation}\label{eq:Verlet_num2}
|\Phi_\delta(z) - \varphi_\delta(z)|  \leqslant     2 \delta^3 (1+L_2)(1+  |z|^2)\,.
\end{equation}
Considering the Markov transition $Q_\delta$ given by $Q_\delta f(z) = f(\varphi_{\delta}(z))$ (which leaves invariant $\mu$), we get that
\[\mathcal W_2 (\mu P_V,\mu) = \mathcal W_2 (\mu P_V,\mu Q_\delta) \leqslant \sqrt{\mathbb E_\mu\po |\Phi_\delta(Z) - \varphi_\delta(Z)|^2\pf }\,.\]
Depending whether we use \eqref{eq:Verlet_num1} or \eqref{eq:Verlet_num2}, we bound the right hand side in terms of  $m_j = \int_{\R^{2d}} |z|^j \mu(\dd z)$ with either $j=2$ or $j=4$. From \cite[Lemma 30]{MonmarcheSplitting}, 
\begin{equation}\label{eq:m2}
m_2 \leqslant \po 1+ \frac1m\pf d\,,\qquad  m_4 \leqslant \po 1+\frac1{m^2}\pf(d+2)^2 \,,
\end{equation}
 from which
 \begin{equation}\label{eq:W2erreur1}
\mathcal W_2 (\mu P_V,\mu) \leqslant   3 \delta^2 \sqrt{\frac{d}{m}}     
 \end{equation}
if we use \eqref{eq:Verlet_num1} and 
 \begin{equation}\label{eq:W2erreur2}
\mathcal W_2 (\mu P_V,\mu) \leqslant   11 \delta^3 (1+L_2) \frac{d}{m} 
 \end{equation}   
if we use \eqref{eq:Verlet_num2} (we used that $d\geqslant 1 \geqslant m$).
Moreover, using that $\Phi_V$ is $(1+2\delta)$-Lipschitz (since $\delta\leqslant 1$), we get
\begin{equation}\label{eq:Wp_controlVerlet}
\forall \nu,\nu'\in\mathcal P(\R^{2d})\,,\qquad \mathcal W_2(\nu P_V,\nu' P_V) \leqslant (1+2\delta) \mathcal W_2(\nu,\nu')\,.
\end{equation}
   Besides, the \new{randomization} step $\new{P_R}$ leaves invariant $\mu$ and is such that $\mathcal W_2 (\nu \new{P_R},\nu' \new{P_R}) \leqslant \mathcal W_2 (\nu,\nu')$ for all $\nu,\nu'$  (with a parallel coupling). Then, by straightforward inductions on $K$ and then $n$,
  \begin{eqnarray*}
  \mathcal W_2 \po \mu P^n,\mu\pf & = & \mathcal W_2\po \mu  P^{n-1} \new{P_R}P_V^K\new{P_R},\mu  \new{P_R}\pf \\
  & \leqslant  & \mathcal W_2\po \mu P^{n-1} \new{P_R}P_V^K,\mu \pf\\
    & \leqslant  & \mathcal W_2\po \mu  P^{n-1} \new{P_R}P_V^K,\mu  P_V\pf +  \mathcal W_2\po \mu   P_V,\mu \pf\\ 
    & \leqslant& (1+2\delta) \mathcal W_2\po \mu  P^{n-1} \new{P_R} P_V^{K-1},\mu  \pf + \mathcal W_2\po \mu   P_V,\mu \pf \\
    & \leqslant& (1+2\delta) \mathcal W_2\po \mu  P^{n-1} \new{P_R}  ,\mu  \pf + \mathcal W_2\po \mu   P_V,\mu \pf \sum_{j=0}^{K-1} (1+2\delta)^j \\
    & \leqslant& (1+2\delta) \mathcal W_2\po \mu  P^{n-1}    ,\mu  \pf + \mathcal W_2\po \mu   P_V,\mu \pf \sum_{j=0}^{K-1} (1+2\delta)^j \\
    & \leqslant &  \mathcal W_2\po \mu   P_V,\mu \pf \sum_{j=0}^{nK-1} (1+2\delta)^j \,.
  \end{eqnarray*}
  Now, using the contraction \eqref{eq:contractMz} yields
  \[\mathcal W_2 \po \tilde\mu,\mu P^n\pf =  \mathcal W_2 \po \tilde\mu P^n,\mu P^n\pf \leqslant |M||M^{-1}|  (1-\rho)^n   \mathcal W_2 \po \tilde\mu  ,\mu \pf\,.  \]
  Combining these two bounds applied with $n = \lceil -\ln(2|M||M^{-1}|)/\ln(1-\rho)\rceil$, we get
  \[\mathcal W_2(\mu,\tilde\mu) \leqslant \mathcal W_2(\mu P^n,\tilde\mu) + \mathcal W_2(\mu,\mu P^n) \leqslant \frac12 \mathcal W_2(\mu,\tilde\mu) + \frac{ (1+2\delta) ^{nK}}{2\delta}\mathcal W_2\po \mu   P_V,\mu \pf  \,,\]
  and thus, using that $\ln(1-\rho)\leqslant -\rho$, $\rho\leqslant 1$ and $1+\ln(2) \leqslant \ln(6)$,
  \[\mathcal W_2(\mu,\tilde\mu) \leqslant   \frac{ e ^{2nK\delta }}{\delta}\mathcal W_2\po \mu   P_V,\mu \pf  \leqslant   \frac{ e ^{2 \ln(6|M||M^{-1}|) nK\delta/\rho }}{\delta}\mathcal W_2\po \mu   P_V,\mu \pf \,.\]
  Conclusion follows either from \eqref{eq:W2erreur1} or \eqref{eq:W2erreur2}.
 \end{proof}

From this, we get the following bound on the second moment of $\tilde \mu$:

\begin{lemma}\label{lem:moment_mu_tilde}
Under Assumption~\ref{Assu:mLb} with $L=1$ $b_\theta = \na U$ for some $U\in\mathcal C^2(\R^d)$, assume that $M,\rho>0$ are such that \eqref{eq:contractMz} holds for the parallel coupling of the \new{g}HMC chain. Then
\[\int_{\R^{2d}} |z|^2 \tilde \mu(\dd z) \leqslant C_2 d \qquad \text{with}\qquad C_2= \frac{ 1}{m} \po 4 + 18 \delta^2  (6|M||M^{-1}|)^{4 K\delta/\rho} \pf \,. \] 
\end{lemma}

\begin{proof}
Given $(Z,\tilde Z)$ a coupling of $\mu$ and $\tilde \mu$,
\[  \mathbb E\po |\tilde Z|^2 \pf \leqslant 2 \mathbb E\po | Z|^2 \pf  +\modifplb{2} \mathbb E\po |\tilde Z-Z|^2 \pf \,. \]  
Taking the infimum over all such couplings,
conclusion follows from \eqref{eq:m2} and \eqref{eq:erreur_num_prop1} (using that $m\leqslant 1$).
\end{proof}

Before analysing the contribution of the stochastic gradient error, we need the following result on its variance.

\begin{lemma}\label{lem:gradsto}
Under Assumption~\ref{assu:gradient_sto}, assuming furthermore that $L=1$ and $\na U(0)=0$, then, for all $x\in\R^d$,
\[\mathbb E\po |b_{\theta}(x) - \na U(x)|^2 \pf \leqslant \frac1p \po \modifplb{6}|x|^2 + 3 \mathbb E\po |\na U_{\theta_1}(0)|^2 \pf\pf  <\infty\,.\]
\end{lemma}
\begin{proof}
By independence and then, using that $b_{\theta_1}$  and $\na U$ are Lipschitz, for any $y\in\R^d$,
\begin{eqnarray*}
\mathbb E\po |b_{\theta}(x) - \na U(x)|^2 \pf & = & \frac1p \mathbb E\po |\na U_{\theta_1}(x) - \na U(x)|^2 \pf \\
& \leqslant &   \frac1p \mathbb E\po \po 2|x-y|+ |\na U_{\theta_1}(y) - \na U(y)|\pf ^2 \pf \\
& \leqslant &  \frac{\modifplb{6}}{p}|x-y|^2  \modifplb{+}\frac3p \mathbb E\po |\na U_{\theta_1}(y) - \na U(y)|^2 \pf \,.
\end{eqnarray*}
 In particular, we get that this variance is in fact finite for all $x\in\R^d$. Taking the previous inequality at $y=0$ concludes.
\end{proof}

 \begin{proof}[Proof of Proposition~\ref{prop:erreur_sto}]
Fix $z_0=(x_0,v_0) \in\R^{2d}$ and assume by rescaling that $L=1$. Given $(\theta_k,G_k,G_k')_{k\in\N}$ an i.i.d. sequence of random variables where $\theta_k$, $G_k$ and $G_k'$ with $\theta_k \sim \omega$ and $G_k,G_k'\sim\mathrm{N}(0,I_d)$,  define $(z_k)_{k\in\N}=(x_k,v_k)_{k\in\N}$ by induction as follows:
\begin{eqnarray*}
v_{k+1/3} & = & \eta_k v_k + \sqrt{1-\eta_k^2} G_k \\
x_{k+1/2} &= & x_k + \frac{\delta}{2}v_{k+1/3}\\
v_{k+2/3} &= & v_{k+1/3} - \delta b_{\theta_k}(x_{k+1/2}) \\
x_{k+1} &= & x_{k+1/2} + \frac{\delta}{2}v_{k+2/3}\\
v_{k+1} & = & \eta_k v_{k+2/3} + \sqrt{1-\eta_k^2} G_k' \,,
\end{eqnarray*}
where $\eta_k = \eta$ if $k=0$ modulo $K$ and $\eta_k=1$ otherwise. This ensures that $(z_{Kn})_{n\in\N}$ is a SG\new{g}HMC chain starting at $z_0$. Define $(z_k')_{k\in\N}$ exactly as $(z_k)_{k\in\N}$  (with the same variables $\theta_k,G_k,G_k'$) except that $b_{\theta_k}$ is replaced by $\na U$, so that $(z_k')_{k\in\N}$ is a (non-stochastic) gHMC chain starting at $z_0$. Denote $\Delta z_{k+i} = z_{k+i} - z_{k+i'}$ (with $i\in\{\emptyset,1/3,1/2,2/3\}$) and similarly for $x,v$. Using that $|a+sb|^2 \leqslant (1+s)|a|^2 + (s+s^2)|b|^2$ for $a,b\in\R, s>0$, 
\begin{eqnarray*}
\mathbb E \po |\Delta v_{k+1/3}|^2 \pf  & \leqslant  & \mathbb E \po |\Delta v_{k}|^2 \pf\\
\mathbb E \po |\Delta x_{k+1/2}|^2 \pf  &\leqslant  & \po 1+ \frac{\delta}{2}\pf \mathbb E \po |\Delta x_{k}|^2 \pf  + \po \frac{\delta}{2} + \frac{\delta^2}{4}\pf \mathbb E \po |\Delta v_{k+1/3}|^2 \pf\,,
\end{eqnarray*}
and we use similar bounds for $\Delta x_{k+1}$ and $\Delta v_{k+1}$. Using that $b_{\theta_k}(x_{k+1/2}) - \na U(x_{k+1/2})$ averages as $0$ conditionally to $\Delta v_{k+1/3},x_{k+1/2},x'_{k+1/2}$ and Lemma~\ref{lem:gradsto},
\begin{eqnarray*}
\lefteqn{\mathbb E \po |\Delta v_{k+2/3}|^2 \pf}\\
 & =  & \mathbb E \po |\Delta v_{k+1/3} + \delta(b_{\theta_k}(x'_{k+1/2}) - b_{\theta_k}(x_{k+1/2})|^2 \pf + \delta^2  \mathbb E \po |b_{\theta_k}(x_{k+1/2}') - \na U(x_{k+1/2}')|^2 \pf \\
 & \leqslant & \po 1+ \delta \pf \mathbb E \po |\Delta v_{k+1/3}|^2 \pf +(\delta +\delta^2 ) \mathbb E \po |\Delta x_{k+1/2}|^2 \pf +    \frac{\delta^2}p  \mathbb E\po   \modifplb{6}|x_{k+1/2}'|^2 + 3 |\na U_{\theta_1}(0)|^2 \pf \,.
\end{eqnarray*}

Combining these bounds leads to
\begin{eqnarray}
\mathbb E \po |\Delta z_{k+1}|^2 \pf  &\leqslant&  \po 1+ 2\delta +  \frac{9\delta^2}{4} + \frac{3\delta^3}{2}  + \frac{5\delta^4}{8} + \frac{\delta^5}{8}\pf  \mathbb E \po |\Delta x_{k}|^2 \pf  \nonumber \\
& & + \po 1+ 2 \delta  + \frac{7\delta^2}{4} + \frac{11\delta^3}{8}+  \frac{3\delta^4}{4}  + \frac{5\delta^5}{16} + \frac{\delta^6}{16} \pf\mathbb E \po |\Delta v_{k}|^2 \pf \nonumber\\
& & + \po 1+ \frac{\delta}{2} + \frac{\delta^2}{4}\pf \frac{\delta^2}p \mathbb E\po   \modifplb{6}|x_{k+1/2}'|^2 + 3 |\na U_{\theta_1}(0)|^2 \pf \nonumber \\
& \leqslant & (1+4\delta) \mathbb E \po |\Delta z_{k}|^2 \pf + 2 \frac{\delta^2}p \mathbb E\po   \modifplb{6}|x_{k+1/2}'|^2 + 3 |\na U_{\theta_1}(0)|^2 \pf \label{eq:combining}
\end{eqnarray}
 using that $\delta\leqslant 1/2$. Besides,
 \[\mathbb E\po   |x_{k+1/2}'|^2 \pf \leqslant \po 1+ \frac{\delta}{2}\pf  \mathbb E\po |z_k'|^2 \pf +  \frac{\delta^2}{4}(1-\eta_k^2) d \,. \]

  Recall that $\tilde \mu$ is the invariant measure of the non-stochastic unadjusted \new{g}HMC chain $P$. In particular, assuming that $z_0 \sim \tilde \mu$, we get that $z'_{Kn}\sim \tilde \mu$ for all $n\in\N$.  Since we assumed that $0$ is the minimizer of $U$, we get that the origin is fixed by the Verlet steps, i.e. $\delta_0 P_V = \delta_0$. Hence, for $j\in\cco 0,K-1\ccf$, using \eqref{eq:Wp_controlVerlet},
 \begin{eqnarray*}
  \mathbb E \po |z_{Kn+j}'|^2 \pf  &=&  \mathcal W_2^2(\tilde \mu \new{P_R} P_V^j,\delta_0 P_V^j) \\
 & \leqslant &   (1+2\delta)^{2 K} \mathcal W_2^2(\tilde \mu \new{P_R}  , \delta_0) \\ 
  & \leqslant &   2(1+2\delta)^{2 K} \po \mathcal W_2^2(\tilde \mu \new{P_R}    , \delta_0 \new{P_R}) + \mathcal W_2^2(\delta_0 \new{P_R} , \delta_0)\pf  \\
  & \leqslant &   2e^{4K\delta} \po \mathcal W_2^2(\tilde \mu     , \delta_0 ) + (1-\eta^2)d \pf  \,,
 \end{eqnarray*}
where we used that $\new{P_R}$ decreases the Wasserstein distance.  Using Lemma~\ref{lem:moment_mu_tilde} to bound $\mathcal W_2(\tilde \mu     , \delta_0 ) $, we gather the last inequalities to bound $\mathbb E(|x_{k+1/2}'|^2 )$ in \eqref{eq:combining}, which then reads
\[\mathbb E \po |\Delta z_{k+1}|^2 \pf \leqslant  (1+4\delta) \mathbb E \po |\Delta z_{k}|^2 \pf + \frac{R \delta^2}p \]
with
\begin{eqnarray*}
R &=&      \modifplb{24}  \po 1+ \frac{\delta}{2}\pf  e^{4K\delta} \po C_2 d + (1-\eta^2)d \pf  + \modifplb{12} \frac{\delta^2}{4}  d + 6 \mathbb E\po  |\na U_{\theta_1}(0)|^2 \pf \\
& \leqslant &     \modifplb{30}   e^{4K\delta} \po C_2  +  2  \pf  \po d   +  \mathbb E\po  |\na U_{\theta_1}(0)|^2 \pf\pf  \,.
\end{eqnarray*}
Hence, by induction, using that $(z_{Kn},z_{Kn}')$ is a coupling of $\tilde \mu P_s^n$ and $\tilde \mu P^n = \tilde \mu $,
\[\mathcal W_2^2 \po \tilde \mu P_s^n,\tilde \mu \pf \leqslant \mathbb E \po |\Delta z_{Kn}|^2 \pf  \leqslant  \frac{(1+4\delta)^{Kn}R \delta}{4p} \,.\] 

The conclusion is then the same as in the proof of Proposition~\ref{prop:numerique}.  Using the contraction \eqref{eq:contractMz} yields
  \[\mathcal W_2 \po \tilde\mu_s,\tilde \mu P_s^n\pf =  \mathcal W_2 \po \tilde\mu_s P_s^n,\tilde \mu P_s^n\pf \leqslant |M||M^{-1}|  (1-\rho)^n   \mathcal W_2 \po \tilde\mu_s  ,\tilde \mu \pf\,.  \]
  Taking $n = \lceil -\ln(2|M||M^{-1}|)/\ln(1-\rho)\rceil$, we get (using that $\delta\leqslant 1/2$)
  \[\mathcal W_2(\tilde \mu_s,\tilde\mu) \leqslant \mathcal W_2(\tilde\mu_s,\tilde \mu P_s^n) + \mathcal W_2(\tilde \mu,\tilde \mu P_s^n) \leqslant \frac12 \mathcal W_2(\mu,\tilde\mu) +  \frac{e^{2K\delta n}\sqrt{R \delta}}{2\sqrt{p}}  \,,\]
  and thus, using that $\ln(1-\rho)\leqslant -\rho$, $\rho\leqslant 1$ and $1+\ln(2) \leqslant \ln(6)$,
  \[\mathcal W_2(\tilde\mu_s,\tilde\mu) \leqslant  \frac{e^{2K\delta n}\sqrt{R \delta}}{\sqrt{p}} 
  \leqslant   \frac{ e ^{2 \ln(6|M||M^{-1}|) K\delta/\rho  }\sqrt{R \delta}}{\sqrt{p}}  \,,\]
  which concludes the proof of Proposition~\ref{prop:erreur_sto}.
 \end{proof}

 }

\section{\rev{Additional numerical experiments}}\label{sec:additional_numerics}

\rev{The settings of these experiments are exactly the same as in Section~\ref{sec:numerique}, except that $d=50$. The target measure is the Gaussian distribution with diagonal covariance matrix with coefficients $1+2k$ for $k \in \cco 0,49\ccf$, which gives (approximately) the same condition number as the experiments of Section~\ref{sec:numerique} (this removes half of the resonances). As we see, the results are exactly the same and leads to the same comments.
 }

\begin{figure}
\begin{subfigure}{0.49\textwidth}
    \centering
    \includegraphics[scale=0.38]{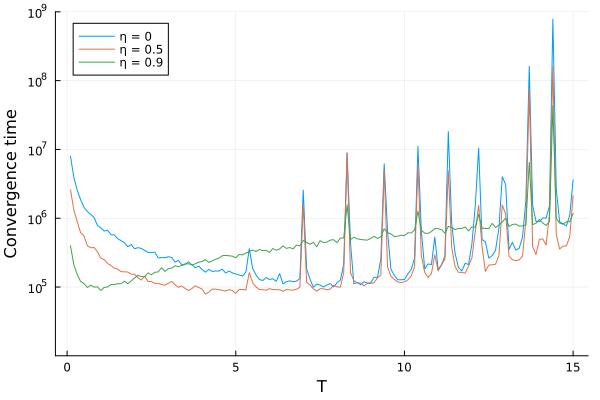}
    \caption{ugHMC ($T=K\delta$, $K$ varies)}
    \label{fig:unadjusted_T50}
    \end{subfigure}
\begin{subfigure}{0.49\textwidth}
    \centering
    \includegraphics[scale=0.38]{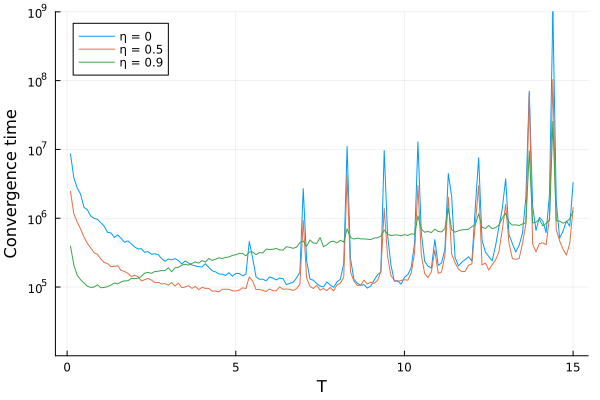}
    \caption{MagHMC ($T=K\delta$, $K$ varies)}
    \label{fig:adjusted_K50}
\end{subfigure}
\begin{subfigure}{0.49\textwidth}
    \centering
    \includegraphics[scale=0.38]{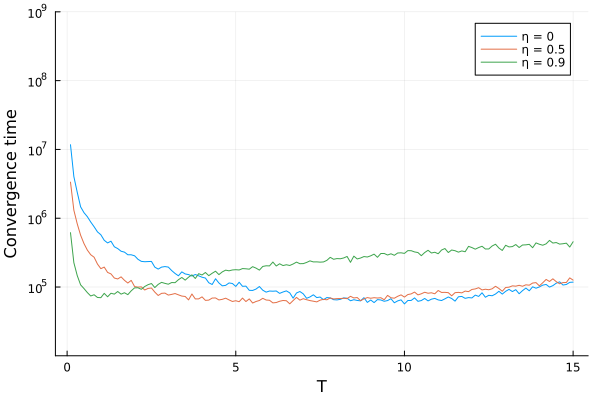}
    \caption{PDMP}
    \label{fig:PDMP-T50}
  \end{subfigure} 
  \begin{subfigure}{0.49\textwidth}
    \centering
    \includegraphics[scale=0.38]{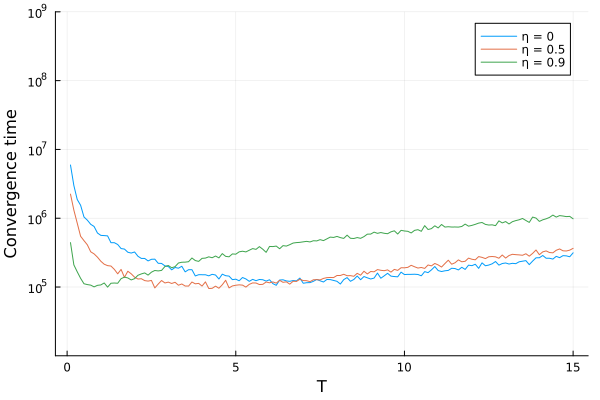}
    \caption{randomized ugHMC ($T=K\mathbb E(\tilde \delta)$, $K$ varies)}
    \label{fig:unadjusted-randomized50.}
    \end{subfigure}
    \caption{Influence of the integration time $T$ (same as Figure~\ref{figure:T} but with $d=50$).} 
\end{figure}

\begin{figure}
\begin{subfigure}{0.49\textwidth}
    \centering
    \includegraphics[scale=0.38]{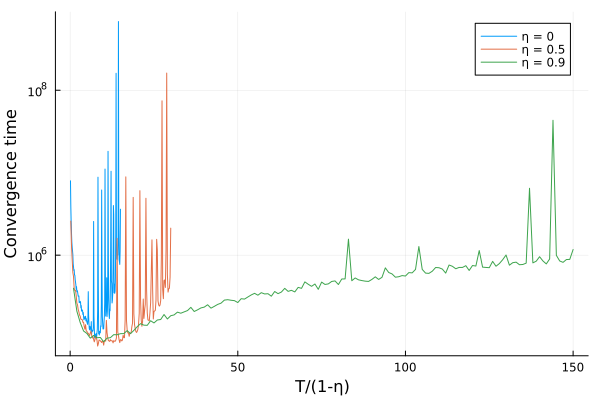}
    \caption{ugHMC}
    \label{fig:uT1-eta50}
     \end{subfigure}
     \begin{subfigure}{0.49\textwidth}
    \centering
    \includegraphics[scale=0.38]{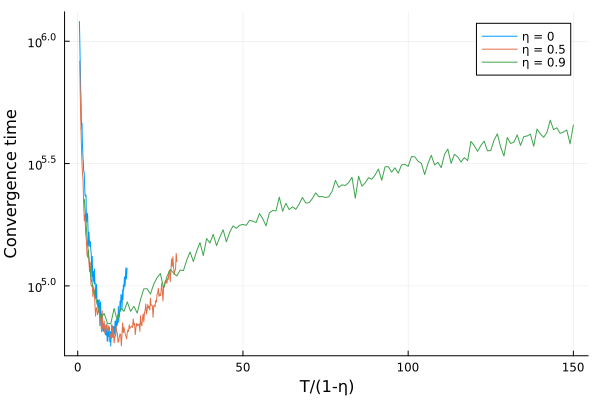}
    \caption{PDMP}
    \label{fig:PDMPT1-eta50}
  \end{subfigure}
  \caption{Influence of $T/(1-\eta)$ (same as Figure~\ref{figure:Tsur1-eta} but with $d=50$).}
\end{figure}

\begin{figure}
\begin{subfigure}{0.49\textwidth}
    \centering
    \includegraphics[scale=0.38]{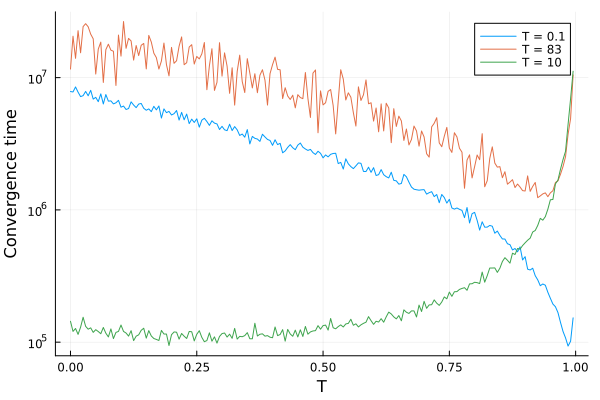}
    \caption{ugHMC}
    \end{subfigure}
\begin{subfigure}{0.49\textwidth}
    \centering
    \includegraphics[scale=0.38]{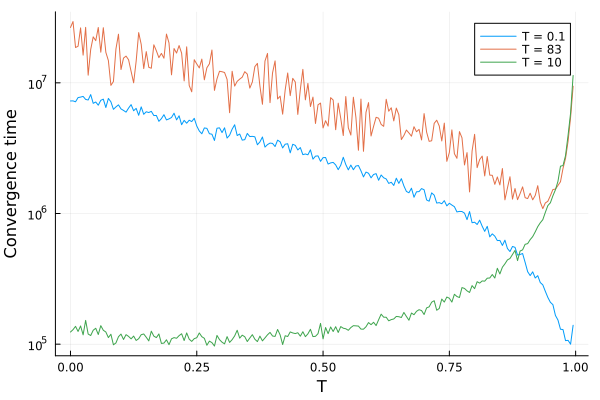}
    \caption{MagHMC}
\end{subfigure}
\centering
\begin{subfigure}{0.49\textwidth}
    \centering
    \includegraphics[scale=0.38]{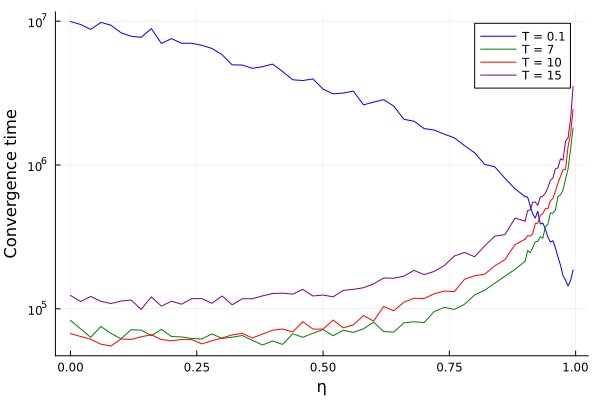}
    \caption{PDMP}
  \end{subfigure}  
  \caption{Influence of $\eta$ (same as Figure~\ref{figure:eta} but with $d=50$). For ugHMC and MagHMC, $T=83$ leads to a resonance in this setting.}
\end{figure}


\end{document}